\documentclass{amsart}
\usepackage{amssymb,amsmath, amsthm,latexsym}
\usepackage{graphics}
\usepackage{amscd}
\usepackage{graphics}
\usepackage{here}
\newcommand{\cal}[1]{\mathcal{#1}}
\theoremstyle{plain}
\newtheorem{coro}{Corollary}
\newtheorem{theorem}{Theorem}
\newtheorem{lemma}{Lemma}[section]
\newtheorem{theo}[lemma]{Theorem}
\newtheorem{proposition}[lemma]{Proposition}
\newtheorem{corollary}[lemma]{Corollary}
\theoremstyle{definition}
\newtheorem{definition}[lemma]{Definition}
\newtheorem{definition*}{Definition A}
\theoremstyle{plain}
\newtheorem{lemma*}[definition*]{Lemma A}
\newtheorem{proposition*}[definition*]{Proposition A}
\parskip=\bigskipamount

\let\egthree=\phi
\let\phi=\varphi
\let\varphi=\egthree

\begin{document}
\title{Geometry of the mapping class groups I:
Boundary amenability}
\author{Ursula Hamenst\"adt}
\thanks
{
Research
partially supported by Sonderforschungsbereich 611}
\date{March 19, 2008}

\begin{abstract}
We construct a geometric model for the mapping class group
${\cal M\cal C\cal G}$ of a non-exceptional oriented surface $S$ of
genus $g$ with $k$ punctures and use it to show that the action of
${\cal M\cal C\cal G}$ on the compact metrizable Hausdorff
space of complete geodesic
laminations for $S$ is topologically amenable. As a consequence,
the Novikov higher signature conjecture holds for every
subgroup of ${\cal M\cal C\cal G}$. 
\end{abstract}

\maketitle

\tableofcontents

\section{Introduction}

A countable group $\Gamma$ is called \emph{boundary amenable} if
it admits a \emph{topologically amenable} action on a compact Hausdorff
space $X$. This means that $\Gamma$ acts on $X$ as a group of
homeomorphisms and that moreover the following holds.
Let ${\cal P}(\Gamma)$ be the convex space of
all probability measures on $\Gamma$. Note that
${\cal P}(\Gamma)$ can be viewed as a subset of the
unit ball in the space $\ell^1(\Gamma)$ of summable functions
on $\Gamma$ and therefore it admits a natural norm $\Vert\,\Vert$.
The group $\Gamma$ acts on
$({\cal P}(\Gamma),\Vert \, \Vert)$ isometrically by left translation. 
We require that
there is a sequence of weak$^*$-continuous
maps $\xi_n:X\to {\cal P}(\Gamma)$ with the property that $\Vert
g\xi_n(x)-\xi_n(gx)\Vert\to 0$ $(n\to \infty)$
uniformly on compact subsets of
$\Gamma\times X$ (see \cite{AR00} for more on
amenable actions).
By a result of Higson \cite{H00} which is based on earlier work 
of Yu \cite{Y98,Y00},
for any countable group $\Gamma$ which is boundary amenable
and for every separable $\Gamma-C^*$-algebra $A$, the Baum-Connes
assembly map
\[\mu:KK_*^\Gamma({\cal E}\Gamma,A)\to
KK(\mathbb{C},C_r^*(\Gamma,A))\] is split injective. As a
consequence, the strong Novikov conjecture holds for
$\Gamma$ and hence the Novikov higher signature conjecture
holds as well \cite{BCH94,MV03}.

Now let $S$ be a non-exceptional 
oriented surface of finite type. This means that
$S$ is a
closed surface of genus $g\geq 0$ from which $k\geq 0$ points,
so-called \emph{punctures}, have been deleted, and where
$3g-3+k\geq 2$. 
The \emph{mapping class group}
${\cal M\cal C\cal G}$ of $S$ is the group 
of all isotopy classes of orientation preserving
self-homeomorphisms of $S$. The mapping class group is
finitely presented. We
refer to \cite{I02} for a summary of
the basic
properties of the mapping class group and for references.

By assumption, 
the Euler characteristic of $S$ 
is negative
and hence $S$ admits a complete hyperbolic
Riemannian metric of finite volume. A \emph{geodesic lamination}
for such a hyperbolic structure is a \emph{compact} subset
of $S$ which is foliated into simple geodesics. Call a
geodesic lamination $\lambda$ \emph{complete} if its complementary
components are all ideal triangles or once punctured monogons
and if moreover $\lambda$ can be approximated in the \emph{Hausdorff
topology} by simple closed geodesics. The space ${\cal C\cal L}$
of complete geodesic laminations on $S$ equipped with the
Hausdorff topology is compact and metrizable (see
Section 2 of this paper). The mapping class group  
naturally 
acts on the space ${\cal C\cal L}$ as
a group of homeomorphisms. In other words, ${\cal C\cal L}$
is a compact metrizable ${\cal M\cal C\cal G}$-space which
we call the \emph{Furstenberg boundary} of ${\cal M\cal C\cal G}$. 
We show.

\begin{theorem}\label{theorem1}
The action of
the mapping class group of a non-exceptional
surface $S$ of finite type on the space of complete
geodesic laminations on $S$ 
is topologically amenable.
\end{theorem}

As a consequence, the mapping class groups are boundary
amenable. 
Since boundary amenability is passed on to subgroups
(see Chapter 5 of \cite{AR00}), 
as explained above the following corollary is immediate
from Theorem \ref{theorem1} and the work of Higson \cite{H00}.

\begin{coro}\label{corollary1}
The Novikov higher order signature
conjecture holds for any subgroup of 
the mapping class group of a non-exceptional
surface of finite type.
\end{coro}

A simple closed curve on the surface $S$ is called 
\emph{essential} if it is neither contractible nor
freely homotopic into a puncture.
The \emph{curve graph} ${\cal C}(S)$ of the surface $S$ is
a locally infinite
metric graph whose vertices are the free
homotopy classes of essential simple closed
curves on $S$ and where two such
vertices are connected by an edge of length one 
if they can be realized disjointly. By an important
result of Masur and Minsky \cite{MM99}, the curve
graph is hyperbolic in the sense of Gromov. The
mapping class group acts on ${\cal C}(S)$ as a group
of isometries.

The Gromov hyperbolic geodesic metric space ${\cal C}(S)$ 
admits a \emph{Gromov boundary} $\partial 
{\cal C}(S)$
which is a metrizable topological space;
however, it is not locally compact. 
The action of the mapping class group ${\cal M\cal C\cal G}$ on 
${\cal C}(S)$ extends to an action on $\partial {\cal C}(S)$ 
by homeomorphisms. 

An action of a group
on a topological space is called \emph{universally
amenable} if it is amenable with respect 
to every invariant Borel measure class.
As another corollary of Theorem \ref{theorem1}
we obtain the
following result.

\begin{coro}\label{corollary2}
The action of ${\cal M\cal C\cal G}$ 
on $\partial {\cal C}(S)$ is universally amenable.
\end{coro}

Now let $n\geq 2$ and for $i\leq n$ let $G_i$ be
a locally compact 
second countable topological group.
A subgroup $\Gamma$
of the group $G=G_1\times\dots\times G_n$ 
is an \emph{irreducible
lattice} in $G$ if the volume of $G/\Gamma$ with
respect to a Haar measure $\lambda$ 
is finite and if moreover
the projection of $\Gamma$ to each factor
is dense. We allow the groups $G_i$ to be discrete.
Let $X$ be a standard Borel $\Gamma$-space and let
$\mu$ be a $\Gamma$-invariant ergodic probability
measure on $X$. The action of $\Gamma$ on $(X,\mu)$
is called \emph{mildly mixing} if there are no
non-trivial recurrent sets, i.e. if for any
measurable set $A\subset X$ and any sequence
$\phi_i\to \infty$ in $\Gamma$, one has
$\mu(\phi_iA\vartriangle A)\to 0$ only
when $\mu(A)=0$ or $\mu(X-A)=0$.
An ${\cal M\cal C\cal G}$-valued \emph{cocycle}
for this action is a measurable map
$\alpha:\Gamma\times X\to {\cal M\cal C\cal G}$ such that
\[\alpha(gh,x)=\alpha(g,hx)\alpha(h,x)\] for 
all $g,h\in \Gamma$ and $\mu$-almost
every $x\in X$. The cocycle $\alpha$
is \emph{cohomologous} to a cocycle
$\beta:\Gamma\times X\to {\cal M\cal C\cal G}$ 
if there is a measurable map $\phi:X\to {\cal M\cal C\cal G}$
such that $\phi(gx)\alpha(g,x)=\beta(g,x)\phi(x)$ for
all $g\in \Gamma$, $\mu$-almost
every $x\in X$. 

For every 
continuous homomorphism $\rho$ of 
a topological group $H$ into a topological
group $L$, the composition with $\rho$ 
of an $H$-valued cocycle is a 
cocycle with values in $L$.
We use the space of complete geodesic
laminations on $S$ to show the following super-rigidity
result for ${\cal M\cal C\cal G}$-valued cocycles.

\begin{theorem}\label{theorem2}
Let $n\geq 2$, 
let $\Gamma<G_1\times\dots\times G_n$
be an irreducible lattice, let $X$ be a mildly 
mixing $\Gamma$-space and let 
$\alpha:\Gamma\times X\to {\cal M\cal C\cal G}$ 
be any cocycle. Then $\alpha$ is cohomologous
to a cocycle $\alpha^\prime$ with values
in a subgroup $H=H_0\times H_1$ of ${\cal M\cal C\cal G}$
where $H_0$ is virtually abelian and where $H_1$
contains a finite normal subgroup $K$ such that
the projection of $\alpha^\prime$ into $H_1/K$
defines a continuous homomorphism 
$G\to H_1/K$. 
\end{theorem}

We also obtain some geometric information on the mapping
class group. Namely,
an \emph{$L$-quasi-isometric embedding} of a metric
space $(X,d)$ into a metric space $(Y,d)$ is a map
$F:X\to Y$ such that
\[d(x,y)/L-L\leq d(Fx,Fy)\leq Ld(x,y)+L\text{ for all }x,y\in X.\]
The map $F$ is a called an \emph{$L$-quasi-isometry} if 
moreover its image $F(X)$ is \emph{$L$-dense} in $Y$, i.e. if
for every $y\in Y$ there is some $x\in X$ such that $d(Fx,y)\leq L$.
An \emph{$L$-quasi-geodesic} in $(X,d)$ is an $L$-quasi-isometric embedding
of a connected subset of the real line.

The mapping class group 
${\cal M\cal C\cal G}$ of $S$ is finitely generated \cite{I02}
and hence 
a finite symmetric generating set ${\cal G}$  
defines a \emph{word norm} $\vert \,\vert$
on ${\cal M\cal C\cal G}$ 
by assigning to an element $g\in {\cal M\cal C\cal G}$ the
minimal length of a word in ${\cal G}$ which represents $g$.
The word norm $\vert\,\vert$ determines a
distance $d$ on ${\cal M\cal C\cal G}$
which is invariant under the action of ${\cal M\cal C\cal G}$
by left translation 
by defining
$d(g,h)=\vert g^{-1}h\vert$. 
Any two such distance functions
are quasi-isometric. Moreover, if $Z$ is any geodesic
metric space on which ${\cal M\cal C\cal G}$ acts isometrically, 
properly and cocompactly, then 
by a well known result of \v{S}varc-Milnor (see \cite{BH99}),
for every $z\in Z$ the 
orbit map ${\cal M\cal C\cal G}\to Z$ which associates to $g\in 
{\cal M\cal C\cal G}$ the point $gz\in Z$ 
is an ${\cal M\cal C\cal G}$-equivariant quasi-isometry.
In particular,
$Z$ can be viewed as a geometric model for ${\cal M\cal C\cal G}$. 

For the proof of the above results, we construct
a new geometric model for the mapping class group of $S$ in 
this sense and investigate its geometric properties. This
model is a locally finite metric graph ${\cal T\cal T}$ 
whose vertices are the isotopy classes of 
\emph{complete train tracks} 
on $S$ (see \cite{PH92}) 
and where such a train track
$\tau$ is connected to a train track $\sigma$ by a directed
edge of length one if $\sigma$ can be obtained from 
$\tau$ by a single \emph{split}. We define this graph in 
Section 3 and show that it is connected. 
In Section 4 we observe that the mapping class group
of $S$ acts properly and cocompactly as a group
of simplicial isometries on ${\cal T\cal T}$.
As a consequence, ${\cal T\cal T}$ is quasi-isometric to 
${\cal M\cal C\cal G}$. 

Define a \emph{splitting arc} in ${\cal T\cal T}$ to 
be a directed simplicial path 
$\gamma:[0,m]\to{\cal T\cal T}$ for the standard 
simplicial structure on $\mathbb{R}$ whose vertices are 
the integers. This means that 
for every $i< m$ the arc $\gamma[i,i+1]$
is an edge in ${\cal T\cal T}$ connecting the train track 
$\gamma(i)$ to a train track
$\gamma(i+1)$ which can be obtained from $\gamma(i)$ 
by a single split. Since ${\cal T\cal T}$ is a geometric
model for ${\cal M\cal C\cal G}$, 
the following result gives some information
on quasi-geodesics in ${\cal M\cal C\cal G}$. 

\begin{theorem}\label{theorem3}
There is a number $L>0$ such that every splitting arc
in ${\cal T\cal T}$ 
is an $L$-quasi-geodesic.
\end{theorem}

Theorem \ref{theorem3} can be used to investigate various geometric
properties of the mapping class group. 
As an illustration, we include here 
a particularly easy corollary.

Namely, a finite symmetric generating set
${\cal G}_\Gamma$ of 
a finitely generated subgroup  
$\Gamma<{\cal M\cal C\cal G}$
defines a distance function $d_\Gamma$ on $\Gamma$. 
Since we can always extend ${\cal G}_\Gamma$ to a finite
symmetric generating set of ${\cal M\cal C\cal G}$, 
for every distance function $d$ on ${\cal M\cal C\cal G}$ defined
by a word norm 
there is a number $q>0$ such that the natural inclusion
$(\Gamma,d_\Gamma)\to ({\cal M\cal C\cal G},d)$ is $q$-Lipschitz.
However, in general the word norm in $\Gamma$ of an element
$g\in \Gamma$ can not be estimated from above by a constant
multiple of its word norm in ${\cal M\cal C\cal G}$.
The group $\Gamma<{\cal M\cal C\cal G}$ is called 
\emph{undistorted}
if there is a constant $c>1$ such that 
$d_\Gamma(g,h)\leq cd(g,h)$
for all $g,h\in \Gamma$. This is equivalent to stating that
the natural inclusion $(\Gamma,d_\Gamma)
\to ({\cal M\cal C\cal G},d)$ is
a quasi-isometric embedding.  

A \emph{pants decomposition} for $S$ is a collection of $3g-3+k$ simple
closed mutually disjoint \emph{pants curves} 
which decompose $S$ into
$2g-2+k$ \emph{pairs of pants}, i.e. three-holed spheres.
Such a pants decomposition determines a free 
abelian subgroup of
${\cal M\cal C\cal G}$ of rank $3g-3+k$ 
which is generated by the
\emph{Dehn twists} about the pants curves. As an immediate
application of Theorem \ref{theorem3} we obtain a new proof of the
following result of 
Farb, Lubotzky and Minsky \cite{FLM01}.

\begin{coro}\label{corollary3}
For every pants decomposition $P$ for $S$, the free abelian
subgroup of ${\cal M\cal C\cal G}$ of rank $3g-3+k$ which 
is generated by the Dehn twists about the pants curves of $P$
is an undistorted subgroup of ${\cal M\cal C\cal G}$.
\end{coro}

The organization of this paper is as follows.
In Section 2 we introduce the space of complete geodesic
laminations. We also
summarize some results on geodesic laminations
and train tracks from the literature which are used  
throughout the paper.
As mentioned above, in Section 3 we define 
the train track complex ${\cal T\cal T}$, and we 
show that ${\cal T\cal T}$ 
is a connected metric graph.
In Section 4 we show that
the mapping class groups acts properly and cocompactly
as a group of isometries on ${\cal T\cal T}$.

In Section 5 we define
a family of connected subgraphs of the train track
complex, one for each complete train track $\tau$ and for
each complete geodesic lamination \emph{carried} by $\tau$. We show
that each of these subgraphs is isometric with respect to
its intrinsic path metric to a cubical graph contained 
in an euclidean space
of fixed dimension. In Section 6 we show that these ``flat cones''
are uniformly quasi-isometrically embedded in ${\cal T\cal T}$.
This immediately implies Theorem \ref{theorem3} and
Corollary \ref{corollary3}. 
In Section 7 we use the results from the
previous sections to
show Theorem \ref{theorem1}. 
In Section 8 we construct an explicit
\emph{strong boundary} for the mapping class group and 
derive Corollary \ref{corollary2}. The proof of Theorem \ref{theorem2}
is contained in Section 9. 

The results in Section 6 depend in an essential way 
on a technical property
on train tracks which is established in the appendix in Section 10. 
The appendix only uses those results from the literature collected
in Section 2 and is independent from the 
rest of the paper

After we completed this work we obtained the preprint of Kida
\cite{Ki05} who shows with an inductive argument
that the action of the mapping class group on its
Stone \v{C}ech compactification is topologically amenable.
By an observation of Higson and Roe \cite{HR00}, this
is equivalent to boundary amenability
for the mapping class groups. Kida also obtains Corollary \ref{corollary2}.

\section{Train tracks and geodesic laminations}

This introductory section is divided into two subsections.
In the first subsection, we introduce the space of
complete geodesic laminations for an oriented
surface $S$ of genus $g\geq 0$ with $k\geq 0$ punctures 
and where $3g-3+k\geq 2$.
The second part summarizes some properties of train tracks
on $S$ which are used throughout the paper.

\subsection{Complete geodesic laminations}

Fix a complete 
hyperbolic metric $h$ of finite volume on the
surface $S$. 
A \emph{geodesic lamination} for the metric $h$ is
a \emph{compact} subset of $S$ which is foliated into simple
geodesics. Particular geodesic laminations are
simple closed geodesics, i.e. laminations which consist of a
single leaf.

A geodesic lamination $\lambda$ is called \emph{minimal} if each
of its half-leaves is dense in $\lambda$. Thus a geodesic lamination
is minimal if it does not contain any proper \emph{sublamination}, 
i.e. a proper closed subset which is a geodesic lamination.
A simple closed
geodesic is a minimal geodesic lamination. A minimal geodesic
lamination with more than one leaf has uncountably many leaves and
is called \emph{minimal arational}. Every geodesic lamination
$\lambda$ is a disjoint union of finitely many
minimal sublaminations
and a finite number of isolated leaves. Each of the isolated
leaves of $\lambda$ either is an isolated closed geodesic and
hence a minimal component, or it \emph{spirals} about one or two
minimal components. We refer to \cite{CB88,CEG87} for 
a detailed discussion of the structure
of a geodesic lamination.

A geodesic lamination is \emph{finite} if it contains only
finitely many leaves, and this is the case if and only if each
minimal component is a closed geodesic. A geodesic lamination is
\emph{maximal} if it is not a proper sublamination of another
geodesic lamination, and this is the case if and only if all
complementary regions are
ideal triangles or once punctured monogons
\cite{CEG87}.
Note that a geodesic lamination can be both minimal and maximal.

Since each geodesic
lamination is a compact subset of the surface $S$, the space of
all geodesic laminations can be equipped with the restriction
of the \emph{Hausdorff topology} for compact subsets of $S$.
Moreover, the tangent lines of a geodesic lamination define a
compact subset of the \emph{projectivized tangent bundle} $PTS$ of
$S$. We therefore can equip the space of geodesic laminations on
$S$ with the Hausdorff topology for compact subsets of $PTS$.
However, these two topologies coincide
\cite{CB88}, and in the sequel we shall freely use
these two descriptions interchangably.
With this topology, the space of all geodesic laminations
is compact, and it contains the space of all maximal
geodesic laminations as a compact subset \cite{CB88}.

\begin{definition}\label{complete1}
A \emph{complete geodesic lamination} is a 
maximal geodesic lamination which can be approximated in the
Hausdorff topology by simple closed geodesics.
\end{definition}

Since every minimal geodesic lamination can be
approximated in the Hausdorff topology by simple closed geodesics
\cite{CEG87}, laminations which are both maximal and minimal are
complete. There are also complete finite geodesic 
laminations. Namely, let
$P$ be a \emph{pants decomposition} of $S$, i.e. $P$ is the union
of $3g-3+m$ pairwise disjoint simple closed geodesics
which decompose $S$ into $2g-2+m$ planar surfaces of Euler characteristic
$-1$. 
A geodesic lamination $\lambda$ with the following properties is
complete.
\begin{enumerate}
\item The pants decomposition $P$ is the union of the minimal components
of $\lambda$.
\item For each component $Q$ of $S-P$ and every pair 
$\gamma_1\not=\gamma_2$
of boundary geodesics of $Q$ there is a leaf of $\lambda$ contained in $Q$
which spirals
in one direction about $\gamma_1$, in the other direction about
$\gamma_2$.
\item Each component $Q$ of $S-P$ containing a puncture of $S$
contains a leaf of $\lambda$ which goes around a 
puncture and spirals in both directions about the same boundary
component of $Q$. 
\item For every component $\gamma$ of $P$, the leaves of $\lambda$
which spiral about $\gamma$ from each side of $\gamma$ define
opposite orientations near $\gamma$ as in Figure A.
\end{enumerate}
\begin{figure}[ht]
\includegraphics{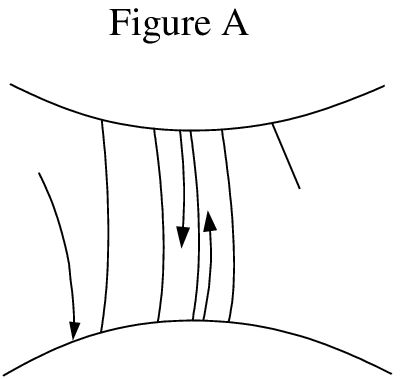}
\end{figure}
The forth condition guarantees that we can approximate
$\lambda$ by smooth simple closed curves which pass from
a leaf spiraling about a component $\gamma$ of $P$
from one side to a leaf 
spiraling about $\gamma$ from the other side and whose
tangents are close to the tangent lines of $\lambda$.

Since maximal geodesic laminations form a closed subset of
the space of all geodesic laminations,
the set ${\cal C\cal L}$ of all
complete geodesic laminations is a closed subset of the space of
all geodesic laminations
with the Hausdorff topology.
Thus ${\cal C\cal L}$ is a compact topological
space. A \emph{simple geodesic multi-curve} is a disjoint
union of simple closed geodesics. We have.

\begin{lemma}\label{compapp}
\begin{enumerate}
\item 
Every geodesic lamination on $S$ which can be approximated
in the Hausdorff topology by simple geodesic multi-curves
is a sublamination of a complete
geodesic lamination.
\item Finite complete geodesic laminations are dense in the
space of all complete geodesic laminations.
\end{enumerate}
\end{lemma}
\begin{proof} Let $\lambda$ be any geodesic
lamination on $S$ which can be approximated
in the Hausdorff topology by a sequence $\{c_i\}$ of simple
geodesic multi-curves. For each $i$ 
let $\mu_i$ be a complete finite
geodesic lamination which contains $c_i$ as a
union of minimal components. Such a geodesic lamination
exists by the above example 
since each simple geodesic multi-curve is a subset
of a pants decomposition for $S$.
By passing to a subsequence we may assume that the
geodesic laminations $\mu_i$ converge as $i\to \infty$ in
the Hausdorff topology to a complete geodesic lamination
$\mu$. Since $c_i\to \lambda$ $(i\to \infty)$ in the
Hausdorff topology, the lamination $\mu$
contains $\lambda$ as a sublamination. This shows the
first part of the lemma.

Since by definition a complete geodesic lamination can be 
approximated in the Hausdorff topology by simple closed geodesics,
the second part follows from the same argument together with the
observation that a geodesic lamination which contains 
a complete geodesic lamination $\lambda$ as a sublamination
coincides with $\lambda$.
\end{proof}

\subsection{Complete train tracks}

A \emph{train track} on $S$ is an embedded
1-complex $\tau\subset S$ whose edges (called \emph{branches}) are
smooth arcs with well-defined tangent vectors at the endpoints. At
any vertex (called a \emph{switch}) the incident edges are
mutually tangent. Through each switch there is a path of class
$C^1$ which is embedded in $\tau$ and contains the switch in its
interior. In particular, the half-branches which are incident on a
fixed switch of $\tau$ are
divided into two classes according to the orientation of the
inward pointing tangent at the switch.
Each closed curve component of $\tau$ has a unique bivalent
switch, and all other switches are at least trivalent. The
complementary regions of the train track have negative Euler
characteristic, which means that they are different from discs
with $0,1$ or $2$ cusps at the boundary and different from annuli
and once-punctured discs with no cusps at the boundary. 
A train track is called \emph{maximal} if each of its
complementary components either is a trigon, i.e. 
a topological disc with three cusps at the boundary,
or a once punctured monogon, i.e. a once punctured
disc with one cusp at the boundary. We always
identify train tracks which are isotopic.

Train tracks were invented by Thurston \cite{T79} and
provide a powerful tool for the investigation of surfaces and
hyperbolic 3-manifolds. A detailed account on train tracks can be
found in the book \cite{PH92} of Penner with Harer which we use as
our main reference. The more recent unpublished manuscript \cite{M03}
of Mosher contains a discussion of train tracks from a somewhat
different viewpoint, however it will not be used in this paper.

A \emph{trainpath} on a train track $\tau$ is a
$C^1$-immersion $\rho:[m,n]\to \tau\subset S$
which maps each interval
$[k,k+1]$ $(m\leq k\leq n-1)$ onto a branch of $\tau$. The integer
$n-m$ is then called the \emph{length} of $\rho$.
We sometimes identify a trainpath on $S$
with its image in $\tau$. Each complementary region of $\tau$
is bounded by a finite number of trainpaths
which either are simple closed curves or terminate
at the cusps of the region.

A train track is called \emph{generic} if all switches are
at most trivalent.
The train track $\tau$ is called \emph{transversely recurrent} if
every branch $b$ of $\tau$ is intersected by an embedded simple
closed curve $c=c(b)\subset S$ which intersects $\tau$
transversely and is such that $S-\tau-c$ does not contain an
embedded \emph{bigon}, i.e. a disc with two corners at the
boundary.

\bigskip

{\bf Remark:} 1) We chose to use transversely recurrent
train tracks for our purpose even though this property is 
nowhere needed.
The main reason for using transversely
recurrent train tracks is convenience of reference
to the existing literature. 

2) Throughout the paper,
we require every train track 
to be generic, and this is indeed necessary for many of our
constructions. Unfortunately this leads to a slight inconsistency
of our terminology with the terminology found in the
literature. 
 
\bigskip

Every generic train track $\tau$ on $S$ is
contained in a closed subset $A$ of $S$ with dense interior and
piecewise smooth boundary which is foliated by smooth arcs
transverse to the branches of $\tau$. These arcs are called
\emph{ties} for $\tau$, and the set $A$ is called a \emph{foliated
neighborhood} of $\tau$ (even though $A$ is not a neighborhood of
$\tau$ in the usual sense). Each of the ties intersects $\tau$ in
a single point, and the switches of $\tau$ are the intersection
points of $\tau$ with the singular ties, i.e. the ties which are
not contained in the interior of a foliated rectangle $R\subset
A$. Collapsing each tie to its intersection point with $\tau$
defines a map $F:A\to \tau$ of class $C^1$ which we call a
\emph{collapsing map}. The map $F$ is the restriction to
$A$ of a map $\tilde F:S\to S$ of class $C^1$ which is homotopic
to the identity.

A train track or a geodesic lamination $\lambda$ 
is \emph{carried} by a transversely
recurrent train track $\tau$ if there is an isotopy $\phi$ of $S$
such that $\lambda$ is contained in a foliated
neighborhood $A$ of $\phi(\tau)$ and is transverse to the ties.
This is equivalent to stating that there is a map
$G:S\to S$ of class $C^1$ which is homotopic to the identity and
which maps $\lambda$ to $\tau$ in such a way that the restriction
of the differential of $G$ to every tangent space of 
$\lambda$ is non-singular. Note that this makes sense since
a train track has a tangent line everywhere 
(compare Theorem 1.6.6 of \cite{PH92}). We then call the
restriction of $G$ to $\lambda$ a \emph{carrying map}.
Every train track $\tau$ which
carries a maximal geodesic lamination is necessarily maximal
\cite{PH92}.
The set of geodesic
laminations which are carried by a transversely recurrent
train track $\tau$ is a closed subset of the
space of all geodesic laminations on $S$ with respect to 
the Hausdorff 
topology (Theorem 1.5.4 of \cite{PH92}, or see \cite{CB88}
for the same result for train tracks which are not necessarily
transversely recurrent).

A \emph{transverse measure} on a train track $\tau$ is a
nonnegative weight function $\mu$ on the branches of $\tau$
satisfying the \emph{switch condition}:
For every switch $s$ of $\tau$, the half-branches incident at $s$\
are divided into two classes according to the orientation
of their inward pointing tangent at $s$. 
We require that
the sums of the weights
over all branches in each of the two classes coincide.
The train track is called
\emph{recurrent} if it admits a transverse measure which is
positive on every branch. We call such a transverse measure $\mu$
\emph{positive}, and we write $\mu>0$. The set $V(\tau)$ of all
transverse measures on $\tau$ is a closed convex cone in a linear
space and hence topologically it is a closed cell. For every
recurrent train track $\tau$, positive measures define the
interior of the convex cone $V(\tau)$.
A train track $\tau$ is called \emph{birecurrent} if
$\tau$ is recurrent and transversely recurrent.

A \emph{measured geodesic lamination} is a geodesic lamination
equipped with a transverse Borel measure of full support which is
invariant under holonomy \cite{PH92}. A measured geodesic
lamination can be viewed as a locally finite Borel measure on the
space $S_\infty$ of unoriented
geodesics in the universal covering ${\bf
H}^2$ of $S$ which is invariant under the action of $\pi_1(S)$ and
supported in the closed set of geodesics whose
endpoints are not separated under the action of 
$\pi_1(S)$. The
weak$^*$-topology on the space of $\pi_1(S)$-invariant locally
finite Borel measures on $S_\infty$ then restricts to a natural
topology on the space ${\cal M\cal L}$ of all measured geodesic
laminations.

If $\lambda$ is a geodesic lamination with transverse measure
$\mu$ and if $\lambda$ is carried by a train track $\tau$, then
the transverse measure $\mu$ induces via a carrying map
$\lambda\to \tau$ a transverse measure on $\tau$, and every
transverse measure on $\tau$ arises in this way (Theorem 1.7.12 in
\cite{PH92}). Moreover, if $\tau$ is maximal and birecurrent, then
the set of measured geodesic laminations which correspond to
positive transverse measures on $\tau$ in this way is an open
subset $U$ of ${\cal M\cal L}$ (Lemma 3.1.2 in \cite{PH92}).

We use measured geodesic laminations to relate maximal
birecurrent train tracks to complete geodesic laminations.
Note that the first
part of the following lemma is well known and reflects the fact
that the space ${\cal C\cal L}$ of complete geodesic laminations
is totally disconnected. We refer to
\cite{Bo97} and \cite{ZB04} for more details and references about
the Hausdorff topology on the space of all geodesic laminations.

\begin{lemma}\label{recchar} 
Let $\tau$ be a maximal
transversely recurrent train track. Then the set of all complete
geodesic laminations on $S$ which are carried by $\tau$ is open
and closed in ${\cal C\cal L}$. This set is non-empty if and only if
$\tau$ is recurrent.
\end{lemma}
\begin{proof}
Let $\tau$ be a maximal transversely recurrent train
track on $S$. Then the subset of ${\cal C\cal L}$ of all complete
geodesic laminations which are carried by $\tau$ is closed 
(Theorem 1.5.4 of \cite{PH92}). 

On the other
hand, if $\lambda\in {\cal C\cal L}$ 
is carried by $\tau$ then after possibly modifying 
$\tau$ with an isotopy we may assume that
there is a foliated
neighborhood $A$ of $\tau$ with collapsing map $F:A\to \tau$
which contains
$\lambda$ in its interior and such that $\lambda$ is transverse
to the ties of $A$. Since the space of tangent lines
$PT\lambda$ 
of $\lambda$ is a compact subset of the projectivized
tangent bundle $PTS$ of $S$, there is a neighborhood 
$U$ of $PT\lambda$ in $PTS$ 
which is mapped by the canonical projection
$PTS\to S$ into $A$ and such that the 
differential of $F$ is non-singular on each line in $U$. 
As a consequence, a geodesic lamination whose 
space of tangent lines is contained in $U$ is carried by $\tau$.
Since the Hausdorff topology on the space of geodesic
laminations coincides with the Hausdorff topology for their
projectivized tangent bundles,
the set of all such geodesic laminations which are moreover complete
is a neighborhood of $\lambda$ in ${\cal C\cal L}$.
This shows that the set of all complete geodesic laminations
which are carried by $\tau$ is open and closed in ${\cal C\cal L}$.

To show the second part of the lemma,
let $\tau$ be a maximal transversely recurrent train
track on $S$ which carries a complete geodesic
lamination $\lambda\in {\cal C\cal L}$. 
To see that $\tau$ is recurrent,
let $\{c_i\}$ be a sequence of simple
closed geodesics which approximate $\lambda$ in the Hausdorff
topology. 
By the above consideration, the curves $c_i$ are carried by
$\tau$ for all sufficiently large $i$. Let $F_i:c_i\to \tau$
be a carrying map. Then for each
sufficiently large $i$, the curve $c_i$ defines a \emph{counting
measure} on $\tau$ by associating to each branch $b$ of $\tau$ the
number of components of $F_i^{-1}(b)$. This counting measure
satisfies the switch condition.
Now a carrying map $F:\lambda\to
\tau$ maps $\lambda$ \emph{onto} $\tau$ since $\lambda$ is
maximal, and hence the same is true for the carrying map $F_i:c_i\to
\tau$ provided that $i$ is sufficiently large.  Thus
for sufficiently large $i$, the counting measure
defined by $c_i$ is positive and therefore 
$\tau$ is recurrent. 

To show that a maximal birecurrent train track carries a complete
geodesic lamination, recall that 
a geodesic
lamination which is both minimal and maximal is complete 
and supports a transverse  measure \cite{CEG87}. 
Since the set of all measured geodesic laminations
carried by a maximal birecurrent train track 
$\tau$ has non-empty interior, it is enough to show 
that the set of measured geodesic
laminations whose support is minimal and maximal is dense in the
space ${\cal M\cal L}$ of all measured geodesic laminations.

However, as was pointed out to me by McMullen, this follows for
example from the work of Kerckhoff, Masur and Smillie
\cite{KMS86}. Namely, the \emph{vertical foliation} of a
\emph{quadratic differential} defines a measured geodesic
lamination, and every measured geodesic lamination is of this form
\cite{T79}. The set of measured geodesic laminations which arise from
quadratic differentials with only simple zeros is dense in the
space of all measured geodesic laminations (compare the discussion
in \cite{KMS86}). Given such a quadratic differential $q$, there
are only countably many $\theta\in [0,2\pi]$ such that the
vertical foliation of $e^{i\theta}q$ contains a compact vertical
arc connecting two zeros of $q$. On the other hand, for each
quadratic differential $q$ the set of all $\theta\in [0,2\pi]$
such that the vertical foliation of $e^{i\theta}q$ is
\emph{uniquely ergodic} and hence minimal has full Lebesgue
measure. Thus there is a dense set of points $\theta\in [0,2\pi]$
with the property that the vertical foliation of $e^{i\theta}q$ is
both minimal and maximal.

As a consequence, a maximal birecurrent train track
carries a geodesic lamination which is both
maximal and minimal and hence complete.
\end{proof}

\begin{definition}\label{complete2}
A train track $\tau$ on $S$ which is generic, maximal and 
birecurrent 
is called \emph{complete}. 
\end{definition}

Note that this definition of a complete train track slightly
differs from the one in \cite{PH92} since we require a complete train
track to be generic. By Lemma \ref{recchar}, 
a generic transversely recurrent
train track is complete if and only if it carries a complete
geodesic lamination.

There are two basic ways to modify a complete train track
to another complete train track. 
Namely, let $\tilde b$ be a
half-branch of a generic train track $\tau$ and let $v$ be
a trivalent switch of $\tau$ on which $\tilde b$ is incident.
Then $\tilde b$ is called \emph{large} if any immersed
arc of class $C^1$ in $\tau$ through $v$ intersects
$\tilde b$. The branch $b$ of $\tau$ containing $\tilde b$ is 
called \emph{large at $v$}. A half-branch which is not large
is called \emph{small}. 
A branch $b$ in a generic train track $\tau$ 
\emph{large} if each of its two half-branches is large. A large
branch $b$
is necessarily incident on two distinct switches, and it
is large at both of them. A branch is called \emph{small} if each
of its two half-branches is small. A branch is called \emph{mixed}
if one of its half-branches is large and the other half-branch is
small (for all this, see \cite{PH92} p.118).

The two basic ways to modify a complete train track $\tau$
to another complete train track are as follows. 
First, we can \emph{shift} $\tau$
along a mixed branch to a train track $\tau^\prime$ as shown in
Figure B. 
\begin{figure}[ht]
\includegraphics{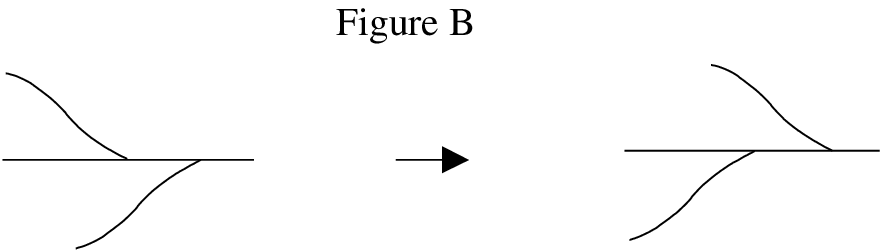}
\end{figure}
If $\tau$ is complete then the same is true for
$\tau^\prime$. Moreover, a train track or a geodesic
lamination is carried
by $\tau$ if and only if it is carried by $\tau^\prime$ (see
\cite{PH92} p.119). In particular, the shift $\tau^\prime$ of
$\tau$ is carried by $\tau$. Note that there is a natural
bijection of the set of branches of $\tau$ onto the set of
branches of $\tau^\prime$.

Second, if $e$ is a large branch of $\tau$ then we can perform a
right or left \emph{split} of $\tau$ at $e$ as shown in Figure C.
\begin{figure}[ht]
\includegraphics{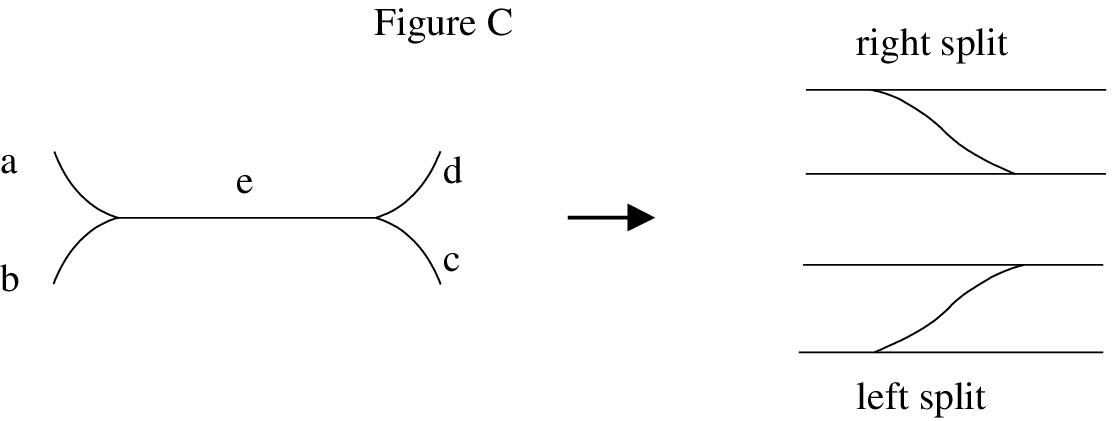}
\end{figure}
Note that a right split at $e$ is uniquely determined by
the orientation of $S$ and does not depend on the orientation of
$e$. Using the labels in the figure, in the case of a right split
we call the branches $a$ and $c$ \emph{winners} of the split, and
the branches $b,d$ are \emph{losers} of the split. If we perform a
left split, then the branches $b,d$ are winners of the split, and
the branches $a,c$ are losers of the split (see \cite{MM99}). The
split $\tau^\prime$ of a train track $\tau$ is carried by $\tau$.
There is a natural bijection of
the set of branches of $\tau$ onto the set of branches of
$\tau^\prime$ which maps the branch $e$ to the diagonal $e^\prime$
of the split. The split of a generic maximal transversely recurrent
train track is generic, maximal and transversely recurrent.
If $\tau$ is complete and if $\lambda\in {\cal C\cal L}$ is
carried by $\tau$, then there is a unique choice of a right or
left split of $\tau$ at $e$ with the property that the split track
$\tau^\prime$ carries $\lambda$. We call $\tau^\prime$ the
\emph{$\lambda$-split} of $\tau$.
By Lemma \ref{recchar}, the train track
$\tau^\prime$ is recurrent and hence complete. In particular, a
complete train track $\tau$ can always be split at any large
branch $e$ to a complete train track $\tau^\prime$; however there
may be a choice of a right or left split at $e$ such that the
resulting track is not recurrent any more (compare p.120 in
\cite{PH92}). The reverse of a split is called a \emph{collapse}.

\section{The complex of train tracks}

Let as before $S$ be 
an oriented surface of
genus $g\geq 0$ with $k\geq 0$ punctures and where $3g-3+k\geq 2$.
Using the notations from Section 2, 
define 
${\cal T\cal T}$ to be the directed graph whose set of vertices
is the set ${\cal V}({\cal T\cal T})$ 
of isotopy classes of complete train tracks on $S$ and whose
edges are determined as follows. The train track $\tau$ is connected
in ${\cal T\cal T}$ to the train track $\tau^\prime$ by a directed edge
if and only if $\tau^\prime$ can be obtained from $\tau$ by a single
right or left split. The goal of this section is to 
show that ${\cal T\cal T}$ is connected.

For this note that 
transversely recurrent train tracks can be viewed as finite
combinatorial 
approximations of geodesic laminations (which in general have
uncountably many leaves). Indeed, Theorem 1.6.5 of \cite{PH92} 
shows that a geodesic lamination $\lambda$ 
can be approximated in the
Hausdorff topology by transversely recurrent train tracks
which moreover carry $\lambda$. This does not imply, however,
that a transversely recurrent 
train track $\tau$ which is close to $\lambda$ in the Hausdorff topology
necessarily carries $\lambda$. For example, 
$\tau$ could have a very
short branch which makes a sharp turn in the 
wrong direction. In fact, it
is not difficult to see that a finite complete geodesic lamination
$\lambda$ can be approximated in the Hausdorff
topology by transversely recurrent train tracks which do not
carry $\lambda$.

Since carrying is a relation determined
by maps of class $C^1$, with restrictions on the tangent map,
to relate an approximation of a geodesic lamination by
train tracks to carrying we have to use 
approximations 
in the Hausdorff topology for compact subsets of the projectivized
tangent bundle $PTS$ of $S$. To this end, we use
geometric realizations of train tracks. Note that this
is consistent with the fact that a geodesic lamination 
is a geometric representative of a purely topological object, defined
independently of the choice of a complete hyperbolic metric
of finite volume on $S$.

The representative of a train track $\tau$ which we use is its
\emph{straightening} with
respect to a hyperbolic metric on $S$. This straightening
is the immersed
graph in $S$ whose vertices are the switches of $\tau$ and whose
edges are the unique geodesic arcs which are homotopic with fixed
endpoints to the branches of $\tau$.
The tangent
lines of the straightening of $\tau$ then define a closed subset
of the projectivized tangent bundle $PTS$ of $S$. The hyperbolic
metric on $S$ naturally induces a Riemannian metric and hence a
distance function on $PTS$. For a number $\epsilon>0$ we say that
the train track $\tau$ \emph{$\epsilon$-follows} the geodesic
lamination $\lambda$ if the tangent lines of the straightening of
$\tau$ are contained in the $\epsilon$-neighborhood of the
projectivized tangent bundle $PT\lambda$ of $\lambda$ and if
moreover the straightening of every trainpath on $\tau$ is a
piecewise geodesic whose exterior angles at the breakpoints are
not bigger than $\epsilon$ (here a vanishing exterior angle means
that the arc is smooth). The train track $\tau$ is \emph{$a$-long}
for a number $a>0$ if the length of every edge of the straightening of 
$\tau$ is at least $a$.

The following technical lemma 
is a strengthening of Theorem 1.6.5 of \cite{PH92} which
follows from the same line of arguments. It 
shows that a geodesic
lamination $\lambda$ can be approximated in the $C^1$-topology by
train tracks which carry $\lambda$.

\begin{lemma}\label{following}
There is a number $a>0$ with the following property.
Let $\lambda$ be any geodesic lamination on $S$.
Then for every $\epsilon>0$ there is an $a$-long 
generic transversely recurrent train
track $\tau$ which carries $\lambda$ and which $\epsilon$-follows
$\lambda$.
\end{lemma}
\begin{proof}
By the collar theorem for hyperbolic surfaces (see \cite{B92}), a complete
simple geodesic on $S$ which is contained in a compact
subset of $S$ does not enter deeply into a cusp. This means that
there is a compact bordered subsurface $S_0$ of $S$ which contains
the $1$-neighborhood of every geodesic lamination on $S$. Let
$k>0$ be the maximal number of branches of any train track on $S$
(this only depends on the topological type of $S$) and let $a\in
(0,1)$ be sufficiently small that $(5k+2)a$ is smaller than
the
smallest length of any non-contractible closed curve in $S_0$.

For $\epsilon >0$ there is a number $\delta_0=\delta_0(\epsilon)\in
(0,a/2)$ with the following property. Let $\gamma$ be a geodesic
line in ${\bf H}^2$ and let $\zeta$ be an arc in ${\bf H}^2$ 
which consists of two geodesic  
segments of length at least $a$. Assume that
$\zeta$ is contained in the $\delta_0$-neighborhood of $\gamma$
and that the length of $\zeta$ does not exceed the sum of the
distance of its endpoints and $4\delta_0$.
Then the tangent lines
of $\zeta$ are contained in the $\epsilon$-neighborhood of the
tangent lines of $\gamma$ as subsets of the projectivized tangent
bundle of ${\bf H}^2$, and the exterior angle at the breakpoint of
$\zeta$ is at most $\epsilon$.

We now use the arguments of Casson and Bleiler as
described in the proof of Theorem 1.6.5 of \cite{PH92}.
Let $\lambda$ be a geodesic lamination on $S$ and  
for $\delta>0$ let $N_\delta$ be the closed $\delta$-neighborhood of
$\lambda$ in $S$. This neighborhood lifts to the
$\delta$-neighborhood $\tilde N_\delta$ of the lift
$\tilde\lambda$ of $\lambda$ to the hyperbolic plane
${\bf H}^2$. For sufficiently small $\delta$, say for
all $\delta<\delta_1$, 
each component of ${\bf H}^2-\tilde
\lambda$ contains precisely one component of ${\bf H}^2-\tilde
N_\delta$. These components are polygons whose sides are arcs of
constant geodesic curvature, and this curvature tends to $0$ with
$\delta$.

Let $\delta<\min\{\delta_0,\delta_1\}$ 
be sufficiently small that $N_\delta$ can be
foliated by smooth vertical arcs which are transverse to the
leaves of $\lambda$ and of length smaller than $\delta_0$. We may
assume that the arcs through the finitely many corners of the
complementary components of $N_\delta$ are geodesics (p.74-75 in
\cite{PH92}). We call this foliation of $N_\delta$ the
\emph{vertical foliation}, and we denote it by ${\cal F}^\perp$.
The singular leaves of ${\cal F}^\perp$, i.e. the leaves through
the corners of $S-N_\delta$, decompose $N_\delta$ into closed
foliated rectangles with embedded interior.
Two opposite sides of these rectangles are
subarcs of singular leaves of the vertical foliation and
hence they are geodesics. The other
two sides are arcs of constant curvature. The intersection of any
two distinct such rectangles is contained in the
singular leaves of ${\cal F}^\perp$. Moreover,
the rectangles are projections to $S$ of \emph{convex} subsets of the
hyperbolic plane. 
After a small adjustment near the corners of the complementary
components of $N_\delta$, collapsing each leaf of
${\cal F}^\perp$ to a suitably chosen point in its interior
defines
a train track $\tau$ on $S$ and a map $F:N_\delta\to \tau$
of class $C^1$ which is homotopic to the identity and whose
restriction to $\lambda$ is a carrying map $\lambda\to \tau$. The
switches of $\tau$ are precisely the collapses of the singular
leaves of $\tau$, and each branch of $\tau$ is a
collapse of one of the foliated rectangles. The set $N_\delta$ is a
foliated neighborhood of $\tau$. 
Every edge of the straightening
of $\tau$ is a geodesic arc which is contained in one of the
rectangles and connects the two sides of the rectangle contained
in leaves of ${\cal F}^\perp$. Via slightly changing $\delta$ we
may assume that $\tau$ is generic; this is equivalent to saying
that each singular leaf of ${\cal F}^\perp$ contains precisely one
corner of a component of $S-N_\delta$. The train track $\tau$ is
transversely recurrent \cite{PH92}.

Define the \emph{length} $\ell(R)$ of
a foliated rectangle $R\subset N_\delta$  
as above to be the \emph{intrinsic} distance
in $R$ between its two sides which are contained
in ${\cal F}^\perp$.  Since a rectangle $R$ lifts to a convex
subset $\tilde R$ of the hyperbolic plane, this length is just the
distance between the two geodesic sides of $\tilde R$.
Let $R_1\not=R_2$ be two rectangles which intersect
along a nontrivial subarc $c$ of a singular leaf of ${\cal F}^\perp$.
Let $a_i$ be the side of $R_i$ opposite to the side containing $c$.
Then $R_1\cup R_2$ contains
a geodesic segment $\gamma$
which connects $a_1$ to $a_2$ and intersects $c$.
The segment $\gamma$ is a subarc of a leaf of $\lambda$, and its
length is not smaller than $\ell(R_1)+\ell(R_2)$. Hence the
length of any curve in $S$ with one endpoint in $a_1$ and
the second endpoint in $a_2$
which is homotopic to $\gamma$ relative
to $a_1\cup a_2$ is not smaller than $\ell(R_1)+ \ell(R_2)-2\delta_0$.
On the other hand, the length of a geodesic arc contained in $R_i$ and
connecting $a_i$ to $c$ which is homotopic to $\gamma\cap R_i$ relative
to $a_i\cup c$ is not bigger than $\ell(R_i)+2\delta_0$.

Let $\zeta$ be the subarc of the straightening
of $\tau$ which consists of the straightening of the collapses of the
rectangles $R_1$ and $R_2$. By convexity,
$\zeta\cap R_i$ is contained in the
$\delta_0$-neighborhood of the geodesic arc $\gamma\cap R_i$, and
the length of $\zeta$
does not exceed the sum of the 
length of the geodesic homotopic to $\zeta$
with fixed endpoints and $4\delta_0$.
Thus by the choice of $\delta_0$, the
train track $\tau$ $\epsilon$-follows $\lambda$ provided that the
length of each edge from the straightening of $\tau$
is at least $a$, and this is the case if the length of each of
the rectangles is at least $a$.

We now successively modify the train track $\tau$ to
a train track $\tau^\prime$ which is embedded in $N_\delta$ and
transverse to the ties, which 
carries $\lambda$ and such that the length of each edge from the
straightening of $\tau^\prime$ is at least $a$ as follows.
Call a rectangle $R$ \emph{short} if its length is at
most $a$. Let $n\geq 0$ be the number of short rectangles in
$N_\delta$. If $n=0$ then
$\tau$ is $a$-long, so
assume that $n>0$. We remove from $N_\delta$ a suitably chosen
geodesic arc $\beta_0$ so that $N_\delta-\beta_0$ is
partitioned into rectangles which are foliated
by the restriction of the vertical foliation ${\cal F}^\perp$
and for which the
number of short rectangles
is at most $n-1$.

This \emph{unzipping} of the train track $\tau$ (\cite{PH92}
p.74-75) is done as follows. Let $R\subset N_\delta$ be a short
rectangle. The boundary $\partial R$ of $R$ contains a
corner $x$ of a complementary component $T$ of $N_\delta$. This
corner projects to a switch $F(x)$ of $\tau$, and the collapse of
$R$ is incident on $F(x)$. 
Let $\tilde T$ be a lift of the complementary component
$T$ to ${\bf H}^2$ and let $\tilde x$
be the lift of $x$ to $\tilde T$. Let $\tilde \lambda$ be the lift
of $\lambda$ to ${\bf H}^2$; then $\tilde T$ is contained in a
unique component $C$ of ${\bf H}^2-\tilde \lambda$. The point
$\tilde x$ is at distance $\delta$ to two frontier
leaves of $C$. Since every complementary
component of $\tilde \lambda$ contains precisely one
complementary component of ${\bf H}^2-\tilde N_\delta$,
these leaves have a common endpoint
$\xi$ in the ideal boundary of ${\bf H}^2$. There is a unique
geodesic ray $\tilde \beta$ in ${\bf H}^2$ 
which connects $\tilde x$ to $\xi$,
and this ray is entirely contained in the intersection of $\tilde
N_\delta$ with the component $C$ of ${\bf H}^2-\tilde \lambda$.
The projection of $\tilde \beta$ to $S$ is a one-sided infinite
simple geodesic $\beta$ beginning at $x$ which is contained in
$N_\delta$ and is disjoint from $\lambda$. We may assume that
$\beta$ is parametrized by arc length and is everywhere transverse
to the leaves of the foliation ${\cal F}^\perp$. In particular,
the collapsing 
map $F$ maps $\beta$ up to parametrization to a one-sided
infinite trainpath on $\tau$.

We claim that the finite subarc $\beta_0$ of $\beta$ of length
$5ka$ which begins at $x$ is mapped by $F$ injectively into
$\tau$. For this note that by the choice of 
the constant $a>0$, the length of
the union of $\beta_0$ with any leaf of the vertical foliation
${\cal F}^\perp$ is smaller than the minimal length of a
non-contractible closed curve in the compact surface $S_0\supset
N_\delta$. If $\beta_0$ intersects the same leaf of 
the vertical foliation twice then there is a subarc
$\beta_1$ of $\beta_0$ whose 
concatenation with a subarc of a leaf of ${\cal F}^\perp$
is a closed curve $c$ of length smaller than $5ka+a$ 
which is contained in $N_\delta\subset S_0$.
By the choice of the constant $a$, the curve $c$ 
is contractible in $S$. On the other hand,
$c$ is freely homotopic to the image of $\beta_1$ under 
the collapsing map $F$ and hence to a subarc of a trainpath
on $\tau$ which begins and ends at the same point (the induced
orientations of its tangent line at that point may
be opposite).
By the definition of a train track,
such a curve is homotopically nontrivial in $S$ which is 
a contradiction.

Therefore $\beta_0$ intersects each leaf of the
vertical foliation ${\cal F}^\perp$ at most once and hence it is
mapped by $F$ injectively into $\tau$. The length of the
intersection of $\beta_0$ with a rectangle $R^\prime$ from our
system of rectangles is at most $\ell(R^\prime)+2\delta_0\leq
\ell(R^\prime)+a$. On the other hand, there are at most $k$
distinct rectangles in our system of rectangles and hence
$\beta_0$ intersects the interior of a rectangle of length at least $4a$.

Let $t\in [0,5ka]$ be the infimum of all numbers $s>0$ such that
$\beta(s)$ is contained in the interior of 
a rectangle $\hat R$ of length at least
$4a$. Then $F(\beta(t))$ is a switch in $\tau$, and $\beta(t+2a)$
is an interior point of $\hat R$. We may assume that the
leaf of the vertical foliation through $\beta(t+2a)$ is a geodesic.
Then this leaf subdivides the rectangle $\hat R$ into
two foliated rectangles of length at least $a$.
Cut $N_\delta$ open along $\beta[0,t+2a)$ as shown in Figure D.
\begin{figure}[ht]
\includegraphics{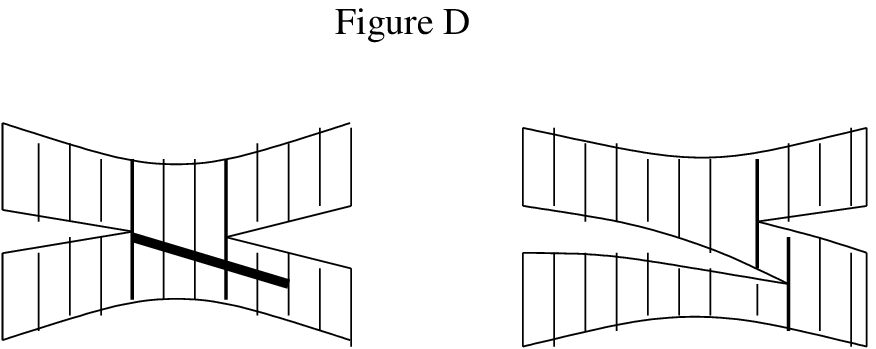}
\end{figure}
Since $\beta$ is transverse to the
vertical foliation ${\cal F}^\perp$, the foliation ${\cal F}^\perp$ 
restricts to a
foliation of $N_\delta-\beta[0,t+2a)$. Now $\beta$ is disjoint
from $\lambda$, and hence 
collapsing each leaf of
${\cal F}^\perp\vert N_\delta-\beta[0,t+2a)$ 
to a suitably chosen
point in its interior yields a train track
$\sigma$ which carries $\lambda$.
The switches of $\sigma$ are the collapses of the
singular leaves of ${\cal F}^\perp\vert N_\delta -\beta[0,t+2a)$.
Thus $N_\delta-\beta[0,t+2a)$ is partitioned into a
finite number of rectangles as before. 
Each of these rectangles lifts to a convex subset of the hyperbolic
plane. The restriction to $\sigma$
of the collapsing map $F$ defines a map
$\sigma\to\tau$ with the property that the preimage of every
switch $w\not= F(x)$ of $\tau$ consists of exactly one switch of
$\sigma$. Note that $\sigma$ is transversely recurrent and carries
$\lambda$ \cite{PH92}.

We claim that the number of short rectangles of $N_\delta-\beta[0,t+2a)$
is at most $n-1$. For this
let $e$ be the branch
of $\tau$ which is incident and large at the switch
$F(x)$ and write $R_0=F^{-1}(e)$.
Also denote by $R_{-1},R_{-2}$ the (not necessarily distinct) 
rectangles which are mapped
by $F$ onto the branches of $\tau$ which are incident and small
at $F(x)$. These rectangles are properly contained
in rectangles $A(R_{-1}),A(R_{-2})$ in $N_\delta-\beta[0,t+2a)$.
There is a short rectangle among the rectangles 
$R_{-2},R_{-1},R_0$.

Assume that for some $q\geq 0$,
the path $F(\beta[0,t+2a))$ on $\tau$
passes through the switches $F(x),v_1,\dots,v_q=F(\beta(t))$ of
$\tau$ in this order. If $q=0$ then 
we have
$\ell(A(R_{-i}))\geq \ell(R_{-i})
+a>a$ for $i=1,2$. Since 
every rectangle in $N_\delta$ different from $R_i$ for $i=-2,-1,0$
is also a rectangle in $N_\delta-\beta[0,t+2a)$, in this
case the claim is obvious. 

In the case $q\geq 1$, we obtain the same conclusion
as follows. For $i=1,2$ let as before $A(R_{-i})$
be the rectangle in $N_\delta-\beta_0[t,2a)$ which contains
$R_{-i}$. Note that by the definition of length and by convexity,
we have
\begin{equation}\label{increase}
\ell(A(R_{-i}))\geq \ell(R_{-i})+\ell(R_0).
\end{equation}
We extend the map $A$ to a 
bijection from the rectangles in $N_\delta$ to the rectangles
in $N_\delta-\beta[0,t+2a)$ and 
compare the length of a rectangle $\tilde R$ in $N_\delta$
to the length of $A(\tilde R)$ in
$N_\delta-\beta_0[0,t+2a)$ as follows.
For $1\leq i\leq q-1$ let
$R_i$ be the rectangle in $N_\delta$ which is mapped by $F$ to the
branch of $\tau$ which is
incident on the switches $v_{i}$ and $v_{i+1}$ and
is crossed through by
the path $F(\beta[0,t+2a))$.
For each $i\in \{1,\dots,q-1\}$, the point
$F^{-1}(v_i)\cap \sigma$ is a switch in $\sigma$. Thus 
there is a unique rectangle $A(R_i)$ in
$N_\delta-\beta[0,t+2a)$ which is mapped by $F$ to a trainpath
$\rho:[0,s]\to \tau$ with $\rho(0)=v_{i}$ and $\rho[0,1]=F(R_i)$.
Also let
$R_q$ be the rectangle $\hat R$ in $N_\delta$ which contains
$\beta(t+2a)$ in its interior, and let $A(R_q)$ 
be the rectangle in 
$N_\delta-\beta[0,t+2a)$ which is mapped by $F$ onto the arc
$F(\beta[t,t+2a])$. Finally, the rectangle $R_0$ is mapped by $A$ to the
rectangle which is contained in $R_q$ and whose collapse 
to a branch of $\sigma$ is large
at the switch which is the collapse of $\beta_0(t+2a)$.

Note that for $i\in 1\leq i\leq q-1$
the length
of $A(R_i)$ is not smaller than the length of $R_i$, 
and the length of $A(R_0)$ and of 
$A(R_q)$ is at least $a$. Hence for every rectangle
$\tilde R$ in $N_\delta$ 
different from $R_0,R_q$ the length
of $A(\tilde R)$ is not smaller than the length of $\tilde R$.
Now if the rectangle $R_0$ is short,
then a rectangle $A(\tilde R)$ can only be short
if $\tilde R$ is short and different from $R_0$. 
On the other hand,
if the length of the rectangle $R_0$ is at least $a$,
then at least one of the rectangles $R_{-1},R_{-2}$ is short.
Moreover, 
by inequality (\ref{increase}), the length of each of the rectangles
$A(R_i)$ for $i=-2,-1,0,q$ is at least $a$ and hence
a rectangle $A(\tilde R)$ can only be short if $\tilde R$
is short and different from both $R_{-1},R_{-2}$.
This shows the above claim.

We can repeat this construction with the train track
$\sigma$ and the foliated set $N_\delta-\beta[0,t+2a)$.
Namely, a short rectangle in $N_\delta-
\beta[0,t+2a)$ does not contain the point $\beta(t+2a)$
in its closure. Thus if $R^\prime$ is any short rectangle in
$N_\delta-\beta[0,t+2a)$, then its boundary contains a
corner of a complementary component of $S-N_\delta$ as before,
and this corner is the starting point of a one-sided infinite
simple geodesic $\beta^\prime$ which is contained in $N_\delta$ and
is disjoint from both $\lambda$ and $\beta$.
As a consequence, $\beta^\prime$ is contained in
$N_\delta-\beta[0,t+2a]$ and we can use our above construction
for $\beta^\prime$ and $N_\delta-\beta[0,t+2a)$ to reduce
the number of short rectangles of $N_\delta-\beta[0,t+2a)$.

In this way we construct
inductively in a uniformly bounded number of steps
a train track $\eta$ which is embedded in $N_\delta$,
which carries $\lambda$
and such that the length of the straightening of
each branch of $\eta$ is at
least $a$. Thus by the choice of $\delta$ and our
construction,  the train track $\eta$
is $a$-long and $\epsilon$-follows $\lambda$. This shows
the lemma.
\end{proof}

The next lemma shows how $C^1$-approximation of geodesic laminations
by train tracks relates to carrying.

\begin{lemma}\label{longcarry}
Let $\lambda$ be a geodesic lamination and 
let $\tau$ be a train track which carries $\lambda$.
Then there is a number $\epsilon >0$ such that $\tau$ carries each
train track $\sigma$
which $\epsilon$-follows $\lambda$.
\end{lemma}
\begin{proof}
Let $\tau$ be a train track which carries the
geodesic lamination $\lambda$. By Theorem 1.6.6 of \cite{PH92},
after possibly changing $\tau$ by an isotopy we may
assume that $\lambda$ is contained in the interior of a foliated
neighborhood $A$ of $\tau$ and is transverse to the ties. The
restriction of the differential $dF$ of the collapsing map $F:A\to
\tau$ to each tangent line of $\lambda$ is nonsingular. Since
$\lambda$ is a compact subset of $S$ and the projectivized tangent
bundle $PT\lambda$ of $\lambda$ is a compact subset of the
projectivized tangent bundle $PTS$ of $S$, there is a number
$\epsilon >0$ such that $A$ contains the $\epsilon$-neighborhood
of $\lambda$ and that moreover the restriction of $dF$ to each
line $z\in PTS$ which is contained in the $\epsilon$-neighborhood
of $PT\lambda$ is nonsingular.

The straightening
$\sigma$ of train track $\zeta$ which
$\epsilon/2$-follows $\lambda$ is embedded in $A$ and is transverse
to the ties.
The collapsing map $F$
restricts to a map on $\sigma$ which maps $\sigma$ to $\tau$.
Its differential maps each tangent line of $\sigma$ onto a
tangent line of $\tau$. If $\rho\subset \sigma$ is the straightening
of any trainpath on $\zeta$, then the exterior angles
at the breakpoints of $\rho$ are at most $\epsilon/2$. Thus
we can smoothen the graph $\sigma$ near its vertices to a train
track $\sigma^\prime$ which is isotopic to $\zeta$ and embedded in
$A$ and such that the restriction of $dF$ to each tangent
line of $\sigma^\prime$ is nonsingular. But this just means that
$\sigma^\prime$ is carried by $\tau$, with carrying map $F\vert \sigma^\prime$.
This shows the lemma.
\end{proof}

The space ${\cal P\cal M\cal L}$ of
\emph{projective measured geodesic laminations} on $S$ is 
the quotient of the space ${\cal M\cal L}$ of measured geodesic
laminations under the natural action of the multiplicative group
$(0,\infty)$. The space ${\cal
P\cal M\cal L}$ will be equipped with the quotient topology. With
this topology, ${\cal P\cal M\cal L}$ is homeomorphic to a sphere
\cite{FLP91}.
The mapping class group acts naturally on 
${\cal P\cal M\cal L}$ as a group of homeomorphisms.
We use this action to show.

\begin{lemma}\label{twocarry}
Let $\lambda,\mu$ be any two
complete geodesic laminations. Then there is a
complete train track
$\tau$ which carries both $\lambda,\mu$.
\end{lemma}
\begin{proof}
The mapping class group ${\cal M\cal C\cal G}$ acts on the
space of isotopy classes of complete train  tracks, and it acts as
a group of homeomorphisms on the space ${\cal P\cal M\cal L}$ of
projective measured geodesic laminations. Every
\emph{pseudo-Anosov} element $g\in {\cal M\cal C\cal G}$ admits a
pair of fixed points $\alpha_+\not=\alpha_-$ in ${\cal P\cal M\cal
L}$ and acts with respect to these fixed points with
\emph{sourth-sink dynamics}: For every neighborhood $V$ of
$\alpha_+$ and every neighborhood $W$ of $\alpha_-$ there is a
number $k>0$ such that $g^k({\cal P\cal M\cal L}-W)\subset V$ and
$g^{-k}({\cal P\cal M\cal L}-V)\subset W$.
Then $\alpha_+$ is called the \emph{attracting fixed point} of
$g$.

The fixed points of pseudo-Anosov elements in ${\cal M\cal C\cal G}$ are
measured geodesic laminations whose support is minimal and
\emph{fills up} $S$. This means that the complementary components
of this support are topological
discs or once punctured topological discs. As a consequence, every
complete geodesic lamination which contains the support of such a
fixed point as a sublamination consists of this support and
finitely many isolated leaves. In particular, there are only
finitely many complete geodesic laminations of this form.

For every complete train track $\tau$ on $S$ there is an open
subset $U$ of ${\cal P\cal M\cal L}$ with the property that the support of
each measured geodesic lamination $\lambda\in U$ is carried by
$\tau$ \cite{PH92}. The set $U$ corresponds precisely to those
projective measured geodesic laminations which are defined by
positive projective transverse measures on $\tau$ (Theorem 1.7.12
in \cite{PH92}).

Let $\lambda,\mu\in {\cal C\cal L}$ be any two complete
geodesic laminations.
Since fixed points of pseudo-Anosov elements in ${\cal M\cal C\cal G}$ are
dense in ${\cal P\cal M\cal L}$ there is a pseudo-Anosov element
$g\in {\cal M\cal C\cal G}$ with the following properties.
\begin{enumerate}
\item The attracting fixed point $\alpha_+\in {\cal P\cal M\cal L}$
of $g$ is contained
in the set $U$.
\item Every leaf of the support of
the repelling fixed point $\alpha_-$ of $g$ intersects
every leaf of $\lambda,\mu$ transversely.
\end{enumerate}

By Lemma \ref{compapp}, there is a complete geodesic lamination 
which contains the support of $\alpha_+$ as a sublamination.
We claim that each such complete 
geodesic lamination $\nu$ is carried by $\tau$. Namely,
choose a sequence of simple closed geodesics
$\{c_i\}$ which converge in the Hausdorff topology to $\nu$.
Then each of the curves $c_i$ defines a projective measured geodesic
lamination. After passing to a subsequence we may assume that
these projective measured geodesic laminations converge as $i\to
\infty$ in ${\cal P\cal M\cal L}$ to a projective measured
geodesic lamination $\zeta\in {\cal P\cal M\cal L}$ whose support
is necessarily a sublamination of $\nu$. Since the support of
$\alpha_+$ is the only minimal component of $\nu$, the support
of $\zeta$ coincides with the support of $\alpha_+$. Now the support
of $\alpha_+$ is \emph{uniquely ergodic} which means that there is
a single projective measured geodesic lamination whose support
equals the support of $\alpha_+$ \cite{FLP91}. Thus
$\zeta=\alpha_+$ and the projective measured 
geodesic laminations defined
by $c_i$ converge as $i\to\infty$ to $\alpha_+$. In particular,
for sufficiently large $i$ these projective measured geodesic
laminations
are contained in the open set $U$ and the curves $c_i$ are carried
by $\tau$. The set of geodesic laminations which are carried by
$\tau$ is closed in the Hausdorff topology and hence $\nu$ is
carried by $\tau$ as well. As a consequence, $\tau$ carries every
complete geodesic lamination which contains the support of
$\alpha_+$ as a sublamination.

Up to passing to a subsequence, 
as $k\to \infty$ the complete geodesic laminations
$g^k(\lambda),g^k(\mu)$ converge
in the Hausdorff topology to complete geodesic
laminations $\tilde \lambda,\tilde \mu$.
Moreover, for any projective transverse measure $\beta,\xi$
supported in $\lambda,\mu$ the projective 
measured geodesic laminations $g^k\beta,g^k\xi$ converge
as $k\to \infty$ to one of the
two fixed points $\alpha_+,\alpha_-$ for
the action of $g$ on ${\cal P\cal M\cal L}$
which moreover is supported in $\tilde \lambda,\tilde \mu$.
By (2) above, the geodesic laminations
$\lambda,\mu$ intersect the support of the lamination $\alpha_-$
transversely and therefore $g^k\beta\to \alpha^+,
g^k\xi\to \alpha^+$ as $k\to \infty$ \cite{FLP91}. 
As a consequence,
the laminations $\tilde \lambda,\tilde\mu$  ´
contain the support of $\alpha_+$ as a sublamination.
Thus by the above consideration, $\tilde
\lambda,\tilde\mu$ are carried by $\tau$. 
By Lemma \ref{recchar},
the set of all 
complete
geodesic laminations which are carried by $\tau$ 
is an open subset of the space
${\cal C\cal L}$ of all complete geodesic laminations and hence
there is some $k>0$ such that the laminations
$g^k\lambda,g^k\mu$ are both carried by $\tau$ as well. Then
$\lambda,\mu$ are carried by the train track $g^{-k}(\tau)$.
This completes the proof of the lemma. \end{proof}

Write $\sigma\prec\tau$ if the train track $\sigma$ 
is carried by the train track $\tau$. We have.

\begin{corollary}\label{threecarry}
For complete train tracks $\tau,\sigma$
on $S$ there are complete train tracks
$\tau_0,\sigma_0,\zeta$ such that $\tau_0\prec\tau,\tau_0\prec
\zeta$ and $\sigma_0\prec\sigma, \sigma_0\prec \zeta$.
\end{corollary}
\begin{proof}
For complete train tracks $\tau,\sigma$ on $S$ 
choose geodesic laminations $\lambda,\mu$
such that $\lambda$ is carried by $\tau$ and $\mu$ is carried by
$\sigma$. By Lemma \ref{twocarry} there is a complete train track $\zeta$
which carries both $\lambda$ and $\mu$.

By Lemma \ref{following}, 
for every $\epsilon >0$ there is a complete
train track $\tau_{\lambda,\epsilon},\sigma_{\mu,\epsilon}$ which
$\epsilon$-follows $\lambda,\mu$. By Lemma \ref{longcarry}, 
for sufficiently
small $\epsilon$ the train tracks $\tau_{\lambda,\epsilon},
\sigma_{\mu,\epsilon}$ are carried by $\zeta$ and moreover
$\tau_{\lambda,\epsilon}$ is carried by $\tau$, $\sigma_{\lambda,
\epsilon}$ is carried by $\sigma$.\end{proof}

A map $\rho$ which assigns to a positive integer 
$k$ contained in an interval $[m,n]\subset \mathbb{R}$
a complete
train track $\rho(k)$ such that $\rho(k+1)$ is obtained from
$\rho(k)$ by a single split will be called a \emph{splitting
sequence}. If there is a splitting sequence connecting a train
track $\tau$ to a train track $\sigma$ then we say that $\tau$ is
\emph{splittable} to $\sigma$. A map $\rho$ which assigns to each
integer $k$ from an interval $[m,n]\subset\mathbb{R}$ a complete train
track $\rho(k)$ such that either $\rho(k+1)$ is obtained from
$\rho(k)$ by a single split or a single collapse is called a
\emph{splitting and collapsing sequence}. The following corollary is a
consequence of Corollary \ref{threecarry} 
and Corollary 2.4.3 of \cite{PH92}.

\begin{corollary}\label{splittingarc}
Any two complete train tracks on $S$ can
be connected by a splitting and collapsing sequence.
\end{corollary}
\begin{proof}
By Corollary 2.4.3 of \cite{PH92}, if a complete
train track $\tau$ is carried by a complete train track $\sigma$
then there is a complete train track $\rho$ which can be obtained
from both $\tau$ and $\sigma$ by a splitting sequence. Namely,
since $\tau$ is birecurrent and maximal it carries a geodesic
lamination which is both minimal and maximal (see the discussion
in the proof of Lemma \ref{recchar}). This geodesic lamination is the
support of a transverse measure and hence it defines a measured
geodesic lamination $\mu$. Corollary 2.4.3 of \cite{PH92} shows
that $\sigma$ and $\tau$ can be split to the same generic
birecurrent
train track $\eta$ which carries the measured geodesic lamination
$\mu$ (where a splitting move in the sense of Penner and Harer may
be a \emph{collision} which is defined to be a 
split followed by the removal of the diagonal of the split).
Since the support of the measured geodesic lamination $\mu$ is a
maximal geodesic lamination, the train track $\eta$ is maximal and
hence complete (and no collision can have occurred in the process).

Let $\tau,\sigma$ be any complete train tracks on $S$, and let
$\tau_0,\sigma_0,\zeta$ be as in Corollary \ref{threecarry}. Then the train
tracks $\tau,\tau_0$ are both splittable to the same complete
train track $\tau_1$, and the train tracks $\sigma,\sigma_0$ are
both splittable to the same complete train track $\sigma_1$. In
particular, the train tracks $\tau_1,\sigma_1$ are carried by
$\zeta$. Using Corollary 2.4.3 of \cite{PH92} once more we deduce
that the train tracks $\tau_1,\zeta$ are splittable to the same
train track $\tau_2$, and the train tracks $\sigma_1,\zeta$ are
splittable to the same train track $\sigma_2$. But this just means
that there is a splitting and collapsing sequence from $\tau$ to $\sigma$
which passes through $\tau_1,\tau_2,\zeta, \sigma_2,\sigma_1$.
\end{proof}

The following
corollary is now immediate from Corollary \ref{splittingarc}.

\begin{corollary}\label{connectedness} The train track complex
${\cal T\cal T}$ is connected.\end{corollary}

As a consequence, if we identify each edge in ${\cal T\cal T}$
with the unit interval $[0,1]$ then this provides
${\cal T\cal T}$ with the structure of a connected locally finite metric
graph. Thus ${\cal T\cal T}$ is a locally compact complete geodesic
metric space. In the sequel we always assume that ${\cal T\cal T}$
is equipped with this metric without further comment.

\section{Train tracks and the  mapping class group}

The purpose of this section is to show that the connected
graph ${\cal T\cal T}$ which we defined in Section 2 is quasi-isometric
to the mapping class group of the surface $S$ of
genus $g$ with $m$ punctures where $3g-3+m\geq 2$.
A quasi-isometry between two metric spaces was defined in the
introduction.

Consider first the \emph{Teichm\"uller space}
${\cal T}_{g,k}$ of complete marked hyperbolic
metrics on $S$ of finite volume.
Let $\epsilon >0$ be smaller than half of a
\emph{Margulis constant} for
hyperbolic surfaces (see \cite{B92}).
Call a surface $S_0\in {\cal T}_{g,k}$
\emph{thick}
if the systole of $S_0$, i.e. the length of the shortest closed
geodesic, is at least $\epsilon$. For sufficiently small
$\epsilon$ the set ${\cal T}_{{\rm thick}}$ of thick surfaces is a
connected closed subset of ${\cal T}_{g,k}$
with dense interior which
is invariant under the natural action of the mapping class group
${\cal M\cal C\cal G}$. The quotient of ${\cal T}_{{\rm thick}}$ 
under the action of
${\cal M\cal C\cal G}$ is compact.

The \emph{Teichm\"uller metric} on ${\cal T}_{g,k}$
is a complete ${\cal M\cal C\cal G}$-invariant
Finsler metric $\Vert\,\Vert$. We
equip ${\cal T}_{{\rm thick}}$ with the path metric defined by the
restriction of this Finsler metric. In other words, the distance
between two points in ${\cal T}_{{\rm thick}}$ 
is defined to be the infimum
of the $\Vert\,\Vert$-lengths of paths in ${\cal T}_{{\rm thick}}$
connecting these two points. With respect to this distance
function, the mapping class group ${\cal M\cal C\cal G}$ acts
properly discontinuously, isometrically and cocompactly
on ${\cal T}_{{\rm thick}}$.

The mapping class group ${\cal M\cal C\cal G}$ is finitely generated
(see \cite{I02}). A 
finite symmetric set ${\cal G}$
of generators defines a word norm and hence an 
${\cal M\cal C\cal G}$-invariant distance on ${\cal M\cal C\cal G}$.
We always assume that ${\cal M\cal C\cal G}$ is equipped with
such a fixed distance. The
next easy lemma is included here as an illustration of the
various ways to understand the geometry of ${\cal M\cal C\cal G}$.

\begin{lemma}\label{quasiiso}
${\cal M\cal C\cal G}$ and
${\cal T}_{{\rm thick}}$ are equivariantly quasi-isometric.
\end{lemma}
\begin{proof}
The mapping class group ${\cal M\cal C\cal G}$ acts on the
length space ${\cal T}_{{\rm thick}}$ isometrically, properly and
cocompactly. By the well known lemma of \v{S}varc-Milnor
(Proposition I.8.19 in \cite{BH99})
this implies that ${\cal T}_{{\rm thick}}$ is
equivariantly quasi-isometric to ${\cal M\cal C\cal G}$.
\end{proof}

The mapping class group also acts naturally
as a group of simplicial
isometries on
the train track complex ${\cal T\cal T}$.
Note that by definition, the train track complex is a
locally finite directed metric graph and hence a complete
locally compact geodesic metric space.
We want to show that ${\cal T\cal T}$ is equivariantly quasi-isometric
to ${\cal M\cal C\cal G}$. By the \v{S}varc-Milnor lemma,
for this it is enough to show that the action of ${\cal M\cal C\cal G}$ on
${\cal T\cal T}$ is proper and cocompact.
We show first that this action is cocompact.
Recall that the vertices of ${\cal T\cal T}$ consist of \emph{generic}
train tracks.

\begin{lemma}\label{cocompact}
${\cal M\cal C\cal G}$ acts cocompactly on ${\cal T\cal T}$.
\end{lemma}
\begin{proof}
The number $k$ of switches of a complete train
track $\tau$
on $S$ is just the number of cusps of the complementary components
of $\tau$, and this number
only depends on the
topological type of $S$. There are only finitely many abstract
trivalent graphs with $k$ vertices. We assign to each
half-edge of such a graph one of the three colors red, yellow, green
in such a way that every vertex is incident on a half-edge of each
color. We do not require that the two half-edges of a single edge
have the same color. We call a trivalent graph with such a
coloring a \emph{colored graph}. Clearly there are only finitely
many colored graphs with $k$ vertices up to isomorphism preserving the
coloring.

To every train track $\tau\in 
{\cal V}({\cal T\cal T})$ we associate
a colored trivalent graph $G$ with $k$ vertices as follows.
The underlying topological graph of $G$ is just the abstract
graph defined by $\tau$.
Thus the half-branches of $\tau$ are in one-one correspondence
to the half-edges of $G$. We color a half-edge of $G$ red
if and only if the corresponding half-branch in $\tau$ is
large.
We orient such a large half-branch $b$ in $\tau$ in
such a way that it ends at the switch $v$ of $\tau$ on which it is
incident.
A neighborhood of $b$ in $S$ is divided by $b$ into two components.
The orientation of $b$ together with
the orientation of $S$ determine the component
to the right and to the left of $b$.
One of
the two small half-branches incident on $v$ lies to the left
of the other, i.e. the union of this half-branch with the
half-branch $b$ is contained in the boundary of the component of
$S-\tau$ to the left of $b$.
This half-branch will be colored yellow, and the
second small half-branch incident on $v$
(which lies to the right of the yellow
half-branch)
will be colored green. In this way we obtain a map $\Psi$
from ${\cal V}({\cal T\cal T})$ to a finite set of colored trivalent
graphs with $k$ vertices.

Let $\tau,\sigma\in {\cal V}({\cal T\cal T})$ be
such that $\Psi(\tau)=\Psi(\sigma)$. Then the map $\Psi$
determines a homeomorphism of $\tau$ onto $\sigma$
which preserves the coloring of the half-branches. We claim
that this homeomorphism maps the sides of the complementary
components of $\tau$ to the sides of the complementary components
of $\sigma$.

For this let $b$ be any branch of $\tau$ and
let $v$ be a switch of $\tau$ on which $b$ is incident.
The branch $b$ is contained in the boundary
of two (not necessarily distinct) complementary components
of $\tau$.
We orient $b$ in such a way that it ends at $v$; then
we can distinguish the complementary
component $T$ of $\tau$ which is to the left of $b$.
We have to show that the finite trainpath $\rho$ on $\tau$
which defines the
side of $T$ containing $b$
is uniquely determined by the colored graph
$\Psi(\tau)=G$.

Assume first
that the half-branch of $b$ incident on $v$ is
colored red. By the
definition of our coloring, the branch $b$ is large at $v$
and $\rho$ necessarily
contains the yellow half-branch incident on $v$.
Similarly, if the half-branch of $b$ incident on $v$ is
green (and hence small) then $\rho$ contains the red (large)
half-branch incident on $v$. If the half-branch of $b$ incident
on $v$ is yellow then $T$ has a cusp at $v$ and the trainpath $\rho$
ends at $v$. But this just means that we can successively
construct the trainpath $\rho$ from the
coloring. In other words, the sides of the complementary
components of $\tau$ are uniquely determined by
the colored graph $G$. As a consequence, whenever $\Psi(\tau)=
\Psi(\sigma)$ then the homeomorphism of $\tau$ onto $\sigma$ induced
by $\Psi$ maps the boundary of each complementary trigon of $\tau$
to the boundary of a complementary trigon of $\sigma$, and
it maps the boundary of a complementary once punctured monogon
of $\tau$ to the boundary of a complementary once punctured
monogon of $\sigma$.
Thus this homeomorphism of $\tau$ onto $\sigma$
can be extended to
the complementary components of $\tau$, and this extension is
an orientation preserving  
homeomorphism of $S$ which maps $\tau$ to $\sigma$.
In other words, $\sigma$ and $\tau$ are contained in
the same orbit of the action of ${\cal M\cal C\cal G}$ on ${\cal T\cal T}$.

As a consequence,
there is a \emph{finite} subset $A$ of
${\cal V}({\cal T\cal T})$ with the property that the
translates of $A$ under ${\cal M\cal C\cal G}$ cover
all of ${\cal V}({\cal T\cal T})$. But this is equivalent to saying
that the action of ${\cal M\cal C\cal G}$ on ${\cal T\cal T}$ is cocompact.
\end{proof}

\begin{lemma}\label{proper}
The action of ${\cal M\cal C\cal G}$ on ${\cal T\cal T}$ is
proper.
\end{lemma}
\begin{proof}
Since ${\cal T\cal T}$ is a locally finite graph, every
compact subset of ${\cal T\cal T}$ contains only finitely many edges.
The action of ${\cal M\cal C\cal G}$ on ${\cal T\cal T}$ is
simplicial and isometric and therefore this action is proper
if and only if the stabilizer in ${\cal M\cal C\cal G}$ of every vertex
$\tau\in {\cal V}({\cal T\cal T})$ is finite (compare \cite{BH99}).

To show that this is the case, let 
$\tau\in {\cal V}({\cal T\cal T})$ and let
$G<{\cal M\cal C\cal G}$ be the stabilizer of $\tau$. Since
the number of branches and switches of $\tau$
only depends on
the topological type of $S$,
the subgroup $G_0$ of $G$ which fixes
every branch and every switch of $\tau$ is of
finite index in $G$.
But a homeomorphism of $S$ which preserves
each branch and switch of $\tau$ is isotopic to a map which
preserves the branches pointwise. Such a map then fixes pointwise
the boundaries of the complementary regions of $\tau$. However,
every complementary region of $\tau$ is a topological disc or a once
punctured topological disc, and every homeomorphism of such a disc
which fixes the boundary pointwise is isotopic to the identity
(see \cite{I02}). This shows that $G_0$ is trivial and therefore the
stabilizer of $\tau\in {\cal V}({\cal T\cal T})$ under the action of 
${\cal M\cal C\cal G}$ is finite.
\end{proof}

\begin{corollary}\label{quasisometric}
${\cal M\cal C\cal G}$ and ${\cal T\cal T}$ are
quasi-isometric.
\end{corollary}
\begin{proof}
By Lemma \ref{cocompact} and 
Lemma \ref{proper}, the mapping class group
${\cal M\cal C\cal G}$ acts properly and cocompactly on
the length space ${\cal T\cal T}$. Thus
the corollary
follows from the lemma of \v{S}varc-Milnor
(Proposition I.8.19 in \cite{BH99}). \end{proof}

\section{Flat cones}

In this section we define a family of connected
subgraphs of ${\cal T\cal T}$, one for every 
complete train track $\tau$ and every complete geodesic
lamination $\lambda$ carried by $\tau$. We show that these
subspaces equipped with their \emph{intrinsic} path-metric
are isometric to cubical graphs contained in an euclidean space of
fixed dimension. Thus these subgraphs can be viewed
as ``flat cones''. In Section 6 we will see that these flat
cones are quasi-isometrically embedded in ${\cal T\cal T}$.
This then immediately implies Theorem \ref{theorem3} and
Corollary \ref{corollary3} from the introduction.

For a complete train
track $\tau\in {\cal V}({\cal T\cal T})$ and a complete geodesic
lamination $\lambda\in {\cal C\cal L}$ which is carried
by $\tau$, define the
\emph{flat cone} $E(\tau,\lambda)\subset {\cal T\cal T}$
to be the full subgraph of ${\cal T\cal T}$ whose vertices
consist of all complete train tracks which can
be obtained from $\tau$ by any $\lambda$-splitting sequence,
i.e. by a sequence of $\lambda$-splits.
By construction, $E(\tau,\lambda)$ is a connected
subgraph of ${\cal T\cal T}$ and hence it can
be equipped with
an intrinsic path metric $d_E$.

Let $m>0$ be the number of branches of
a complete train track on $S$. Let $e_1,\dots,e_m$ be the
standard basis of $\mathbb{R}^m$.
Define a \emph{cubical graph}
in $\mathbb{R}^m$ to be an embedded graph whose vertices
are points with integer coordinates (i.e. points contained
in $\mathbb{Z}^m$) and whose edges are line segments of length
one connecting two of these vertices $v_1,v_2$
with $v_1-v_2=\pm e_i$ for some $i\leq m$. Note that
each such edge has a natural direction.

As in the introduction,
define a \emph{splitting arc} in ${\cal T\cal T}$ to be
a map $\gamma:[0,n]\to {\cal T\cal T}$ with the property
that for any integer $i\in [1,n]$, the arc
$\gamma[i-1,i]$ is a directed edge in ${\cal T\cal T}$.
In other words, $\{\gamma(i)\}$ is a splitting sequence.
We have.

\begin{lemma}\label{flatstrip}
For every flat cone $E(\tau,\lambda)\subset
{\cal T\cal T}$ there is an isometry $\Phi$ of $E(\tau,\lambda)$
equipped with the intrinsic path metric $d_E$ 
onto a cubical graph in $\mathbb{R}^m$ which maps any splitting
arc in $E(\tau,\lambda)$ to a directed edge-path in
$\Phi(E(\tau,\lambda))$. 
\end{lemma}
\begin{proof}
Let $\lambda\in {\cal C\cal L}$, let 
$\tau\in {\cal V}({\cal T\cal T})$
be a complete train track which carries $\lambda$, let
$\sigma\in E(\tau,\lambda)$
and let $e$ be a large branch in $\tau$.
We claim that whether or not a splitting
sequence connecting $\tau$ to $\sigma$  contains a split
at the large branch $e$ only depends on $\sigma$
but not on the choice of the splitting sequence.

For this note first that
$\sigma$ is carried by $\tau$.
For a large branch
$e$ of $\tau$ define
$\nu(e,\sigma)$ to be the minimal cardinality
of the preimage of a point
$x$ in the interior of $e$
under \emph{any} carrying map $F:\sigma\to\tau$.
Let $\tau^\prime$ be obtained from $\tau$ by a
$\lambda$-split at a
large branch $e^\prime\not=e$ and assume that $\tau^\prime$ is
splittable to $\sigma$. Then
$\tau^\prime$ carries $\sigma$, and the branch
in $\tau^\prime$ correponding to $e$ under the natural
identification of the branches of $\tau$ with
the branches of $\tau^\prime$ is large. We denote it
again by $e$. Let $\nu^\prime(e,\sigma)$ be the
minimal cardinality of the preimage
of a point $x$ in the interior of $e$ under
any carrying map
$\sigma\to\tau^\prime$; we claim that
$\nu^\prime(e,\sigma)\geq \nu(e,\sigma)$.

To see this, simply observe that there are
disjoint neighborhoods $U^\prime$ of $e^\prime$ and $U$ of $e$ in
$S$ and there is a carrying
map $G:\tau^\prime\to\tau$ which equals the identity outside
$U^\prime$.
Every carrying map $F^\prime:\sigma\to\tau^\prime$ can be
composed with $G$ to a carrying map
$G\circ F^\prime:\sigma\to\tau$.
For a point $x$ in the interior of $e$
the cardinality of the preimage of $x$ under
the carrying map $G\circ F^\prime:\sigma\to\tau$ coincides with
the cardinality of the preimage of $x=G^{-1}(x)$
under the carrying map $F^\prime:\sigma\to\tau^\prime$.
Thus we have $\nu(e,\sigma)\leq \nu^\prime(e,\sigma)$.

As a consequence, if $e\subset \tau$ is a large branch with
$\nu(e,\sigma)\geq 2$, then every splitting sequence connecting $\tau$ to
$\sigma$ has to contain a split at $e$, and the choice of a right
or left split is determined by the requirement that the split
track carries $\lambda$. On the other hand, if 
$\nu(e,\sigma)=1$ then by
Lemma \ref{foliated} from the appendix, 
$\sigma$ is not carried by a split at $e$. Thus whether
or not a splitting sequence connecting $\tau$ to $\sigma$ contains
a split at $e$ is independent of the splitting sequence.

As above, let
$m$ be the number of branches of a complete train track on $S$.
Denote by $e_1,\dots,e_m$ the standard basis
of $\mathbb{R}^m$ and 
choose any point $q\in \mathbb{Z}^m$.
Number the branches of $\tau$ in an arbitrary way. Note that
this numbering induces a natural numbering of the branches
on any train track which can be obtained from $\tau$ by
a single split. Let
$\alpha:[0,\infty)\to E(\tau,\lambda)$ be 
any splitting arc with $\alpha(0)=\tau$.
Define a map $\Phi_\alpha:\alpha[0,\infty)\to \mathbb{R}^m$
inductively as follows. 
Let $\Phi_\alpha(\tau)=q$ and assume by induction that
$\Phi_\alpha$ has been defined on $\alpha[0,\ell-1]$
for some $\ell\geq 1$.
Let $p\leq m$ be the number 
of the large branch $e$ of $\alpha(\ell-1)$ induced from
the numbering of the branches of $\tau$ via $\alpha$ so that
$\alpha(\ell)$ is obtained
from $\alpha(\ell -1)$ by a single
$\lambda$-split at $e$. Define
$\Phi_\alpha(\alpha[\ell-1,\ell])$ to be the line segment
in $\mathbb{R}^m$ connecting $\Phi_\alpha(\alpha(\ell-1))$ to
$\Phi_\alpha(\alpha(\ell-1))+e_p$.
In this way we obtain for every splitting arc 
$\alpha:[0,\infty)\to
E(\tau,\lambda)$ a map $\Phi_\alpha:\alpha[0,\infty)\to
\mathbb{R}^m$.

We claim that for every train track 
$\tau\in {\cal V}({\cal T\cal T})$, for every complete
geodesic lamination $\lambda\in {\cal C\cal L}$ carried
by $\tau$, for every
$\sigma\in E(\tau,\lambda)$ and for every splitting arc
$\alpha:[0,\infty)\to E(\tau,\lambda)$ issuing from 
$\alpha(0)=\tau$ and
passing through
$\alpha(\ell)=\sigma$, the image of $\sigma$ under the map
$\Phi_\alpha$ which is determined as above by 
a numbering of the branches of $\tau$, a point $q\in
\mathbb{Z}^m$ and by
$\alpha$ is in fact independent of 
the splitting arc $\alpha$ connecting $\tau$ to $\sigma$.

For this we proceed by
induction on the length of the \emph{shortest} splitting arc
connecting $\tau$ to $\sigma$. The case that this length vanishes
is trivial, so assume that for some $k\geq 1$ the above claim
holds for all $\tau,\lambda$ and all
$\sigma\in E(\tau,\lambda)$ which can be obtained
from $\tau$ by a splitting arc of length at most
$k-1$. Let $\sigma\in E(\tau,\lambda)$ be such that
there is a splitting arc $\alpha:[0,k]\to 
E(\tau,\lambda)$ of length $k$ connecting $\alpha(0)=\tau$ to $\sigma$
and let $\beta:[0,p]\to E(\tau,\lambda)$ be a splitting arc 
connecting $\beta(0)=\tau$ to $\sigma$ of length $p\geq k$.
For a fixed choice of a numbering of the branches
of $\tau$ and a fixed point $q\in\mathbb{Z}^m$, these
splitting arcs determine maps 
$\Phi_\alpha:\alpha[0,k]\to \mathbb{R}^m,
\Phi_\beta:\beta[0,p]\to\mathbb{R}^m$
with $\Phi_\alpha(\tau)=\Phi_\beta(\tau)=q$.

Let $b\subset \tau$ be
the large branch with the property that $\alpha(1)$ is obtained
from $\tau$ by a split at $b$. Let $s\leq n$ be the number of $b$
with respect to our numbering of the branches of $\tau$; then
$\Phi_\alpha(\alpha(1))=q+e_s$. By the discussion in the
beginning of this proof, the
splitting sequence $\{\beta(j)\}_{0\leq j\leq p}$ also contains a
split at the branch $b$. If this is the $j$-th split in this
splitting sequence, then the first $j-1$ splits of the sequence
commute with the split at $b$. Let $\beta^\prime:[0,p]\to
E(\tau,\lambda)$ be the splitting arc 
which we
obtain from $\beta$ by exchanging the orders of the 
first $j$ splits in such
a way that $\beta^\prime(0)=\tau$ 
and that for
every $0\leq i\leq j-1$ the train track $\beta^\prime(i+1)$ is
obtained from $\beta(i)$ by a $\lambda$-split at $b$.
This splitting arc
then determines a map $\Phi_{\beta^\prime}:\beta^\prime[0,p]\to
\mathbb{R}^m$.
By construction of the maps
$\Phi_\beta,\Phi_{\beta^\prime}$ we have
$\Phi_\beta(\beta(j))=\Phi_{\beta^\prime}(\beta^\prime(j))$ and
$\Phi_\beta(\sigma)=\Phi_{\beta^\prime}(\sigma)$, moreover
$\alpha(1)=\beta^\prime(1)$ and hence $\Phi_\alpha(\alpha(1))=
\Phi_{\beta^\prime}(\beta^\prime(1))$. Therefore we can apply the
induction hypothesis to the splitting arcs $\alpha[1,k]$ and
$\beta^\prime[1,p]$ issuing from $\alpha(1)$,
the numbering of the branches
of $\alpha(1)$ inherited from the numbering of the branches
of $\tau$ and the point $\Phi_\alpha(\alpha(1))=q+e_s
\in \mathbb{Z}^m$ 
to conclude that the images of $\sigma$ under
the maps $\Phi_\alpha$ and $\Phi_\beta$ coincide. 

By induction,
this construction defines a path-isometric map $\Phi$ of
$E(\tau,\lambda)$ into $\mathbb{R}^m$ whose image is a cubical
graph in $\mathbb{R}^m$. The above discussion shows that
this map is uniquely determined by the
choice of a numbering of the branches of $\tau$ and the choice of
$\Phi(\tau)\in \mathbb{Z}^m$. This shows the lemma. \end{proof}

Since directed edge-paths in a cubical graph in $\mathbb{R}^m$ 
are geodesics,
we obtain as an immediate corollary.

\begin{corollary}\label{splitine}
Splitting arcs are geodesics in
$(E(\tau,\lambda),d_E)$.
\end{corollary}

In the remainder of this section we describe the intrinsic
geometry of the flat cones $E(\tau,\lambda)$ more explicitly.

The \emph{Hausdorff distance} between two subsets
$A,B$ of a metric space $(X,d)$ 
is the infimum of all numbers $r>0$ such that
$A$ is contained in the $r$-neighborhood of $B$ and 
$B$ is contained in the $r$-neighborhood of $A$. If the diameter
of $A,B$ is infinite then 
the Hausdorff distance between $A$ and $B$ 
may be infinite.

A connected subspace
$Y$ of a geodesic metric space $(X,d)$ is called
\emph{strictly convex} if for any two points
$y,z\in Y$, every geodesic in $(X,d)$ 
connecting $y$ to $z$ is entirely
contained in $Y$. The next lemma
is a first easy step toward an understanding of the intrinsic
geometry of a flat cone. For this note that 
for every vertex $\sigma\in E(\tau,\lambda)$ the
flat cone $E(\sigma,\lambda)$ is a complete subgraph
of $E(\tau,\lambda)$.

\begin{lemma}\label{convex}
For $\sigma\in E(\tau,\lambda)$, the subspace
$E(\sigma,\lambda)$ of $(E(\tau,\lambda),d_E)$ is strictly convex,
and its Hausdorff distance to $E(\tau,\lambda)$ does not exceed 
$d_E(\tau,\sigma)$.
\end{lemma}
\begin{proof}
Since a strictly convex subspace $A$ of a strictly convex
subspace $B$ of a geodesic metric space $X$ is strictly convex
in $X$, it suffices to show the following.
If $\sigma\in E(\tau,\lambda)$ is obtained from 
$\tau$ by a single split at a large branch $e$ then
$E(\sigma,\lambda)$ is a strictly convex subspace of 
$E(\tau,\lambda)$ whose Hausdorff distance to 
$E(\tau,\lambda)$ equals one.

For this note that by Lemma \ref{flatstrip} and its proof,  
if $\sigma\in E(\tau,\lambda)$ is a train track which
can be obtained from $\tau$ by a single split at a large branch $e$
and if $\eta\in E(\tau,\lambda)$ is \emph{not} contained in 
$E(\sigma,\lambda)$, then a splitting sequence connecting $\tau$ to 
$\eta$ does not contain a split at $e$. Moreover, $\sigma$ is
splittable to a train track $\eta^\prime\in E(\tau,\lambda)$ which
can be obtained from $\eta$ by a single split at $e$. In other words,
there is a natural retraction 
$R:E(\tau,\lambda)\to E(\sigma,\lambda)$
which equals the identity on $E(\tau,\sigma)$ and maps a train track
$\eta\in E(\tau,\lambda)-E(\sigma,\lambda)$ to the 
train track obtained from $\eta$ by a $\lambda$-split at $e$.
This shows that the Hausdorff distance between
$E(\tau,\lambda)$ and $E(\sigma,\lambda)$ does not exceed
$1$. 

Now let $\zeta_1,\zeta_2$ be any vertices in $E(\tau,\lambda)$ which
are connected by an edge. We may assume that $\zeta_2$ can
be obtained from $\zeta_1$ by a single split at a large branch
$e^\prime$. If both $\zeta_1,\zeta_2$ are contained in 
$E(\tau,\lambda)-E(\sigma,\lambda)$ then since $\lambda$-splits
at distinct large branches commute, the train tracks
$R(\zeta_1),R(\zeta_2)$ are connected by an edge in 
$E(\sigma,\lambda)$. On the other hand, if $\zeta_1\in
E(\tau,\lambda)-E(\sigma,\lambda)$ and 
$\zeta_2\in E(\sigma,\lambda)$ then $R(\zeta_1)=R(\zeta_2)=\zeta_2$.
As a consequence, the retraction $R$ is distance non-increasing.
Moreover, any simplicial path in $E(\tau,\lambda)$ 
connecting two points 
$\sigma,\eta\in E(\sigma,\lambda)$ and which is not entirely
contained in $E(\sigma,\lambda)$ passes through an edge which
is mapped by $R$ to a single point. 
This implies strict convexity of $E(\sigma,\lambda)\subset
E(\tau,\lambda)$.
\end{proof}

The next lemma can be used to calculate distances
in a flat cone $(E(\tau,\lambda),d_E)$ explicitly.

\begin{lemma}\label{distinflat}
Let $\sigma,\eta\in E(\tau,\lambda)$ be any two vertices.
Then there are unique train tracks $\Theta_-(\sigma,\eta),
\Theta_+(\sigma,\eta)\in E(\tau,\lambda)$ 
with the following properties.
\begin{enumerate}
\item $\sigma,\eta\in E(\Theta_-(\sigma,\eta),\lambda)
\subset E(\tau,\lambda)$, and
there is a geodesic in $(E(\tau,\lambda),d_E)$ connecting
$\sigma$ to $\eta$ which passes through $\Theta_-(\sigma,\eta)$.
\item $E(\Theta_+(\sigma,\eta),\lambda)=
E(\sigma,\lambda)\cap E(\eta,\lambda)$, and
there is a geodesic in 
$(E(\tau,\lambda),d_E)$ connecting $\sigma$ to 
$\eta$ which passes
through $\Theta_+(\sigma,\eta)$.
\item $d_E(\sigma,\Theta_-(\sigma,\eta))=
d_E(\eta,\Theta_+(\sigma,\eta))$.
\end{enumerate}
\end{lemma}
\begin{proof}
Let $\sigma,\eta$ be vertices in $E(\tau,\lambda)$.
Then $\sigma,\eta$ are complete train tracks, and 
$\tau$ is splittable
to both $\sigma$ and $\eta$. 
Let $A\subset E(\tau,\lambda)$ be the set of all complete train tracks
which can be obtained from $\tau$ by a splitting sequence
and which are splittable to both $\sigma$ and $\eta$. Note that
$A$ is a \emph{finite} set of vertices of $E(\tau,\lambda)$.
For $\beta,\beta^\prime\in A$ write $\beta<\beta^\prime$
if $\beta$ is splittable to $\beta^\prime$. Then $<$ is a partial order
on $A$. 

Let $\Theta_-(\sigma,\eta)$ 
be a maximal element for this partial order.
Then $\sigma,\eta$ are both contained in
$E(\Theta_-(\sigma,\eta),\lambda)$ and
hence by Lemma \ref{convex}, every geodesic in $E(\tau,\lambda)$ 
connecting $\sigma$ to $\eta$ is contained in 
$E(\Theta_-(\sigma,\eta),\lambda)$.
If $e$ is any large branch of $\Theta_-(\sigma,\eta)$ 
and if the train track
obtained from $\Theta_-(\sigma,\eta)$ by a 
$\lambda$-split at $e$ is splittable to 
$\sigma$ then by maximality of $\Theta_-(\sigma,\eta)$, 
it is not splittable to $\eta$.

Let $\Phi:E(\Theta_-(\sigma,\eta),\lambda)\to 
\mathbb{R}^m$ be an isometry
of $E(\Theta_-(\sigma,\eta),\lambda)$ 
onto a cubical graph in $\mathbb{R}^m$
defined as in the proof of Lemma \ref{flatstrip} by the 
choice of the basepoint $\Phi(\Theta_-(\sigma,\eta))=0$ 
and a numbering of the branches of 
$\Theta_-(\sigma,\eta)$.

We claim that up to a permutation of the standard basis of
$\mathbb{R}^m$, there is a number $\ell \geq 1$ such that for the
standard direct orthogonal decomposition
$\mathbb{R}^m=\mathbb{R}^\ell\oplus \mathbb{R}^{m-\ell}$ we have
$\Phi(\sigma)\in \mathbb{R}^\ell\times \{0\}$ and $\Phi(\eta)\in
\{0\}\times \mathbb{R}^{m-\ell}$. 
Namely, by the choice of the train track
$\Theta_-(\sigma,\eta)$ 
and the fact that $\sigma,\eta$ both carry the complete
geodesic lamination $\lambda$, the set of large branches of 
$\Theta_-(\sigma,\eta)$ can be partitioned into disjoint subsets
${\cal E}^+,{\cal E}^-$ such that a splitting sequence connecting
$\Theta_-(\sigma,\eta)$ 
to $\sigma$ does not contain any split at a large branch
branch $e\in {\cal E}^+$ and that a splitting sequence connecting
$\Theta_-(\sigma,\eta)$ 
to $\eta$ does not contain any split at a large branch
$e\in {\cal E}^-$.

Following \cite{PH92}, we call a trainpath 
$\rho:[0,p]\to \Theta_-(\sigma,\eta)$
\emph{one-sided large} if for every $i<p$ the half-branch
$\rho[i,i+1/2]$ is large and if $\rho[p-1,p]$ is a large branch. A
one-sided large trainpath 
$\rho:[0,p]\to \Theta_-(\sigma,\eta)$ is embedded
\cite{PH92}, and for every $i\in \{1,\dots,p-1\}$ the branch
$\rho[i-1,i]\subset \Theta_-(\sigma,\eta)$ 
is mixed. For every large half-branch $\hat
b$ of $\Theta_-(\sigma,\eta)$ there is 
a unique one-sided large trainpath issuing
from $\hat b$. Define ${\cal A}_0^+,{\cal A}_0^-$ to be the set of
all branches of $\Theta_-(\sigma,\eta)$ 
contained in a one-sided large trainpath
ending at a branch in ${\cal E}^{+},{\cal E}^-$. Then the sets
${\cal A}_0^+,{\cal A}_0^-$ are disjoint, and a branch of 
$\Theta_-(\sigma,\eta)$
is \emph{not} contained in ${\cal A}_0^+\cup {\cal A}_0^-$ if and
only if it is small. Each endpoint of a small branch is a starting
point of a one-sided large trainpath. Define ${\cal A}^{\pm}$ to
be the union of ${\cal A}_0^{\pm}$ with all small branches $b$ of
$\Theta_-(\sigma,\eta)$ 
with the property that both large half-branches incident
on the endpoints of $b$ are contained in ${\cal A}_0^{\pm}$. If
$b\not\in{\cal A}^+ \cup {\cal A}^-$ then $b$ is a small branch
incident on two distinct switches, and one of these switches is
the starting point of a one-sided large trainpath in ${\cal
A}_0^+$,the other is the starting point of a one-sided large
trainpath in ${\cal A}_0^-$.

The map $\Phi$ is determined
by a numbering of the branches of $\Theta_-(\sigma,\eta)$.
We may assume that this numbering is such that for the cardinality
$\ell$ of ${\cal A}^-$, the set ${\cal A}^-$ consists of
the branches with numbers $1,\dots,\ell$.
A splitting sequence connecting $\Theta_-(\sigma,\eta)$ to
$\sigma$ does not contain any split at a large branch
$e\in {\cal E}^+$ by assumption.
Therefore, such a splitting sequence only
contains splits at the branches in ${\cal A}^-$.
By the choice
of our numbering, the image of any such splitting
sequence under the map
$\Phi$ is contained in the linear subspace spanned
by the first $\ell$ vectors of the standard basis of
$\mathbb{R}^m$.
Similarly, the image under $\Phi$ of a
splitting sequence connecting $\Theta_-(\sigma,\eta)$ 
to $\eta$ is contained
in the subspace
$\mathbb{R}^{m-\ell}\subset \mathbb{R}^m$ spanned
by the last $m-\ell$ vectors of the standard basis. This shows our claim.

The image under $\Phi$ of any splitting arc in 
$E(\Theta_-(\sigma,\eta),\lambda)$
is a geodesic edge path in the standard cubical graph
${\cal G}\subset \mathbb{R}^m$ 
with vertex set $\mathbb{Z}^m$ and where two vertices are
connected by a straight line segment if their distance
in $\mathbb{R}^m$ equals one.
If $\gamma_\sigma,\gamma_\eta$
are such geodesic edge-paths connecting $0=\Phi(\Theta_-(\sigma,\eta))$ to
$\Phi(\sigma),\Phi(\eta)$ which are images of 
splitting arcs then
$\gamma_\sigma\subset \mathbb{R}^\ell\times \{0\}, 
\gamma_\eta\subset
\{0\}\times \mathbb{R}^{m-\ell}$ and hence $\gamma_\eta\circ
\gamma_\sigma^{-1}$ is a geodesic in ${\cal G}$
connecting $\Phi(\sigma)$ to $\Phi(\eta)$. 
Since $\Phi$ is an isometry of 
$E(\Theta_-(\sigma,\eta),\lambda)$ 
onto a connected subgraph of ${\cal G}$ 
and $E(\Theta_-(\sigma,\eta),\lambda)$ is 
a strictly convex subset of $E(\tau,\lambda)$, 
we conclude that there is a geodesic in 
$E(\tau,\lambda)$ connecting $\sigma$ to $\eta$
which passes through $\Theta_-(\sigma,\eta)$.
This shows the first two statements of the lemma.

On the other hand, the considerations in the previous three
paragraphs of this proof also show that there
is a unique vertex 
$\Theta_+(\sigma,\eta)\in E(\Theta_-(\sigma,\eta),\lambda)$
such that $\Phi(\Theta_+(\sigma,\eta))=\Phi(\sigma)+\Phi(\eta)$.
The train track $\Theta_+(\sigma,\eta)$ satisfies the
properties in statement 2) and 3) of the lemma.
This completes the proof of the lemma.
\end{proof}

In the sequel we call $\tau$ the \emph{basepoint} of the flat
strip $E(\tau,\lambda)$. Another immediate consequence of
Lemma \ref{flatstrip} is the following growth control.

\begin{corollary}\label{growth}
For every $k>0$
the number of vertices in $E(\tau,\lambda)$
whose intrinsic distance to the basepoint is at most $k$
is not bigger than $(k+1)^m$.\end{corollary}
\begin{proof} Let $\Phi:E(\tau,\lambda)\to \mathbb{R}^m$
be an embedding of $E(\tau,\lambda)$ onto a cubical
graph in $\mathbb{R}^m$
as in Lemma \ref{flatstrip}. Assume that $\Phi$ maps the basepoint
of $E(\tau,\lambda)$ to $0$. 
By construction, if we denote
by $\vert \,\vert $ the norm on $\mathbb{R}^m$ defined
by $\vert x\vert =\sum_i\vert x_i\vert$
then $\Phi$ maps the ball of radius $k$ about $0$ in
$E(\tau,\lambda)$ into the intersection of the 
ball of radius $k$ about $0$ in $(\mathbb{R}^m,\vert\,\vert)$
with the cone $\{x\mid x_i\geq 0\}$.
Thus the image under $\Phi$ of the set of vertices in
$E(\tau,\lambda)$ whose distance to the basepoint is at most
$k$ is not bigger than the cardinality of the set
$\{x\in \mathbb{Z}^m\mid x_i\geq 0,\vert x\vert \leq k\}$,
and this cardinality
is not bigger than $(k+1)^m$.
\end{proof}

Call a train track $\sigma\in {\cal V}({\cal T\cal T})$ 
a \emph{full split}
of a train track $\tau\in {\cal V}({\cal T\cal T})$ 
if $\sigma$ can be
obtained from $\tau$ by splitting $\tau$ at 
\emph{each} large branch
precisely once. For a complete geodesic lamination $\lambda\in
{\cal C\cal L}$ which is carried by $\tau$ 
we call a full split $\sigma$ of $\tau$ a
\emph{full $\lambda$-split} if $\lambda$ is carried by $\sigma$.
Note that a full $\lambda$-split of $\tau$ is uniquely
determined by $\lambda$ and $\tau$: There are no choices involved.
A \emph{full $\lambda$-splitting sequence} of length $k\geq 0$ is a
sequence $\{\eta(i)\}_{0\leq i\leq k}\subset 
{\cal V}({\cal T\cal T})$ of
train tracks with the property that for every $i<k$ the train
track $\eta(i+1)$ is a full $\lambda$-split of $\eta(i)$.
A full $\lambda$-splitting sequence of length $k\geq 0$
issuing from $\tau$ is unique. We call its endpoint the \emph{full
$k$-fold $\lambda$-split} of $\tau$. The following corollary is
immediate from Lemma \ref{flatstrip}.

\begin{corollary}\label{fullsplit}
Let $\{\tau(i)\}_{0\leq i\leq p}$ be a
splitting sequence of length $p$. If $\lambda\in {\cal C\cal L}$
is carried by $\tau(p)$ then $\tau(p)$ is splittable to the full
$p$-fold $\lambda$-split of $\tau(0)$.
\end{corollary}
\begin{proof}
If $\{\tau(i)\}_{0\leq i\leq p}$ is any splitting
sequence of length $p$ 
and if $\lambda\in {\cal C\cal L}$ is carried by $\tau(p)$
then $\tau(p)\in E(\tau(0),\lambda)$ and the same is true for the
full $p$-fold $\lambda$-split $\sigma$ of $\tau(0)$. Thus the
corollary is immediate from Lemma \ref{flatstrip}. \end{proof}

\section{Quasi-geodesics}

In Section 5 we defined
for a complete train track $\tau$ which
carries a complete geodesic lamination $\lambda$
a flat cone $E(\tau,\lambda)\subset
{\cal T\cal T}$, and we investigated its intrinsic path metric 
$d_E$. We showed that
$(E(\tau,\lambda),d_E)$
is isometric to a cubical graph in $\mathbb{R}^m$ where $m>0$
is the number of branches of a complete train track on $S$.
By definition, the inclusion $(E(\tau,\lambda),d_E)\to {\cal T\cal T}$
is a one-Lipschitz map. 

The goal of this section is to relate the intrinsic path
metric on $E(\tau,\lambda)$ to the restriction of the 
metric $d$ on ${\cal T\cal T}$. For this
recall from the introduction that 
a map $F:(X,d)\to (Y,d)$ 
between two metric spaces $(X,d)$ and $(Y,d)$
is an \emph{$L$-quasi-isometric embedding}
if
\[d(x,y)/L-L\leq d(Fx,Fy)\leq Ld(x,y)+L\forall x,y\in X.\]
If moreover $F(X)$ is \emph{$L$-dense} in $Y$, i.e. if   
for every $y\in Y$ there is some $x\in X$ with 
$d(Fx,y)\leq L$, then $F$ is called an \emph{$L$-quasi-isometry}.
We show.

\begin{theo}\label{qiembedding}
There is a number $L>0$ such that for every $\tau\in {\cal V}({\cal T\cal T})$
and every complete geodesic lamination $\lambda$ carried by $\tau$ the
inclusion $(E(\tau,\lambda),d_E)\to {\cal T\cal T}$ is an
$L$-quasi-isometric embedding.
\end{theo}

Since by Corollary \ref{splitine} 
splitting arcs in $E(\tau,\lambda)$
are geodesics, Theorem \ref{theorem3} from the introduction is an immediate
consequence of Theorem \ref{qiembedding}

The proof of Theorem \ref{qiembedding} consists of two steps.
In the first step we consider for a complete geodesic lamination
$\lambda$ on $S$ the metric graph ${\cal S}(\lambda)$ 
whose vertex set
${\cal V}({\cal S}(\lambda))\subset{\cal V}({\cal T\cal T})$ 
is the set of all complete train tracks on $S$ which carry
$\lambda$ and where two such vertices $\tau,\sigma$ are
connected by an edge of length one if and only if either
$\tau$ and $\sigma$ are connected by an edge 
in ${\cal T\cal T}$ or if $\tau$ can be obtained from
$\sigma$ by a single shift. Lemma \ref{following} and 
Lemma \ref{longcarry} 
together with Proposition \ref{shiftsplitcarry} immediately
imply that ${\cal S}(\lambda)$ is connected. In other words,
${\cal S}(\lambda)$ with its
intrinsic path metric $d_\lambda$ 
is a geodesic metric space. We show that there 
is a number $q>1$ not depending on $\lambda$ and there
is a $q$-quasi-isometry  
${\cal S}(\lambda)\to {\cal T\cal T}$ whose restriction
to the vertex set of ${\cal S}(\lambda)$ is just the inclusion.
In a second step, we then establish that the 
inclusion $E(\tau,\lambda)\to {\cal S}(\lambda)$ is a 
uniform quasi-isometric
embedding. 

We begin with analyzing the intrinsic geometry of 
the metric graph ${\cal S}(\lambda)$.
For this define a 
\emph{splitting and shifting
sequence} to be a sequence $\{\alpha(i)\}\subset {\cal V}({\cal T\cal T})$ 
such that 
for each $i$, the train track $\alpha(i+1)$ is obtained
from $\alpha(i)$ either by a single
split or a single shift. 
The following important result of Penner and Harer
(Theorem 2.4.1 of \cite{PH92}) relates splitting and
shifting of train tracks to carrying.

\begin{proposition}\label{shiftsplitcarry}\cite{PH92}
If $\sigma\in {\cal V}({\cal T\cal T})$ is carried
by $\tau\in {\cal V}({\cal T\cal T})$ then $\tau$
can be connected to $\sigma$ by a splitting and shifting
sequence.
\end{proposition}

We also need the following
local version of Proposition \ref{shiftsplitcarry}.

\begin{lemma}\label{local2}
For every $k>0$ there is a number $p_1(k)>0$ with the following
property. Let $\sigma\in {\cal V}({\cal T\cal T})$ be carried by
$\tau\in {\cal V}({\cal T\cal T})$.
\begin{enumerate}
\item
If $d(\tau,\sigma)\leq k$
then $\tau$ can be connected to $\sigma$ by a splitting and shifting
sequence of length at most $p_1(k)$. 
\item If $\tau$ can be connected to $\sigma$ by a splitting and
shifting sequence of length at most $k$ then
$d(\tau,\sigma)\leq p_1(k)$.
\end{enumerate}
\end{lemma}
\begin{proof}
Up to the action of the mapping class group, for every $k>0$ 
there are only finitely many pairs $\sigma\prec\tau$ 
of complete train tracks whose 
distance is at most $k$ or such that $\sigma$ can
be obtained from $\tau$ by a splitting and shifting sequence
of length at most $k$. Thus the lemma follows from
Proposition \ref{shiftsplitcarry} and invariance
under the action of the mapping class group.
\end{proof}

To compare the metric space 
$({\cal S}(\lambda),d_\lambda)$ to the train track complex
${\cal T\cal T}$, 
we use the following uniform
local control on the graphs ${\cal S}(\lambda)$.

\begin{lemma}\label{local}
For every $k>0$ there is a number $p_2(k)>0$
with the following property. Let $\tau,\eta
\in {\cal V}({\cal T\cal T})$ be two complete train
tracks of distance at most $k$ which carry a common complete geodesic
lamination $\lambda$. Then there is a train track $\sigma$ which
carries $\lambda$, which is carried by both
$\tau,\eta$ and whose distance to $\tau,\eta$ is at most
$p_2(k)$.
\end{lemma}
\begin{proof}
Let $\alpha,\beta$ be two complete train tracks on $S$
and let ${\cal
C\cal L}(\alpha,\beta)$ be the set of all complete geodesic
laminations which are carried by both $\alpha$ and $\beta$. 
By Lemma \ref{recchar}, ${\cal C\cal
L}(\alpha,\beta)$ is an open subset of ${\cal
C\cal L}$. By Lemma \ref{following} and Lemma \ref{longcarry},
for every $\nu\in {\cal C\cal L}(\alpha,\beta)$ there is a
complete train track $\eta$ which carries $\nu$ and is carried by
both $\alpha,\beta$. 

By Lemma \ref{recchar},
for every complete train track $\eta$, 
the set ${\cal C\cal L}(\eta)$ of all complete geodesic
laminations which are carried by
$\eta$ is an open
subset of ${\cal C\cal L}$. By
compactness of ${\cal C\cal L}(\alpha,\beta)$, there are
\emph{finitely} many train tracks $\eta_1,\dots,\eta_k$ carried by
both $\alpha$ and $\beta$ and such that ${\cal C\cal
L}(\alpha,\beta)=\cup_{i=1}^k{\cal C\cal L}(\eta_i)$. 
In other words, there is
a number $\ell(\alpha,\beta)>0$ with the following property.
For every $\nu\in {\cal
C\cal L}(\alpha,\beta)$ 
there is a train track $\eta\in {\cal V}({\cal S}(\nu))$ which
carries $\nu$, which 
is carried by both $\alpha$ and $\beta$ and whose
distance to $\alpha,\beta$ is 
at most $\ell(\alpha,\beta)$. 

For $k>0$, up to the action of the mapping
class group there are only finitely many pairs $(\alpha,\beta)$ of
complete train tracks whose distance in ${\cal T\cal T}$ is at most
$k$ and which carry a common complete geodesic
lamination. 
By invariance under the action of the
mapping class group, this implies that there is a number
$p_2(k)>0$ such that the following holds true. For every complete
geodesic lamination $\lambda$ and for
every pair $\alpha,\beta\in {\cal V}({\cal S}(\lambda))$ 
with $d(\alpha,\beta)\leq k$ 
there is a train track $\eta\in {\cal V}({\cal S}(\lambda))$ 
with $\max\{d(\eta,\alpha),d(\eta,\beta)\}\leq p_2(k)$ and which
is carried by both $\alpha,\beta$.
This shows the lemma. \end{proof}

The set ${\cal V}({\cal S}(\lambda))$ of 
vertices of the metric graph ${\cal S}(\lambda)$ is contained
in the set ${\cal V}({\cal T\cal T})$ of vertices of 
the train track complex ${\cal T\cal T}$.
The inclusion ${\cal V}({\cal S}(\lambda))\to {\cal V}({\cal T\cal T})$ 
can be extended to a map
${\cal S}(\lambda)\to {\cal T\cal T}$. Since we are only
interested in the large-scale geometric properties of such a map, 
we do not need to give a precise definition. We only
require
that a point on an edge of ${\cal S}(\lambda)$
is mapped to a point in ${\cal T\cal T}$ of uniformly
bounded distance
to the endpoints of the edge, viewed as vertices in 
${\cal T\cal T}$. 
Since edges in ${\cal S}(\lambda)$ which
are not edges in ${\cal T\cal T}$ connect two 
complete train tracks which can be obtained from
each other by a single shift and hence
whose distance
in ${\cal T\cal T}$ is uniformly bounded, such a map clearly exists.
We call such a map
\emph{natural}. 
The next proposition is the first step toward Theorem 
\ref{qiembedding}.

\begin{proposition}\label{quasiisometry} 
There is a number $q_1>0$ such that
for every $\lambda\in {\cal C\cal L}$
a natural map $({\cal S}(\lambda),d_\lambda)\to {\cal T\cal T}$
is a $q_1$-quasi-isometry.
\end{proposition}
\begin{proof} Let ${\cal S\cal T}$ be the metric
graph whose vertices are the complete train tracks
on $S$ and where two vertices $\tau,\sigma$ are connected
by an edge of length one if and only if either they are connected
by an edge in ${\cal T\cal T}$ or if $\tau$ can be obtained from
$\sigma$ by a single shift. Since ${\cal T\cal T}$ is
connected, the same is true for ${\cal S\cal T}$.
Then ${\cal S\cal T}$ is
a locally finite metric graph which admits a
properly discontinuous cocompact 
isometric action of ${\cal M\cal C\cal G}$.
As a consequence, ${\cal S\cal T}$ is
equivariantly quasi-isometric
to ${\cal T\cal T}$. For every $\lambda\in 
{\cal C\cal L}$ the graph ${\cal S}(\lambda)$ is a complete
subgraph of ${\cal S\cal T}$. It is now enough to
show that there is a number $L>1$ not depending
on $\lambda$ such that
the natural inclusion ${\cal S}(\lambda)\to
{\cal S\cal T}$ is an $L$-quasi-isometry.
For simplicity of notation, for the remainder of this proof
we denote by $d$ the distance function of the metric graph
${\cal S\cal T}$.

By Lemma \ref{recchar}, 
for every $\tau\in {\cal V}({\cal T\cal T})$ 
the set of all complete
geodesic laminations which are carried by $\tau$ is
an open subset of ${\cal C\cal L}$. Since ${\cal C\cal L}$
is a compact space, there is a finite set ${\cal E}\subset
{\cal V}({\cal T\cal T})$ so that every complete geodesic lamination
$\lambda$ is carried by a train track $\tau\in {\cal E}$.
On the other hand, the mapping class group acts cocompactly
on ${\cal S\cal T}$, and it acts as a group of homeomorphisms
on ${\cal C\cal L}$. Thus by equivariance under the
action of the mapping class group, 
there is a number $D>0$ and for
every $\sigma\in {\cal S\cal T}$ there is some 
$\tau\in {\cal V}({\cal S}(\lambda))$ with 
$d(\tau,\sigma)\leq D$.
Since ${\cal S}(\lambda)$
is a complete subgraph of ${\cal S\cal T}$
this means that the inclusion ${\cal S}(\lambda)\to {\cal
S\cal T}$ is a $1$-Lipschitz map with $D$-dense image.

Let again $d_\lambda$ be the intrinsic
distance on ${\cal S}(\lambda)$. We
have to show that there is a number $L>1$ 
not depending on $\lambda$ 
such that
$d_\lambda(\tau,\sigma)\leq Ld(\tau,\sigma)$ for all
$\tau,\sigma\in {\cal V}({\cal S}(\lambda))$.
For this let $\tau,\sigma$
be any two vertices of
${\cal S}(\lambda)$. Let $\gamma:[0,m]\to {\cal S\cal T}$ be a
simplicial geodesic connecting $\tau=\gamma(0)$ to
$\sigma=\gamma(m)$ (i.e. $\gamma$ maps integer points in
$\mathbb{R}$ to vertices of ${\cal S\cal T}$). 
By the above consideration, for every $i\leq
m$ there is a train track $\zeta(i)\in {\cal V}({\cal S}(\lambda))$
with $d(\zeta(i),\gamma(i))\leq
D$ and where
$\zeta(0)=\tau,\zeta(m)=\sigma$.
Then the distance in ${\cal S\cal T}$
between $\zeta(i)$ and $\zeta(i+1)$ 
is at most $2D+1$. 

By Lemma \ref{local} and Lemma \ref{local2},
there is a constant $\kappa=p_1(p_2(D))>0$ only depending on
$D$ and for every $i\leq m$ there is a train track
$\beta(i)\in {\cal V}(\cal S(\lambda))$ 
which carries $\lambda$, which is carried by both 
$\zeta(i)$ and $\zeta(i+1)$ and
which can
be obtained from both $\zeta(i)$ and $\zeta(i+1)$
by a splitting and shifting
sequence of length at most $\kappa$. 
As a consequence, we have
$d_\lambda(\zeta(i-1),\zeta(i))\leq
2\kappa$ for all $i\leq m$. 
But this just means that the distance
$d_\lambda(\tau,\sigma)$
in ${\cal S}(\lambda)$ between $\tau$ and $\sigma$ is not bigger
than $2\kappa d(\tau,\sigma)$.
In other words, the inclusion $({\cal S}(\lambda),d_\lambda)\to
{\cal S\cal T}$ is a $2\kappa$-quasi-isometry for
the constant $2\kappa>0$ not depending on $\lambda$.
This shows the proposition. \end{proof}

For the proof of Theorem \ref{qiembedding} we are left 
with showing that for every $\lambda\in {\cal C\cal L}$
and every complete train track $\tau$ which carries $\lambda$,
the inclusion $E(\tau,\lambda)\to {\cal S}(\lambda)$
is an $L$-quasi-isometric embedding for a number $L>1$ which 
does not depend on $\lambda$. For this denote
for a complete geodesic 
lamination $\lambda\in {\cal C\cal L}$
and for a complete train track
$\tau\in {\cal V}({\cal S}(\lambda))$ which carries $\lambda$
by $C(\tau,\lambda)\subset {\cal S}(\lambda)$  
the complete subgraph of ${\cal S}(\lambda)$ 
whose vertex set is the set of 
all complete train tracks $\eta$ which are carried
by $\tau$ and which carry $\lambda$. 
By Proposition \ref{shiftsplitcarry}, 
Lemma \ref{following} and Lemma \ref{longcarry},
$C(\tau,\lambda)$ is 
connected and hence a geodesic metric space in its own right.
Moreover, it contains the flat cone
$E(\tau,\lambda)$ as a connected
subgraph.
We denote the intrinsic metric on $C(\tau,\lambda)$ by
$d_C$, and we denote as before by $d_E$ 
the intrinsic metric on $E(\tau,\lambda)$.
We first show that the natural inclusion
$(E(\tau,\lambda),d_E)\to (C(\tau,\lambda),d_C)$ is a uniform
quasi-isometry. In a second step, we establish that
the natural inclusion $(C(\tau,\lambda),d_C)
\to ({\cal S}(\lambda),d_\lambda)$
is a uniform quasi-isometric embedding. This then completes
the proof of Theorem \ref{qiembedding}. 

The proof of the following lemma 
relies on Proposition \ref{carrynear} 
from the appendix.

\begin{lemma}\label{quasiisom}
There is a number $q_2>0$ such that for every 
complete geodesic lamination $\lambda\in {\cal C\cal L}$
and every complete train track 
$\tau$ which carries $\lambda$
the inclusion $(E(\tau,\lambda),d_E)\to (C(\tau,\lambda),d_C)$ is
a $q_2$-quasi-isometry.
\end{lemma}
\begin{proof} We begin with showing that there is a constant
$k_1>0$ such that the subgraph
$E(\tau,\lambda)$ is $k_1$-dense in $C(\tau,\lambda)$
with respect to the intrinsic metric
$d_C$ on $C(\tau,\lambda)$.

For this let
$\chi>0$ be as in Proposition \ref{carrynear}.
Then by the definition of the graph ${\cal S}(\lambda)$,
the distance in ${\cal S}(\lambda)$ between a vertex
$\eta$ of the graph $C(\tau,\lambda)$ 
and its subgraph $E(\tau,\lambda)$ is at most
$\chi$. Thus by 
Lemma \ref{local} and Lemma \ref{local2},
there is
a complete train track $\sigma\in E(\tau,\lambda)$
and there is a complete train track 
$\zeta\in C(\tau,\lambda)$
which is carried by both $\sigma$ and $\eta$ and such that moreover
$\sigma,\eta$ can be connected to $\zeta$ by a 
splitting and shifting sequence of uniformly bounded length.
As a consequence, the distance in $C(\tau,\lambda)$ between
$\eta$ and $\sigma$ is uniformly bounded. This shows that
there is a constant $k_1>0$ such that
$E(\tau,\lambda)$ is $k_1$-dense in $(C(\tau,\sigma),d_C)$.

Since the inclusion $(E(\tau,\lambda),d_E)\to (C(\tau,\lambda),d_C)$
is clearly one-Lipschitz, for the proof of the lemma 
we are left with showing the existence of 
a universal constant
$k_2>0$ such that $d_E(\xi,\eta)\leq k_2d_C(\xi,\eta)$ for
all vertices $\xi,\eta\in E(\tau,\lambda)$.

Now both $E(\tau,\lambda)$
and $C(\tau,\lambda)$ are geodesic metric spaces and
$E(\tau,\lambda)$ is $k_1$-dense in $C(\tau,\lambda)$. Therefore
it is enough
to show the existence of a constant $k_3>0$ with the 
following property. If 
$\xi,\eta\in E(\tau,\lambda)$ are any two vertices
with $d_C(\xi,\eta)\leq 3k_1$
then $d_E(\xi,\eta)\leq k_3$.

Namely, assume that this property holds true. Let
$\xi,\eta\in E(\tau,\lambda)$
be arbitrary vertices with $d_C(\xi,\eta)=d>0$. Let  
$\gamma:[0,d]\to C(\tau,\lambda)$ be a simplicial geodesic
connecting 
$\gamma(0)=\xi$ to 
$\gamma(d)=\eta$. Since $E(\tau,\lambda)$ is $k_1$-dense in 
$C(\tau,\lambda)$ we can replace $\gamma$ 
by a simplicial path $\tilde \gamma:[0,\tilde d]\to C(\tau,\lambda)$ 
of length $\tilde d\leq 3d$  
with the same endpoints and the additional property that
$\gamma(3\ell k_1)\in E(\tau,\lambda)$ for all
integers $\ell\leq d/k_1$. The arc $\tilde \gamma[0,3k_1]$
is obtained by concatenation of $\gamma[0,2k_1]$ with an arc
of length at most $k_1$ which connects $\gamma(2k_1)$
with a point in $E(\tau,\lambda)$. Inductively,
for each $\ell\leq d/k_1$ the arc 
$\tilde \gamma[3(\ell-1)k_1,3k_1]$ is 
up to parametrization a concatentation of 
a segment of length at most
$k_1$ connecting a point in $E(\tau,\lambda)$ to 
$\gamma((\ell-1) k_1)$, the arc $\gamma[(\ell-1)k_1,\ell k_1]$
and an arc of length at most $k_1$ connecting
$\gamma(\ell k_1)$ to a point in $E(\tau,\lambda)$.
Replace each of the arcs
$\tilde\gamma[3(\ell-1)k_1,3\ell k_1]$ of length at most
$3k_1$ with endpoints
in $E(\tau,\lambda)$ by an arc of length at most 
$k_3$ which is contained in $E(\tau,\lambda)$.
The resulting path is contained in $E(\tau,\lambda)$, it connects
$\xi$ to $\eta$ and its length
does not exceed 
$k_3d_C(\xi,\eta)/k_1$. 

To show the existence of a constant $k_3>0$ with the above
properties,
let $\xi,\eta\in E(\tau,\lambda)$ be vertices 
such that $d_C(\xi,\eta)\leq 
3k_1$ and 
let $\zeta\in C(\tau,\lambda)$ be a 
complete train track which carries $\lambda$ and which
can be obtained from 
both $\xi,\eta$ by a splitting and shifting sequence 
whose length is bounded from above by a universal constant $p>0$.
Such a complete train track exists
by Lemma \ref{local} and Lemma \ref{local2} and the fact
that the distance in ${\cal T\cal T}$ between
$\xi,\eta$ is uniformly bounded.

Since $\zeta$ is carried by $\tau$, there is a train track
$\beta\in E(\tau,\lambda)$
which carries $\zeta$ 
and such that no split of 
$\beta$ carries $\zeta$. It follows from Lemma \ref{specialtrain}
in the apppendix 
that $\beta$ is unique.
Now $\xi\in E(\tau,\lambda)$ 
carries $\zeta$ and therefore by Lemma \ref{flatstrip},
$\xi$ is splittable to $\beta$. The same argument
also shows that 
$\eta$ is splittable to $\beta$. 

On the other hand, up to the action of the
mapping class group, there are only finitely many pairs of 
complete train tracks $(\alpha,\rho)$ such that
$\rho$ can be obtained from $\alpha$ by a splitting and shifting
sequence of length at most $p$. Since $\xi$ is splittable
to $\beta$, $\beta$ carries $\zeta$ and $\zeta$ can be obtained
from $\xi$ by a splitting and shifting sequence of length at most
$p$, the distance in ${\cal T\cal T}$ 
between $\beta$ and $\xi,\zeta$ is uniformly
bounded. The same argument shows that the distance
between $\beta$ and $\eta$ is uniformly bounded.
Then the
distance in $E(\tau,\lambda)$ between $\beta$ and $\xi,\eta$ is 
uniformly bounded as well.
As a consequence, $\xi$ can be
connected to $\eta$ by a path in $E(\tau,\lambda)$ 
of uniformly bounded length which is the concatenation of a
splitting sequence connecting $\xi$ to $\beta$ and a collapsing
sequence connecting $\beta$ to $\eta$.
This completes the proof of the lemma.
\end{proof}

Lemma \ref{local}, 
Lemma \ref{quasiisom} and Lemma \ref{convex}
are used to show the following.

\begin{lemma}\label{nesting}
For every $k>0$ there is a number $p_3(k)>0$ with the following
property. Let $\tau,\sigma$ be complete train tracks
which carry a common complete geodesic lamination $\lambda$ and
such that $d(\tau,\sigma)\leq k$. Then for every
$\tau^\prime\in C(\tau,\lambda)$ there is a
complete train track 
$\sigma^\prime\in C(\tau^\prime,\lambda)\cap C(\sigma,\lambda)$ with
$d(\tau^\prime,\sigma^\prime)\leq p_3(k)$.
\end{lemma}
\begin{proof} 
Let $k>0$, let $\lambda\in {\cal C\cal L}$
and let $\tau,\sigma\in {\cal V}({\cal T\cal T})$ be complete train
tracks which carries $\lambda$. By Lemma \ref{local} 
and Lemma \ref{local2} there is 
a complete train track $\xi\in C(\tau,\lambda)\cap C(\sigma,\lambda)$
which can be obtained from $\tau$ by a splitting and shifting
sequence of length at most $p_1(p_2(k))$. Since
$C(\xi,\lambda)\subset C(\sigma,\lambda)$, it is enough
to show the lemma under the additional assumption that
$\sigma\in C(\tau,\lambda)$.

Since $C(\eta,\lambda)=C(\eta^\prime,\lambda)$ if 
$\eta,\eta^\prime$ are shift equivalent, this means that for 
the proof of the lemma, it is in fact enough to show the existence
of a constant $L>0$ with
the following property.
If the complete train track 
$\sigma$ carries $\lambda\in {\cal C\cal L}$
and can be obtained from a complete train track $\tau$ by
a single split at a large branch $e$ then the 
Hausdorff distance between $C(\tau,\lambda)$ and
$C(\sigma,\lambda)$ as subsets of ${\cal T\cal T}$ is at most $L$.

Now by Proposition \ref{carrynear} and 
Lemma \ref{local2}, the Hausdorff distance 
in ${\cal T\cal T}$ between
$C(\tau,\lambda)$ and $E(\tau,\lambda)$ is at most
$p_1(\chi)$ where $\chi>0$ is as in Proposition \ref{carrynear},
and similarly for $C(\sigma,\lambda)$ and $E(\sigma,\lambda)$.
Moreover, by Lemma \ref{convex}, the Hausdorff distance
between $E(\tau,\lambda)$ and $E(\sigma,\lambda)$ equals $1$.
Together we conclude that the Hausdorff distance between
$C(\tau,\lambda)$ and $C(\sigma,\lambda)$ does not exceed
$2p_1(\chi)+1$. From this the lemma follows.
\end{proof}

The following lemma is the main remaining step for the proof
of Theorem \ref{qiembedding}.
For its formulation, denote as before by
$d_E$ the intrinsic path metric on a flat cone $E(\tau,\lambda)$.
Recall moreover the definition of the metric graph
${\cal S}(\lambda)$ for a complete geodesic lamination $\lambda$.

\begin{lemma}\label{pushforward} 
There is a number $p_4>0$ with
the following property. Let $\lambda$ be a complete geodesic
lamination, let $\tau$ be a complete train track
which carries $\lambda$ and let
$\sigma,\eta\in E(\tau,\lambda)$. 
If $\gamma:[0,m]\to {\cal S}(\lambda)$ is
any simplicial path 
connecting $\gamma(0)=\sigma$ to
$\gamma(m)=\eta$ then the
length of $\gamma$ is not smaller than $d_E(\sigma,\eta)/p_4$.
\end{lemma}
\begin{proof} 
Let $\sigma,\eta\in E(\tau,\lambda)$ and let 
$\xi=\Theta_-(\sigma,\eta)$ be as in Lemma \ref{distinflat}.
Let $\ell_1=d_E(\xi,\sigma),\ell_2=d_E(\xi,\eta)$ and assume
that $\ell_1\leq \ell_2$. By Lemma \ref{distinflat},
there is a train track $\Theta_+(\sigma,\eta)$ which can
be obtained from $\sigma$ by a splitting sequence of 
length $\ell_2$ and which can be obtained from 
$\eta$ by a splitting sequence of length $\ell_1$.

Let $p_3(1)>0$ be as in Lemma \ref{nesting} and let
$p=p_1(p_3(D))$ be as in Lemma \ref{local2}.
Let $\gamma:[0,n]\to {\cal S}(\lambda)$ be
a simplicial path connecting
$\sigma=\gamma(0)$ to $\eta=\gamma(n)$.
We construct inductively a sequence 
$\{\sigma(i)\}_{0\leq i\leq n}
\subset {\cal V}({\cal S}(\lambda))$ with the following properties.
\begin{enumerate}
\item For every $i\leq n$, $\sigma(i)$ is carried by $\gamma(i)$
and carries $\lambda$.
\item For every
$i<n$ the train track $\sigma(i+1)$ can be obtained from
$\sigma(i)$ by a splitting and shifting
sequence of length at most $p$.
\end{enumerate}

For the construction of the sequence $\{\sigma(i)\}$,
we first define $\sigma(0)=\sigma$. Assume by induction that
we constructed already the train tracks $\sigma(i)$ for
all $i<i_0$ and some $i_0>0$. Consider the train track
$\gamma(i_0-1)$; by assumption, it carries the
train track $\sigma(i_0-1)$. If the train track $\gamma(i_0)$
is obtained from $\gamma(i_0-1)$ by a collapse or a shift
then $\gamma(i_0)$ carries $\sigma(i_0-1)$ and we
define $\sigma(i_0)=\sigma(i_0-1)$. Otherwise
$\gamma(i_0)$ is obtained from $\gamma(i_0-1)$
by a single $\lambda$-split.
By Lemma \ref{nesting} and Lemma \ref{local2},
there is a train track
$\sigma(i_0)$ which carries $\lambda$,
is carried by $\gamma(i_0)$ and which
can be obtained from $\sigma(i_0-1)$
by a splitting and shifting sequence of length at most $p$.
The inductively defined
sequence $\{\sigma(i)\}$ has the required properties.

By construction, the train track $\sigma=\sigma(0)$ can
be connected to $\sigma(n)$ by a splitting
and shifting sequence of length at most $pn$. In particular,
we have $\sigma(n)\in C(\sigma,\lambda)$. 
By Lemma \ref{quasiisom}, Lemma \ref{flatstrip} and
Corollary \ref{splitine}, 
the train
track $\sigma$ is splittable with a sequence of $\lambda$-splits
of length at most $q_2pn+q_2$ to a train track
$\nu\in E(\sigma,\lambda)\subset E(\tau,\lambda)$ 
which is contained in the $q_2$-neighborhood
of $\sigma(n)$ in $C(\tau,\lambda)$. 
Since $\sigma(n)$ is carried by $\gamma(n)=\eta\in E(\tau,\lambda)$, 
via replacing $\nu$ by a train track in a uniformly bounded
neighborhood we may assume that $\eta$ is splittable to $\nu$ as
well. Lemma \ref{distinflat} then shows that
the train track $\Theta_+(\sigma,\eta)\in E(\tau,\lambda)$
is splittable to $\nu$. 

From Lemma \ref{flatstrip} we deduce that
the length $\ell_2\geq \ell_1$ of a splitting sequence connecting
$\sigma$ to $\Theta_+(\sigma,\eta)$ 
is not bigger than $q_2pn+L$ where $L>q_2$ is a universal
constant. 
Since the length of a geodesic in $E(\tau,\lambda)$ connecting
$\sigma$ to $\eta$ equals $\ell_1+\ell_2\leq 2\ell_2$, 
we conclude that the intrinsic distance in 
$E(\tau,\lambda)$ between $\sigma$ and 
$\eta$ does not exceed $2(q_2pn+L)$. 
This is just
the statement of the lemma.
\end{proof}

Theorem \ref{qiembedding} is now an immediate consequence
of Lemma \ref{pushforward}. Namely, by Lemma \ref{quasiisometry}
we only have to show that there is a number $L>0$ such that
for every $\lambda\in {\cal C\cal L}$ and every 
train track $\tau\in {\cal V}({\cal S}(\lambda))$, 
if $\sigma,\eta\in E(\tau,\lambda)$ 
can be connected in ${\cal S}(\lambda)$ by a simplicial path
of length $n\geq 0$ then there is an arc in $E(\tau,\lambda)$
connecting $\sigma$ to $\eta$ whose length does not exceed
$nL+L$. However, this was shown in Lemma \ref{pushforward}.
This completes the proof of Theorem \ref{qiembedding}. 
\qed

\bigskip

Theorem \ref{qiembedding} implies the result of Farb, Lubotzky
and Minsky \cite{FLM01}.

\begin{corollary} If $P$ is a pants
decomposition for $S$ and if $\Gamma$
is the free abelian group of rank $3g-3+m$
generated by the Dehn twists about the
pants curves of $P$ then $\Gamma<{\cal M\cal C\cal G}$
is undistorted.\end{corollary}
\begin{proof} Let $P=\{\gamma_1,\dots,
\gamma_{3g-3+m}\}$ be any pants decomposition
for $S$. For each $i\leq 3g-3+m$ let
$\phi_i\in {\cal M\cal C\cal G}$ be a simple (positive or negative) Dehn
twist about $\gamma_i$. The 
elements $\phi_1,\dots,\phi_{3g-3+m}$ generate
a free abelian subgroup $\Gamma$ of ${\cal M\cal C\cal G}$.
We equip $\Gamma$ with the word norm $\vert \,\vert$
defined by the generators $\phi_i,\phi_i^{-1}$.
To show that $\Gamma$ is undistorted,
it suffices to show that there is a $\Gamma$-equivariant
quasi-isometric
embedding of the semi-group
$\Gamma_+=\{\phi_1^{\ell_1}\circ\cdots\circ 
\phi_{3g-3+m}^{\ell_{3g-3+m}}\mid \ell_i\geq 0\}$ into
a flat cone $E(\tau,\lambda)$ for some $\lambda\in {\cal C\cal L}$. 
Namely, if this holds true then the corollary
follows from Theorem \ref{qiembedding} and 
the fact that ${\cal T\cal T}$
is ${\cal M\cal C\cal G}$-equivariantly quasi-isometric
to ${\cal M\cal C\cal G}$.

For this
choose a train track $\tau\in {\cal V}({\cal T\cal T})$ 
which is obtained by collapsing a small tubular neighborhood of 
a finite complete geodesic lamination $\lambda$ on $S$ as
in Section 2 whose minimal components are precisely the components of $P$.
Up to replacing $\tau$ by a shift equivalent train track, every component
$\gamma_i$ of $P$ is carried by an embedded trainpath of length $2$ in 
$\tau$ which consists of a large branch and a small branch
(note that $\tau$ is just 
a complete train track \emph{in standard form}
for the pants decomposition $P$ of $S$ as defined in 
\cite{PH92}). Moreover, every large branch of $\tau$ is of this form.
For a suitable
choice of the spiraling directions of $\lambda$ about the components
of $P$, for every $i\in \{1,\dots,3g-3+m\}$ the train
track $\tau$ is splittable to 
$\phi_i\tau$ with a splitting sequence
of length $2$ (with two splits at a large branch
contained in the embedded trainpath $\gamma_i$,
with the small branch in $\gamma_i$ as a winner).
Moreover, for $i\not=j$ these splitting sequences commute.
As a consequence, the flat cone $E(\tau,\lambda)$ is 
invariant under the semi-group $\Gamma_+$, and the 
map which associates to
an element $\phi\in \Gamma_+$ the train track
$\phi(\tau)\in E(\tau,\lambda)$ is an equivariant 
quasi-isometry
between $\Gamma_+$ and $E(\tau,\lambda)$.
\end{proof}

\section{Boundary amenability}

The mapping class group naturally acts on
the compact Hausdorff space ${\cal C\cal L}$ 
of all complete geodesic laminations on $S$ 
as a group of
homeomorphisms. We show in this section that this
action is topologically amenable.

For this we use the assumptions and notations from the previous
sections. In particular, for
$\tau\in {\cal V}({\cal T\cal T})$ and a complete geodesic
lamination $\lambda$ carried by $\tau$ let as before
$C(\tau,\lambda)$ be the graph whose vertices are the
complete train tracks which are carried by $\tau$ and
carry $\lambda$ and where two such vertices $\sigma,\eta$
are connected
by an edge of length one if and only if either $\sigma$
can be obtained from $\eta$ by a single shift or
$\sigma,\eta$ are connected in ${\cal T\cal T}$ by an edge
of length one. The intrinsic path metric on $C(\tau,\lambda)$ is
denoted as before by $d_C$.

Similarly, the flat cone $E(\tau,\lambda)$ is
the full subgraph of ${\cal T\cal T}$ whose vertices
are the complete train tracks which carry
$\lambda$ and which can be obtained from $\tau$ by
a splitting sequence. We call $\tau$ the
\emph{basepoint} of $E(\tau,\lambda)$. The intrinsic path metric
on $E(\tau,\lambda)$ is denoted as before by $d_E$.

Choose a
finite subset ${\cal E}$ of ${\cal V}({\cal T\cal T})$
such that $\cup_{\phi\in {\cal M\cal C\cal G}}\phi{\cal E}=
{\cal V}({\cal T\cal T})$ and that moreover 
for every $\lambda\in {\cal C\cal L}$ there is a train track
$\eta\in {\cal E}$ which carries $\lambda$.
By equivariance under the action of the
mapping class group, this implies that 
for every $\phi\in {\cal M\cal C\cal G}$
and every $\lambda\in {\cal C\cal L}$ there is 
a train track $\eta\in \phi{\cal E}$ which carries $\lambda$.
As in the  
proof of Proposition \ref{quasiisometry},
such a set exists since ${\cal C\cal L}$ is compact and
since the set of all complete geodesic laminations which are
carried by a complete train track $\tau$ is open in
${\cal C\cal L}$.

Let ${\cal G}$ be a finite
symmetric set of generators for ${\cal M\cal C\cal G}$ containing
the identity. For an element $\phi\in {\cal M\cal C\cal G}$
let $\vert \phi\vert $ be the word norm of $\phi$ with
respect to the generating set ${\cal G}$.
We have.

\begin{lemma}\label{appro}
There are numbers $0<\kappa_1<\kappa_2$
with the following property.
Let $\phi\in {\cal M\cal C\cal G}$ be
such that $\vert \phi\vert =k$,
let $\lambda\in {\cal C\cal L}$ and 
let $\sigma\in
\phi{\cal E},\tau\in {\cal E}$ be train tracks which carry $\lambda$.
Then for $n\geq k$
the distance in ${\cal T\cal T}$ between
the full $n\kappa_1$-fold $\lambda$-split of $\sigma$ and the
ball in $(E(\tau,\lambda),d_E)$ of radius $n\kappa_2$
about the basepoint is at most
$\kappa_2$.
\end{lemma}
\begin{proof}
Since $\cup_{g\in {\cal G}}g{\cal E}$ is finite, its diameter $D$ 
in ${\cal T\cal T}$ is finite as well. Now ${\cal M\cal C\cal G}$
acts on ${\cal T\cal T}$ as group of isometries and consequently
the following holds true. 
Let $\phi\in {\cal M\cal C\cal G}$ and let $\psi\in {\cal G}$. If 
$\sigma\in \phi{\cal E}$ and if 
$\eta\in \phi\psi{\cal E}$ then we have
$d(\sigma,\eta)\leq D$.

Let $\lambda\in {\cal C\cal L}$ be a complete geodesic lamination
on $S$, let $\tau\in {\cal E}$ be a train track which
carries $\lambda$ and let 
$\phi\in {\cal M\cal C\cal G}$
with $\vert \phi\vert =k$. Then $\phi$ can be
represented in the form
$\phi=\phi_{1}\dots \phi_{k}$ with $\phi_i\in {\cal G}$.
Let $\sigma\in \phi{\cal E}$ be such that $\sigma$ carries
$\lambda$.
For $0\leq i\leq  k$ let $\tau_i\in \phi_{1}\dots \phi_{k-i}{\cal E}$
be a train track which carries $\lambda$ and such that 
$\tau_0=\sigma,\tau_k=\tau$.
Then for each $i$, 
the distance between $\tau_i$ and $\tau_{i-1}$
is bounded from above by $D$.

Let $p_3(D)>0$ be as in Lemma \ref{nesting}.
We construct inductively a  
sequence $(\eta_i)_{0\leq i\leq k}$ of complete
train tracks with the following properties.
\begin{enumerate}
\item $\eta_0=\sigma$.
\item For each $i$, the train track 
$\eta_i$ carries $\lambda$.
\item For each $i$, $\eta_i$ is carried by $\eta_{i-1}$ and $\tau_i$.
\item The distance between $\eta_i$ and $\eta_{i-1}$ 
does not exceed $p_3(D)$.
\end{enumerate}

For the construction, 
assume by induction that the train tracks
$\eta_j$ are already determined for all $j\leq i-1$
and some $i\geq 1$. To construct $\eta_i$, note that
since $\eta_{i-1}$ is carried by $\tau_{i-1}$ and since
$d(\tau_i,\tau_{i-1})\leq D$, 
by Lemma \ref{nesting} 
there is a complete train track $\eta_i$ 
which carries $\lambda$, which
is carried by both $\eta_{i-1}$ and $\tau_i$
and whose distance to
$\eta_{i-1}$ is a most $p_3(D)$. 
This yields the construction.

Since $\eta_k$ is carried by $\sigma$, i.e. 
$\eta_k\in C(\sigma,\lambda)$,
by Lemma \ref{quasiisom}
there is train
track $\sigma^\prime\in E(\sigma,\lambda)$ whose
distance to $\eta_k$ is at most $q_2$.
Then the distance between $\sigma^\prime$ and $\sigma$
does not exceed $kp_3(D)+q_2$ and therefore
by Theorem \ref{qiembedding} and Corollary \ref{splitine},
the length of a splitting sequence
connecting $\sigma$ to $\sigma^\prime$ is at most
$Lp_3(D) k +L(q_2+1)$ where $L>1$ is a universal constant.
By Corollary \ref{fullsplit}, $\sigma^\prime$ is splittable to
the train track $\sigma^{\prime\prime}$ 
obtained from $\sigma$ by a full splitting
sequence of length $Lp_3(D)k+L(q_2+1)$. 
If $m>0$ is the number of branches of a complete train
track on $S$, then a full split consists of at most 
$m$ single splits and hence
$\sigma^{\prime\prime}$ is obtained from
$\sigma$ by a splitting sequence in the usual
sense of length at most $mLp_3(D)k+m(Lq_2+1)$.

Since $\eta_k$ is carried by $\tau$ and the distance
between $\eta_k$ and $\sigma^\prime$ does not exceed $q_2$, 
Lemma \ref{quasiisom} shows 
that the train track $\sigma^{\prime\prime}$
is contained in the $2q_2$-neighborhood of 
the flat cone $E(\tau,\lambda)$. 
On the other hand, we have $d(\tau,\sigma^{\prime\prime})\leq 
d(\tau,\sigma)+d(\sigma,\sigma^{\prime\prime})\leq kD
+mLp_3(D)k+m(Lq_2+1)$ by the above consideration.
In other words, the distance between $\tau$ and
$\sigma^{\prime\prime}$ 
is bounded from above by $\nu k+\nu$ 
for a universal constant $\nu>0$.
Since $\sigma^\prime$ is splittable to $\sigma^{\prime\prime}$,
Lemma \ref{nesting} and 
Theorem \ref{qiembedding} 
show the existence of a universal constant $\kappa_2>0$
such that $\sigma^{\prime\prime}$ is contained 
in the $\kappa_2$-neighborhood of the 
ball of radius $\kappa_2 k$ in $E(\tau,\lambda)$. 
This shows the lemma.
\end{proof}

We are now
ready to complete the proof of Theorem \ref{theorem1} from the
introduction.

\begin{proposition}\label{topam}
The action of ${\cal M\cal C\cal G}$
on ${\cal C\cal L}$ is topologically amenable.
\end{proposition}
\begin{proof}
Let ${\cal E}\subset 
{\cal V}({\cal T\cal T})$ be as in Lemma \ref{appro}, and 
for $\lambda\in {\cal C\cal L}$ and 
$\phi\in {\cal M\cal C\cal G}$ let
$\tau(\lambda,\phi)\in \phi{\cal E}$ be a train track
which carries $\lambda$. 
Following Kaimanovich \cite{Ka04},
for $n\geq 1,k\leq 2n$ and $\phi\in {\cal M\cal C\cal G}$
define $Y(\phi,\lambda,n,k)$ to be
the set of all complete train tracks $\sigma\in 
{\cal V}({\cal T\cal T})$
with the following property. There is an element
$\psi\in {\cal M\cal C\cal G}$ 
with $\vert \psi\phi^{-1}\vert \leq k$ and there
is a complete
train track $\sigma_0\in \psi{\cal E}$
which carries $\lambda$
and such that $\sigma$ can be obtained from $\sigma_0$
by a full $2n\kappa_1$-fold $\lambda$-split where
$\kappa_1>0$ is as in Lemma \ref{appro}.

By Lemma \ref{appro} and invariance under the action of the
mapping class group,
the set $Y(\phi,\lambda,n,k)$ is contained in the
$\kappa_2$-neighborhood in ${\cal T\cal T}$
of the ball of radius $2n\kappa_2$ about the basepoint
in the flat cone
$E(\tau(\lambda,\sigma),\lambda)$ with respect to the intrinsic metric.
Since the number of vertices of ${\cal T\cal T}$ contained
in any ball of radius $\kappa_2$ is uniformly bounded,
by Corollary \ref{growth}
there is a number $C>0$ such that 
the number of points in
the set $Y(\phi,\lambda,n,k)$ is bounded from above
by $Cn^C$. In particular, the cardinality of the
sets $Y(\phi,\lambda,n,k)$ is bounded from above by
a universal polynomial in $n$.

By construction, for any $\phi,\phi^\prime\in 
{\cal M\cal C\cal G}$ with $q=\vert \phi^\prime\phi^{-1}\vert$ and 
for every $\lambda\in {\cal C\cal L}$  
the sets
$Y(\phi,\lambda,n,k)$,
$Y(\phi^\prime,\lambda,n,k^\prime)$
satisfy the following \emph{nesting condition} 
from Lemma 1.35 of \cite{Ka04}
(see the proof of Corollary 1.37 in \cite{Ka04}):
For every $n>q$ and every
$k\leq 2n-q$ we have
\begin{equation}\label{nesting5}
Y(\phi,\lambda,n,k)\subset Y(\phi^\prime,\lambda,n,k+q)\text{ and } 
Y(\phi^\prime,\lambda,n,k)\subset Y(\phi,\lambda,n,k+q). 
\end{equation}

Let ${\cal P}({\cal T\cal T})$
be the space of all probability measures on the set 
${\cal V}({\cal T\cal T})$ of all 
vertices of ${\cal T\cal T}$.
A probability measure on a countable
set $Z$ can be viewed as a non-negative function
$f$ on $Z$ with $\sum_{z\in Z}f(z)=1$.
In other words, such a probabiliy measure is a
point in the space $\ell^1(Z)$ of integrable functions
on $Z$ and hence 
the space of probability measures
on $Z$ can be equipped with
the usual Banach norm
on $\ell^1(Z)$. Thus ${\cal P}({\cal T\cal T})$
is equipped with a natural norm $\Vert\,\Vert$.

Let $m_{Y(\phi,\lambda,n,k)}\in {\cal P}({\cal T\cal T})$ 
be the normalized
counting measure for $Y(\phi,\lambda,n,k)$,
i.e. the probability
measure which is the normalization of the sum of the
Dirac measures on the points in 
$Y(\phi,\lambda,n,k)$. Define
\[\nu_{n}(\phi,\lambda)=\frac{1}{n}\sum_{k=1}^n
m_{Y(\phi,\lambda,n,k)}.\] 
Denote by $\vert B\vert$ the cardinality of a finite
subset of ${\cal V}({\cal T\cal T})$.
Since for each fixed $n$ 
the family of sets $Y(\phi,\lambda,n,k)$ satisfies
the nesting condition (\ref{nesting5}) above, 
Lemma 1.35 of \cite{Ka04}
shows that whenever $\vert \phi^\prime\phi^{-1}\vert \leq q$
and $n> q$ then we have 
\[\Vert \nu_{n}(\phi,\lambda)-\nu_{n}(\phi^\prime,\lambda)\Vert
\leq \frac{2q}{n}+
\frac{4(n-q)}{n}
\bigl[1-
\bigl(\frac{{\rm const}}
{\vert Y(\phi,\lambda,n,n+q)\vert}\bigr)^{2q/(n-q)}
\bigr].\]

The consideration in the second paragraph of this proof shows that
for $n>q$ the cardinality
of the set $Y(\phi,\lambda,n,n+q)$ is bounded from above
by a fixed polynomial in $n$.
As a consequence, we obtain that
\begin{equation}\label{convergence}
\Vert \nu_n(\phi,\lambda)-\nu(\phi^\prime,\lambda)\Vert \to 0
\end{equation}
as $n\to \infty$ and locally uniformly on ${\cal M\cal C\cal G}$.

For $n\geq 0$ and
$\lambda\in {\cal C\cal L}$ define 
\[\mu_n(\lambda)=\nu_{n}(e,\lambda).\] 
By (\ref{convergence}) above, 
the measures $\mu_n(\lambda)\in {\cal P}({\cal T\cal T})$
satisfy
$\Vert \mu_n(g\lambda)-g\mu_n(\lambda)\Vert \to 0$
$(n\to \infty)$ 
uniformly on compact subsets of 
${\cal C\cal L}\times {\cal M\cal C\cal G}$. 
Since
the action of ${\cal M\cal C\cal G}$ on ${\cal T\cal T}$
is isometric and properly
discontinuous, the measures $\mu_n(\lambda)$ on ${\cal T\cal T}$
can be lifted to measures on ${\cal M\cal C\cal G}$ with the same
property (see \cite{AR00} and compare \cite{Ka04}).
As a consequence, the action of ${\cal M\cal C\cal G}$ on
${\cal C\cal L}$ is topologically amenable. 
\end{proof}

{\bf Remark:} An action of a countable group $\Gamma$ on 
a compact Hausdorff space can only be topologically
amenable if the point stabilizers of this action
are amenable subgroups of $\Gamma$. The point stabilizers
of the action of ${\cal M\cal C\cal G}$ on ${\cal C\cal L}$ are
virtually abelian. For example, if $\lambda\in {\cal C\cal L}$
is a complete geodesic lamination whose minimal components
are simple closed curves, then the stabilizer of $\lambda$
in ${\cal M\cal C\cal G}$ contains the free abelian subgroup of
${\cal M\cal C\cal G}$ generated by the Dehn twists about these
components as a subgroup of finite index (this fact can
easily be deduced from the results in \cite{MP89}).
On the other hand, it is well known that an amenable
subgroup $\Gamma$ of ${\cal M\cal C\cal G}$ is virtually
abelian (see \cite{I02} for references). This implies that
every amenable subgroup of ${\cal M\cal C\cal G}$ has a
subgroup of finite index which fixes a point in
${\cal C\cal L}$. This corresponds to the properties
of the action of a simple Lie group of higher rank
on its Furstenberg boundary and hence we call ${\cal C\cal L}$
the \emph{Furstenberg boundary} of ${\cal M\cal C\cal G}$. 
In contrast, the
action of ${\cal M\cal C\cal G}$ on the Thurston boundary of
Teichm\"uller space is not topologically amenable.

\section{A strong boundary for ${\cal M\cal C\cal G}$}

A \emph{strong boundary} for a locally compact second
countable topological
group $\Gamma$ is a standard probability space $(X,\mu)$
with a measure class preserving action of $\Gamma$ and the
following two additional properties.
\begin{enumerate}
\item The $\Gamma$-space $(X,\mu)$ is amenable.
\item \emph{Double ergodicity}:
Let $(E,\pi)$ be any coefficient module for $\Gamma$;
then every $\Gamma$-equivariant weak$^*$-measurable map
$f:(X\times X,\mu\times\mu)\to E$ is constant almost
everywhere.
\end{enumerate}
Kaimanovich \cite{Ka03}
showed that for every locally compact second countable
topological group $\Gamma$, the Poisson boundary of every \'etal\'e
non-degenerate symmetric probability measure on $\Gamma$
is a strong boundary for $\Gamma$.

The mapping class group acts on the Teichm\"uller space
equipped with the \emph{Teichm\"uller metric} as
a group of isometries. Hence the restriction of
the Teichm\"uller distance to an orbit of ${\cal M\cal C\cal G}$
defines a ${\cal M\cal C\cal G}$-invariant distance function
on ${\cal M\cal C\cal G}$ (which however is not quasi-isometric
to the distance function defined by a word norm).
In \cite{KM96}, Kaimanovich and Masur investigate
the Poisson boundary for a symmetric probability measure
$\mu$ of finite entropy on the mapping class group ${\cal M\cal C\cal G}$.
They show that if the first logarithmic
moment of $\mu$ with respect to the Teichm\"uller distance
is finite and if the
support of $\mu$ generates ${\cal M\cal C\cal G}$, then the
Poisson boundary of $\mu$ can be viewed as a measure on
the space ${\cal P\cal M\cal L}$
of all projective measured geodesic laminations on $S$. The measure
class of this measure is ${\cal M\cal C\cal G}$-invariant and
gives full mass to the subset of ${\cal P\cal M\cal L}$
of projective measured geodesic laminations
whose support is uniquely ergodic.

There is a particular ${\cal M\cal C\cal G}$-invariant measure class on
${\cal P\cal M\cal L}$, the \emph{Lebesgue measure class}, which
can be obtained from the family of local linear structures on
${\cal P\cal M\cal L}$ defined by complete train tracks. Namely,
the transverse measure of 
every measured geodesic lamination whose support is carried by a
complete train track $\tau$ defines a nonnegative weight
function on the branches of $\tau$ which satisfies the switch
conditions. Vice versa, every nonnegative weight function on
$\tau$ satisfying the switch conditions determines a measured geodesic
lamination whose support is carried by $\tau$ \cite{PH92}.
The switch
conditions are a system of linear equations with integer
coefficients and hence the solutions of these equations have the
structure of a linear space of dimension $6g-6+2k$. The standard
Lebesgue measure on $\mathbb{R}^{6g-6+2k}$ then induces via the
thus defined coordinate system a Lebesgue measure on 
the closure of an open  
subset of the space of all measured geodesic laminations 
which is invariant under scaling. This measure 
projects to a measure class on the space of projective measured
geodesic laminations. Since the transformations of weight functions
induced by carrying maps are linear,
these locally defined measure classes do not depend
on the choice of the train track used to define the coordinates
and hence they define a ${\cal M\cal C\cal G}$-invariant measure class
on ${\cal P\cal M\cal L}$.

If $\nu$ is a measured geodesic lamination whose support is not
maximal, then this support 
is carried by a
birecurrent generic train track
which is not maximal. Thus the set
of all measured geodesic laminations whose support is not maximal
is contained in a countable union of "hyperplanes" in ${\cal P\cal
M\cal L}$ and hence has vanishing Lebesgue measure. In other
words, the Lebesgue measure class gives full measure to the set of
projective measured geodesic laminations whose support 
is maximal. Now
the support of a measured geodesic lamination is 
a disjoint union of minimal components. Such a geodesic lamination
can only be maximal if it consists of a single minimal component
and hence if it is 
minimal as well. Moreover, by a result of Masur \cite{M82}, the
Lebesgue measure class gives full mass to the set of
projective measured geodesic laminations whose support is 
\emph{uniquely ergodic}, i.e. it admits a unique transverse
measure up to scale. 
In other words, the measure-forgetting map restricted to
a Borel subset ${\cal D}$ of ${\cal P\cal M\cal L}$ of full Lebesgue
measure defines an ${\cal M\cal C\cal G}$-equivariant 
measurable bijection of ${\cal
D}$ onto a Borel subset of ${\cal C\cal L}$ and hence determines an
${\cal M\cal C\cal G}$-invariant measure class $\lambda$ on ${\cal C\cal
L}$ which we call the \emph{Lebesgue measure class}. Since by
Proposition \ref{topam} 
the action of ${\cal M\cal C\cal G}$ on ${\cal C\cal L}$
is topologically amenable, the action of ${\cal M\cal C\cal G}$ 
on ${\cal C\cal L}$ with respect to the Lebesgue measure class
is amenable \cite{AR00}. 
Then the action of ${\cal M\cal C \cal L}$ 
on ${\cal P\cal M\cal L}$ with
respect to the Lebesgue measure class is amenable as well. 
This is used to
show.

\begin{proposition}\label{strongbound}
${\cal P\cal M\cal L}$ 
equipped with the Lebesgue measure
class is a strong boundary for ${\cal M\cal C\cal G}$.
\end{proposition}
\begin{proof}
By our above discussion, we only have to show double
ergodicity for the action of ${\cal M\cal C\cal G}$ with respect to the
Lebesgue measure class. To see this let again ${\cal T}_{g,k}$ be the
Teichm\"uller space of all marked complete hyperbolic metrics on
$S$ of finite volume. The \emph{Teichm\"uller geodesic flow} is a
flow acting on the bundle $Q^1\to {\cal T}_{g,k}$ of area one
\emph{quadratic differentials} over Teichm\"uller space. 
This flow
projects to a flow on the quotient 
$Q^1/{\cal M\cal C\cal G}=Q$. 
The flow on ${\cal Q}$ preserves a Borel probability 
measure $\mu$ in the Lebesgue measure class which lifts to an 
${\cal M\cal C\cal G}$-invariant 
Lebesgue measure $\lambda_0$ on $Q^1$
\cite{M82}. 

By the
Hubbard-Masur theorem \cite{HM79}, for every point $z\in {\cal
T}_{g,k}$ and every projective measured geodesic lamination
$[\nu]$ 
on $S$ there is a unique quadratic differential $q\in {\cal
Q}^1_z$ whose \emph{vertical measured geodesic lamination} is 
contained in the class $[\nu]$ 
(where we identify a projective measured geodesic lamination
with an equivalence class of projective measured foliations on $S$
in the usual way). As $z$ varies
over ${\cal T}_{g,k}$, the set of
all these quadratic differentials defines a submanifold $W^s([\nu])$
of ${\cal Q}^1$ which projects homeomorphically onto ${\cal
T}_{g,k}$.

If the support of $[\nu]$ is uniquely ergodic and fills up $S$
then more can be said. Namely,
for every projective measured geodesic lamination
 $[\nu]^\prime\not=[\nu]$
there is a unique Teichm\"uller geodesic line in
${\cal T}_{g,k}$ generated by a quadratic differential whose
vertical measured geodesic lamination is contained in the class 
$[\nu]$ and whose horizontal measured geodesic lamination
is contained in the class
$[\nu^\prime]$. Moreover, 
if $\gamma,\gamma^\prime:\mathbb{R}\to {\cal T}_{g,k}$
are Teichm\"uller geodesics with vertical measured geodesic lamination 
contained in the class $[\nu]$, then
there is a unique number $a\in \mathbb{R}$ such that
$\lim_{t\to\infty} d(\gamma(t),\gamma^\prime(t+a))=0$
\cite{M82}. Thus if we denote by ${\cal U\cal L}$ the set of all
projective measured geodesic laminations 
on $S$ whose support is uniquely ergodic and fills up $S$, 
then the set of unit
area quadratic differentials in 
${\cal Q}^1$ with horizontal
and vertical measured geodesic lamination 
whose support is uniquely ergodic and 
fills up $S$
is (non-canonically) homeomorphic to
$({\cal U\cal L}\times {\cal U\cal L}-\Delta)\times \mathbb{R}$
(where $\Delta$ is the diagonal)
by associating to a pair of distinct points in ${\cal U\cal L}$
and the origin in $\mathbb{R}$ 
a point on the geodesic
determined by the pair of projective measured geodesic
laminations and extending this map to $({\cal U\cal L}\times
{\cal U\cal L}-\Delta)\times \mathbb{R}$ in such a way 
that the Teichm\"uller geodesic flow acts on
$({\cal U\cal L}\times {\cal U\cal L}-\Delta)\times \mathbb{R}$
via $\Phi^t([\nu],[\nu^\prime],s)=([\nu],[\nu^\prime], t+s)$. 

The Lebesgue measure 
$\lambda_0$ is locally of the form 
$d\lambda_0=d\mu^+\times d\mu^-\times dt$ where $\mu^{\pm}$
is a measure on the space of projective 
measured geodesic laminations which gives full
mass to the projective measured geodesic laminations whose
support is uniquely ergodic and fills up $S$ and which 
is uniformly exponentially expanding (contracting)
under the Teichm\"uller geodesic flow. The projection $\mu$ 
of $\lambda_0$ to the moduli space 
${\cal Q}$ of area one quadratic differentials is
ergodic \cite{M82,V86}. 
The Lebesgue measure class on ${\cal P\cal M\cal
L}\times {\cal P\cal M\cal L}$ can be obtained by
desintegration of the Lebesgue measure on 
$({\cal U\cal F}\times {\cal U\cal F}-\Delta)\times \mathbb{R}$
and therefore this measure class is doubly ergodic
under the action of ${\cal M\cal C\cal G}$ 
by the usual Hopf argument (see \cite{V86}).
This shows the proposition. \end{proof}

The \emph{curve graph} ${\cal C}(S)$ of $S$ is the metric
graph whose vertices are the free homotopy 
classes of essential simple closed curves on $S$ and where two such
curves are joined by an edge of length one if and
only if they can be realized disjointly. The 
curve graph is a hyperbolic geodesic metric graph
\cite{MM99}, and 
the mapping class group acts on ${\cal C}(S)$ 
as a group of simplicial isometries.

The \emph{Gromov boundary} $\partial {\cal C}(S)$ 
of the curve graph of $S$ can be identified with the space
of minimal geodesic laminations which fill up $S$, equipped with the
\emph{coarse Hausdorff topology}: A sequence $\{\lambda_i\}\subset
\partial {\cal C(S)}$ converges to a lamination $\lambda$ if and only
if $\lambda$ is the minimal component of every accumulation point
of this sequence with respect to the usual Hausdorff topology on 
the space of geodesic laminations \cite{Kl99,H05}.
This description of $\partial{\cal C}(S)$ is used to derive
the second corollary from the introduction.

\begin{corollary}\label{univam}
The action of ${\cal M\cal C\cal G}$ on the Gromov
boundary of the curve graph is universally
amenable.
\end{corollary}
\begin{proof}
Let ${\cal A}\subset
{\cal C\cal L}$ be the set of all complete geodesic laminations which
contain a minimal component filling up $S$. Then ${\cal A}$ is 
a Borel subset of ${\cal C\cal L}$.
Namely, for a simple closed curve $c$ on $S$ let
$B(c)$ be the set
of all complete geodesic laminations which
either contain $c$ as a minimal component or
are such that they intersect $c$ transversely in
finitely many points. Clearly $B(c)$ is a countable
union of closed subsets of ${\cal C\cal L}$ 
and therefore $B(c)$ is a Borel set. On the other hand,
the complement of ${\cal A}$
in ${\cal C\cal L}$ is the countable union of
the Borel sets $B(c)$ where 
$c$ ranges over the simple closed curves
on $S$. Then the set ${\cal A}$ 
is a Borel set as well.
Moreover, ${\cal A}$ is invariant under the action of the mapping
class group.

There is a natural continuous
${\cal M\cal C\cal G}$-equivariant map $\phi:{\cal A}\to \partial{\cal C}(S)$
which maps a lamination $\lambda\in {\cal A}$ 
to its (unique) minimal component.
By Lemma \ref{compapp}, the map $\phi$ is surjective. 
Moreover, $\phi$ is finite-to-one, which
means that the cardinality of the preimage of a point in 
$\partial {\cal C}(S)$ is bounded from 
above by a universal  number $\ell >0$.
Since the action of ${\cal M\cal C\cal G}$ on 
${\cal C\cal L}$ is topologically amenable,
the action of ${\cal M\cal C\cal G}$ on ${\cal A}$ is universally amenable
\cite{AR00}. Now the map $\phi$ is finite-to-one and
therefore by Lemma 3.6 of 
\cite{A96} the action of ${\cal M\cal C\cal G}$ 
on $\partial{\cal C}(S)$ is universally
amenable as well. This shows the corollary.
\end{proof}

Since the set of projective measured geodesic laminations
whose supports are uniquely ergodic, maximal and minimal embeds
into the Gromov boundary of the curve graph \cite{Kl99,H05},
the Lebesgue measure class on ${\cal P\cal M\cal L}$ induces
an ${\cal M\cal C\cal G}$-invariant measure class
on $\partial {\cal C}(S)$.
The following corollary is immediate from Proposition
\ref{strongbound}.

\begin{corollary}\label{strong2}
The Gromov boundary of the curve
graph equipped with the Le\-bes\-gue measure class is a strong
boundary for ${\cal M\cal C\cal G}$.
\end{corollary}

\section{Super-rigidity of cocycles}

In this section we show Theorem 2 from the
introduction. We begin with having a closer
look at the action of ${\cal M\cal C\cal G}$ on the
curve graph ${\cal C}(S)$ of $S$.

An action of a group $\Gamma$ on a Borel
space $X$ is called \emph{tame} if 
there exists a countable collection $\Theta$ of
$\Gamma$-invariant Borel subsets of $X$ which
\emph{separates orbits}. This means that  
if $X_1,X_2$ are
$\Gamma$-orbits and if $X_1\not=X_2$ then
there is some $X_0\in \Theta$ such that
$X_1\subset X_0$ and $X_1\cap X_0=\emptyset$
\cite{Z84,A96}.

Let ${\cal M\cal C}(S)$ be the countable set of multi-curves
on $S$ equipped with the discrete topology
and let ${\cal P\cal M\cal C}(S)$ be the space of probability
measures on ${\cal M\cal C}(S)$.
The mapping class group acts naturally on ${\cal M\cal C}(S)$ and
hence on ${\cal P\cal M\cal C}(S)$. 
We have.

\begin{lemma}\label{stabilizer}
The action of 
${\cal M\cal C\cal G}$ on ${\cal P\cal M}{\cal C}(S)$ is tame,
and the stabilizer of a point either is a finite
subgroup of ${\cal M\cal C\cal G}$ or it is contained
in the stabilizer of a multi-curve.
\end{lemma}
\begin{proof}
Since the set ${\cal M\cal C}(S)$ of all multi-curves on $S$ 
is countable, the action of
${\cal M\cal C\cal G}$ on ${\cal M\cal C}(S)$ is tame.
Let ${\cal Q}=\{(v,q_1,q_2)\in {\cal M\cal C}(S)
\times \mathbb{Q}\times \mathbb{Q}\mid
0\leq q_1<q_2\leq 1\}.$ 
The mapping class group acts on ${\cal Q}$, and this action is tame.
Moreover, each point $(v,q_1,q_2)\in {\cal Q}$ defines a
Borel subset $P(q)$
of ${\cal P\cal M\cal C}(S)$ by $P(q)=\{\mu\mid
\mu(v)\in [q_1,q_2]\}$.
The ${\cal M\cal C\cal G}$-orbits of these sets together with
their finite intersections define a countable collection
of ${\cal M\cal C\cal G}$-invariant 
Borel subsets of ${\cal P\cal M\cal C}(S)$ which clearly separates orbits.
This yields that
the action of ${\cal M\cal C\cal G}$ on ${\cal P\cal M}{\cal C}(S)$
is tame.  

We are left with showing that
the stabilizer of a probability measure
$\mu\in {\cal P\cal M}{\cal C}(S)$ either is a finite
subgroup of ${\cal M\cal C\cal G}$ or
is contained in the stabilizer
of a multi-curve. For this note that
a mapping class which preserves
the measure $\mu$ preserves the subset $C$
of ${\cal M\cal C}(S)$ 
on which the value of $\mu$ is maximal.
Since this set is finite and since a \emph{pseudo-Anosov} element
of ${\cal M\cal C\cal G}$ acts on the curve graph with 
unbounded orbits, the stabilizer of $\mu$ does not contain
any pseudo-Anosov element. By a
structural
result of McCarthy and Papadopoulos \cite{MP89} for
subgroups of ${\cal M\cal C\cal G}$, this implies that
either this stabilizer is finite or that it is contained
in the stabilizer of a multi-curve as claimed.
\end{proof}

The Gromov boundary $\partial {\cal C}(S)$ of
the curve graph ${\cal C}(S)$ of $S$
can be identified with the
space of all minimal geodesic laminations
on $S$ which fill up $S$, equipped with the coarse
Hausdorff topology (see \cite{Kl99,H05} and Section 8). In particular,
$\partial {\cal C}(S)$ is non-compact. 
Denote by ${\cal P}(\partial {\cal C}(S))$ the space
of probability measures on the Gromov boundary
$\partial {\cal C}(S)$ of ${\cal C}(S)$ 
and let ${\cal P}_{\geq 3}$ be the subspace of
all measures whose support contains at least
three points. The mapping class group
acts on ${\cal P}_{\geq 3}$ as a group of 
automorphisms. 
We have.

\begin{lemma}\label{tame}
The action of ${\cal M\cal C\cal G}$ on 
${\cal P}_{\geq 3}$ is tame, and the stabilizer
of a point is a finite subgroup of ${\cal M\cal C\cal G}$.
\end{lemma}
\begin{proof}
Let $T$ be the
space of triples of pairwise distinct
points in $\partial {\cal C}(S)$. Every measure
$\mu\in {\cal P}_{\geq 3}$ induces a non-trivial 
finite Borel measure $\mu^3=\mu\times\mu\times \mu$ 
on $T$. The mapping class group
acts diagonally on $T$ as a group
of homeomorphisms. 
Every point $x\in T$ 
has a neighborhood $N$ in $T$ such that
$g(N)\cap N\not=\emptyset$ only
for finitely many $g\in {\cal M\cal C\cal G}$ \cite{H08}.
Since 
$T$ is second countable, this 
immediately implies that the action of
${\cal M\cal C\cal G}$ on the space of probability
measures on $T$ is tame. But then
the action of ${\cal M\cal C\cal G}$ on 
${\cal P}_{\geq 3}$ is tame as well.

To show that the stabilizer of 
$\mu\in {\cal P}_{\geq 3}$ under the action of 
${\cal M\cal C\cal G}$ is finite, choose a 
subset $V$ of $T$ whose closure
in $\partial {\cal C}(S)\times \partial {\cal C}(S)\times
\partial {\cal C}(S)$ does not contain points in the
complement of $T$ and such that
$\mu^3(V)>3\mu^3(T)/4$. Such a set exists since
$\mu$ is Borel and $T$ is an open subset of the metrizable
space $\partial {\cal C}(S)\times \partial{\cal C}(S)\times
\partial{\cal C}(S)$.
Now if $g\in {\cal M\cal C\cal G}$ stabilizes $\mu$ 
then $gV\cap V\not=\emptyset$. On the other hand, 
there are only finitely many 
elements $g\in {\cal M\cal C\cal G}$ such that
$gV\cap V\not=\emptyset$ \cite{H08}. Therefore the
stabilizer of $\mu$ is finite.
\end{proof}

As in Section 8, let ${\cal A}\subset {\cal C\cal L}$
be the Borel subset of all complete geodesic laminations
which contain a minimal component filling up $S$.
Write ${\cal B}={\cal C\cal L}-{\cal A}$. We have.

\begin{lemma}\label{borelmap}
There is a ${\cal M\cal C\cal G}$-equivariant
Borel map $\Phi:{\cal B}\to {\cal M\cal C}(S)$.
\end{lemma}
\begin{proof} 
For every $\lambda\in {\cal B}$, the minimal 
arational components
of $\lambda$ do not fill up $S$. Thus 
a minimal arational component $\tilde \lambda$ of
$\lambda$ fills a proper connected 
subsurface $\tilde S$ of $S$ which is 
determined
by requiring that it 
contains $\tilde \lambda$ and
that every non-peripheral simple closed curve on 
$\tilde S$ intersects $\lambda$ transversely. 
The boundary of $\tilde S$ consists of a collection
of essential simple closed curves 
which do not have an essential intersection
with any minimal component of $\lambda$. 
As a consequence, the union of these 
(geodesic) boundary
curves with the closed curve components of $\lambda$
defines a non-trivial multi-curve on $S$. We obtain
a natural ${\cal M\cal C\cal G}$-equivariant
map $\Phi:{\cal B}\to 
{\cal M\cal C}(S)$ by associating to $\lambda\in {\cal B}$
this multi-curve.

We claim that the map $\Phi$ is Borel.
For this let 
$c$ be any multi-curve and let
$S_1,\dots,S_k$ be the components of $S-c$,
ordered in such a way that the pairs
of pants among these components are
precisely the surfaces $S_{\ell+1},\dots,S_k$.
Let $J$ be any (possibly empty) subset of 
$\{1,\dots,\ell\}$ and let $C_0(J)$ be the set
of all geodesic 
laminations which consist of $\vert J\vert $
minimal arational components and such that
each of these
components fills one of the subsurfaces
$S_j$ for $j\in J$. 
Let $C(J)\subset \Phi^{-1}(c)$ be the set
of all complete geodesic laminations
$\zeta\in \Phi^{-1}(c)$ 
whose minimal arational components fill the
surfaces $S_j$ for $j\in J$. Then 
$\Phi^{-1}(c)$ is a finite 
disjoint union of the sets $C(J)$ and hence $\Phi$ is
Borel if each of the sets $C(J)$ is a Borel set.

To show that this is the case, note that there
is a natural map $\Phi_0:C(J)\to C_0(J)$.
For a fixed geodesic lamination
$\nu\in C_0(J)$, 
the set $\Phi_0^{-1}(\nu)\subset C(J)\subset \Phi^{-1}(c)$
is countable.  
Namely, every lamination $\zeta\in \Phi_0^{-1}(\nu)$ 
is determined by the union of its minimal components
not contained in $\nu$ which are components of $c$ 
and by additional finitely many isolated leaves which
spiral about the minimal components. 
Each of the spiraling leaves is determined
by its closure which is the union of the leaf with 
one or two minimal components and by the homotopy class relative
to the boundary of its intersections with the components $S_i$.

On the other hand, every
complete
geodesic lamination $\xi\in C(J)$
can be obtained from a point $\zeta\in \Phi_0^{-1}(\nu)$ by replacing
some of the minimal components which fill up the
surfaces $S_j$ for $j\in J$ by another such minimal
geodesic lamination.
As a consequence, the set $C(J)$ is a  
countable union $\cup_{\zeta\in \Phi_0^{-1}(\nu)}Q(\zeta)$
where for $\zeta\in \Phi_0^{-1}(\nu)$ the set
$Q(\zeta)$ consists
of all complete geodesic laminations which can
be obtained from $\zeta$ by exchanging some of
the minimal arational components by some other
minimal arational components filling up the same subsurface.
Now to show that $\Phi^{-1}(c)$
is a Borel set, it suffices to show that
each of the sets $Q(\zeta)$ $(\zeta\in \Phi_0^{-1}(\nu))$
is Borel.

Thus let $\zeta\in \Phi_0^{-1}(\nu)$ for $\nu\in C_0(J)$
be fixed.
The set $Q^\prime(\zeta)$ 
of all complete geodesic laminations whose
restriction to $S-\cup_{j\in J}S_j$ coincides with
$\zeta$ is a Borel subset of ${\cal C\cal L}$,
and hence by the
consideration in the proof of Corollary \ref{univam},
the subset $Q(\zeta)$ of $Q^\prime(\zeta)$ 
is a Borel subset of ${\cal C\cal L}$ as well.
This shows that indeed
$\Phi^{-1}(c)$
is a Borel subset of ${\cal B}$ and completes the
proof of the lemma.
\end{proof}

Now let $n\geq 2$ and for $i\leq n$ let $G_i$ be
a locally compact 
second countable topological group. Let
$G=G_1\times \dots\times G_n$
and let $\Gamma<G$ be an irreducible 
lattice. This means that 
$\Gamma$ is a discrete subgroup of $G$ such that
the volume of $G/\Gamma$ with respect
to a Haar measure $\lambda$ on $G$ is finite
(in particular, $G$ is unimodular),
and that for every $i$ the projection of $\Gamma$ to
$G_i$ is dense.
The group $G$ acts on
$G/\Gamma$ by left translation
preserving the projection of the Haar measure.

Let $(X,\mu)$ be a standard probability space
with a measure preserving mildly mixing action of $\Gamma$. 
An ${\cal M\cal C\cal G}$-valued \emph{cocycle} for the action
of $\Gamma$ on $(X,\mu)$ 
is a measurable map 
$\alpha:\Gamma\times X\to {\cal M\cal C\cal G}$ such that
$\alpha(gh,x)=\alpha(g,hx)\alpha(h,x)$ for almost all $x\in X$,
all $g,h\in \Gamma$. The cocyle $\alpha$ is cohomologous to
a cocycle $\alpha^\prime$ if there is a measurable map
$\phi:X\to {\cal M\cal C\cal G}$ such that
$\phi(gx)\alpha(g,x)=\alpha^\prime(g,x)\phi(x)$ for every
$g\in \Gamma$ and almost every $x\in X$.
We show.

\begin{lemma}\label{cocycle}
Let $\alpha$ be an ${\cal M\cal C\cal G}$-valued
cocycle for the action of $\Gamma$ on $(X,\mu)$.
Then $\alpha$ is cohomologous to a cocycle $\alpha^\prime$
which satisfies one of the following
four possibilites.
\begin{enumerate}
\item There is some $i\leq n$ and there is
a finite subgroup $K$ of ${\cal M\cal C\cal G}$ with
normalizer $N$ such that 
$\alpha^\prime$ is induced from a continuous homomorphism
$G\to N/K$ which factors through $G_i$. 
\item The image of $\alpha^\prime$ is contained in
a finite subgroup of ${\cal M\cal C\cal G}$.
\item The image of $\alpha^\prime$ is contained in
a virtually abelian subgroup of ${\cal M\cal C\cal G}$.
\item The image of $\alpha^\prime$ is contained in
the stabilizer of a multi-curve.
\end{enumerate}
\end{lemma}
\begin{proof} Let $\Gamma$ be an irreducible
lattice in $G=G_1\times\dots\times G_n$ and
let $(X,\mu)$ be a standard probability space with a mildly
mixing measure preserving action of $\Gamma$.
Since $\Gamma<G$ is a lattice, the quotient space
$(G\times X)/\Gamma$ can be viewed as a bundle over
$G/\Gamma$ with fiber $X$.

Let $\Omega\subset G$ be a Borel fundamental domain for
the right action of $\Gamma$ on $G$. Then $\Omega\times X\subset
G\times X$ is a finite measure Borel fundamental domain 
for the action of $\Gamma$ on $G\times X$. Thus up to normalization,
the product $\lambda\times \mu$ of a Haar measure $\lambda$
on $G$ and the measure $\mu$ on $X$ projects to a probability
measure $\nu$ on $(G\times X)/\Gamma$.
The action of $G$ on $G\times X$ by left translation commutes
with the right action of $\Gamma$ and hence it projects
to an action on $(G\times X)/\Gamma$ preserving the 
measure $\nu$ (see p.75 of \cite{Z84}).

For $z\in \Omega$ and $g\in G$ 
let $\eta(g,z)\in \Gamma$ be the unique
element such that $gz\in \Omega\eta(g,z)^{-1}$.
Clearly $\eta$ is a Borel map. 
We then obtain
a $\lambda\times \nu$-measurable map $\beta:G\times (G\times X)/\Gamma\to 
{\cal M\cal C\cal G}$ by defining
$\beta(g,(z,\sigma))=\alpha(\eta(g,z),\sigma)$. 
By construction, $\beta$ 
satisfies the cocycle equation for the action of $G$ on
$(G\times X)/\Gamma$.

The group $G$ admits
a \emph{strong boundary} which is
a Borel $G$-space $B$ with an 
invariant doubly ergodic Borel measure class  
$\zeta$ and such that the action of $G$
with respect to this measure class is
amenable \cite{Ka03}. We may assume that $B$ is
also a strong boundary for $\Gamma$. 
Let ${\cal P}{\cal C\cal L}$ be the
space of all Borel probability measures on ${\cal C\cal L}$. 
Since ${\cal C\cal L}$ is a compact metric 
${\cal M\cal C\cal G}$-space, by Proposition 4.3.9 of \cite{Z84}
there is a $\beta$-equivariant Furstenberg map 
$f:(G\times X)/\Gamma\times B  
\to {\cal P}{\cal C \cal L}$.
Now the action of $G$ on $(G\times X)/\Gamma$ is ergodic
and the action of $G$ on
$B$ is doubly ergodic. Therefore
if there is a subset $D$ of $(G\times X)/\Gamma\times B$
of positive measure such that for every $x\in D$ 
the support of $f(x)$ meets the ${\cal M\cal C\cal G}$-invariant
Borel subset
${\cal A}$ of ${\cal C\cal L}$ of all complete geodesic
laminations which contain a minimal component filling up $S$, 
then this is the case for almost every
$x\in (G\times X)/\Gamma\times B$. 
Since every measure can be written as a sum of 
a measure on ${\cal A}$ and a measure on ${\cal B}={\cal C\cal L}-
{\cal A}$,
we may assume without
loss of generality that either for almost every $x$ 
the measure $f(x)$ gives full mass to the set
${\cal A}$ or
that for almost all $x$ the measure
$f(x)$ gives full mass to the set ${\cal B}$.

We distinguish now two cases.

{\sl Case 1:} For almost all $x\in (G\times X)/\Gamma\times B$, 
the measure
$f(x)$ gives full mass to the set ${\cal B}$.

By Lemma \ref{borelmap}, there is an equivariant Borel map
$\Phi:{\cal B}\to {\cal M\cal C}(S)$. The map $\Phi$ induces
a ${\cal M\cal C\cal G}$-equivariant map $\Phi_*$ of the space
of Borel probability measures on ${\cal B}$ into the space
${\cal P\cal M\cal C}(S)$ of probability measures on 
${\cal M\cal C}(S)$. 
Since by Lemma \ref{stabilizer} 
the action of ${\cal M\cal C\cal G}$ on
${\cal P\cal M\cal C}(S)$ is tame, with stabilizers which
either are finite or preserve a multi-curve,
there is 
a $\beta$-equivariant map 
$\tilde f:(G\times X)/\Gamma\times B 
\to {\cal M\cal C\cal G}/{\cal M}_0$ where
${\cal M}_0<{\cal M\cal C\cal G}$ either is a finite group or
is the stabilizer of a multi-curve. 

Since the diagonal action of $G$ on $(G\times X)/\Gamma
\times B\times B$ is also ergodic, if ${\cal M}_0$ is 
a finite group then
Lemma 5.2.10 of \cite{Z84} implies that
$\beta$ and hence $\alpha$ 
is cohomologous to a cocycle ranging in a finite 
subgroup of ${\cal M\cal C\cal G}$. Thus in this case
the second alternative in the lemma holds.

If ${\cal M}_0$ is the stabilizer of a multi-curve $c$,
then as before, $\alpha$ is cohomologous to a cocycle
with image in the stabilizer of a multi-curve, i.e. the forth
alternative in the lemma holds true. To see that this is 
indeed the
case, it is enough to inspect the proof of Lemma 5.2.10 of
\cite{Z84}. Namely, in this case the
quotient space ${\cal M\cal C\cal G}/{\cal M}_0$ can be identified with
the space of multi-curves which are
topologically equivalent to $c$, i.e. which
have the same number of components as $c$ and
which decompose the surface $S$ into connected
components of the same topological types.
Since ${\cal M\cal C\cal G}$ is countable, 
the left action of ${\cal M\cal C\cal G}$ on 
${\cal M\cal C\cal G}/{\cal M}_0\times {\cal M\cal C\cal G}/{\cal M}_0$
is tame.
For $s\in (G\times X)/\Gamma$ let $\tilde f_s:B\times B\to 
{\cal M\cal C\cal G}/{\cal M}_0\times {\cal M\cal C\cal G}/{\cal M}_0$ 
be given by $\tilde f_s(y_1,y_2)=(f(s,y_1),f(s,y_2))$.
Since the action of $G$ on $(G\times X)/\Gamma
\times B\times B$ is ergodic, 
there is an ${\cal M\cal C\cal G}$-orbit ${\cal O}$ for the
action of ${\cal M\cal C\cal G}$ on 
${\cal M\cal C\cal G}/{\cal M}_0\times {\cal M\cal C\cal G}/{\cal M}_0$ 
such that for $\nu$-almost every $s\in (G\times X)/\Gamma$ 
and for $\zeta$-almost all $y_1,y_2\in B$
we have $\tilde f_s(y_1,y_2)\in {\cal O}$.
On the other hand, since ${\cal M\cal C\cal G}/{\cal M}_0$ 
is countable, the
measure $(\tilde f_s)_*(\zeta\times \zeta)$
is a purely atomic product measure and therefore
it does not vanish on the diagonal 
of ${\cal M\cal C\cal G}/{\cal M}_0\times 
{\cal M\cal C\cal G}/{\cal M}_0$. This
diagonal is invariant under the action of 
${\cal M\cal C\cal G}$ and hence it coincides with the orbit
${\cal O}$. Therefore for $\nu$-almost
every $s\in (G\times X)/\Gamma$ the measure
$(\tilde f_s)_*\zeta$ is supported on a single point.
This implies that
$\tilde f$ induces a measurable $\beta$-equivariant
map $(G\times X)/\Gamma\to {\cal M\cal C\cal G}/{\cal M}_0$.
By the cocycle reduction
lemma (Lemma 5.2.11 of \cite{Z84}), 
this just means that the cocycle $\beta$
is cohomologous to a cocycle into 
${\cal M}_0$ and hence the same is true for $\alpha$.

{\sl Case 2:} For almost all 
$x\in (G\times X)/\Gamma\times B$ the measure
$f(x)$ gives full mass to the invariant Borel set  
${\cal A}$.

Via composing $f$ with the natural
${\cal M\cal C\cal G}$-equivariant map ${\cal A}\to
\partial {\cal C}(S)$ (compare Section 8)
we may assume that
$f$ maps $(G\times X)/\Gamma\times B$ into the space
of probability measures on $\partial {\cal C}(S)$.
By Lemma \ref{tame} the action of ${\cal M\cal C\cal G}$
on ${\cal P}_{\geq 3}$ is tame, with
finite point stabilizers. Thus
if there is a set of positive measure in
$(G\times X)/\Gamma\times B$ such that for every point
$x$ in this set the support of $f(x)$ contains 
at least 3 points, then 
we conclude as in Case 1 (and following \cite{Z84}) that
$\alpha$ is cohomologous to a cocycle ranging 
in a finite subgroup of ${\cal M\cal C\cal G}$, i.e. the second
possibility in the lemma is satisfied.

In the case that for almost all $x$ the measure
$f(x)$ is supported on two distinct points, we follow
Lemma 23 in \cite{MMS04}. Namely, let 
${\cal D}_2$ be the space of subsets of cardinality
two in $\partial{\cal C}(S)$ and let 
${\cal D}$ be the set of pairs of distinct (not necessarily disjoint)
points $(x,y)\in {\cal D}_2\times {\cal D}_2$.
The mapping class group naturally acts on ${\cal D}$ as
a group of transformations.

By hyperbolicity, there is a family of distance
functions $\delta_z$ on $\partial {\cal C}(S)$ 
$(z\in {\cal C}(S))$ of uniformly bounded diameter \cite{BH99}.
For two disjoint
elements $x=\{x_1,x_2\}$ and $y=\{y_1,y_2\}$ of ${\cal D}_2$ 
let \[\rho^\prime(x,y)=\vert \log \frac{\delta_z(x_1,y_1)
\delta_z(x_2,y_2)}{\delta_z(x_1,y_2)\delta_z(y_2,x_1)}\vert^{-1}.\]
This function extends to a continuous function
$\rho^\prime:{\cal D}\to [0,\infty)$ such that
$\rho^\prime(x,y)=0$ 
if and only if $x\cap y\not=\emptyset$ (see \cite{MMS04}.
The metric space $(\partial {\cal C}(S),\delta_z)$
is not locally compact, but the action
of ${\cal M\cal C\cal G}$ on the space
of triples of pairwise distinct points is metrically proper
\cite{H08} (this means that for every
open subset $V$ of $T$ whose distance to 
$\partial{\cal C}(S)\times \partial{\cal C}(S)\times 
\partial{\cal C}(S)-T$ with respect to a product metric
is positive,  
we have $gV\cap V\not=\emptyset$ only
for finitely many $g\in {\cal M\cal C\cal G}$).
Therefore there is a continuous function
$h:{\cal D}\to [0,\infty)$ with the following property.
\begin{enumerate}
\item For all $(x,y)\in {\cal D}$ we have
$\sum_{g\in {\cal M\cal C\cal G}}h(g^{-1}(x,y))=1.$
\item For every closed set $L\subset D$ 
whose distance to
$(\rho^\prime)^{-1}(0)$ with respect to 
the product metric is positive, the same is true
for the intersection of the support of $h$ 
with $\cup_{g\in {\cal M\cal C\cal G}}gL$.
\end{enumerate}

The function $\rho(x,y)=\sum_{g\in {\cal M\cal C\cal G}}
h(g^{-1}(x,y))\rho^\prime(g^{-1}(x,y))$ is ${\cal M\cal C\cal G}$-invariant
and continuous, and $\rho(x,y)>0$ if 
$x,y$ are disjoint. 
Using the function $\rho$
we argue as in 
Lemma 3.4 of \cite{MS04} and its proof. We conclude
that either 
the image of the cocycle $\beta$ stabilizes a
pair of distinct points in $\partial {\cal C}(S)$, or
the image of $f$ is contained in the space of
Dirac measures on $\partial {\cal C}(S)$, i.e.
$f$ is an equivariant map $(G\times X)/\Gamma \times B\to 
\partial {\cal C}(S)$. On the other hand,
if the image of $\beta$ stabilizes a pair of distinct points
in $\partial {\cal C}(S)$, then this image
is contained in a virtually abelian subgroup of 
${\cal M\cal C\cal G}$ and therefore the third alternative of 
the lemma is satisfied. 

Now we are left with the case that for almost
every $x\in (G\times X)/\Gamma\times B$ 
the measure $f(x)$ is supported in a 
single point, i.e that
$f$ is a ${\cal M\cal C\cal G}$-equivariant
map $(G\times X)/\Gamma\times B\to \partial{\cal C}(S)$.
Let $H<{\cal M\cal C\cal G}$ be the subgroup of ${\cal M\cal C\cal G}$
generated by the image of the cocycle
$\beta$ and let $\Lambda\subset
\partial {\cal C}(S)$ be the \emph{limit set} for 
the action of $H$ on ${\cal C}(S)$, i.e.
$\Lambda$ is the set of accumulation points of 
an $H$-orbit on ${\cal C}(S)$. 
If this limit set consists of at most
two points, then either $H$ is virtually abelian
and the third alternative in the lemma holds true,
or an $H$-orbit in ${\cal C}(S)$ is bounded. 
However, if an $H$-orbit on ${\cal C}(S)$ is bounded then
the structure result of \cite{MP89}
implies that either $H$ is finite or $H$  
stabilizes a multi-curve (compare
the proof of Lemma \ref{stabilizer}). 
Thus for the
purpose of the lemma, we may assume that
$\Lambda$ contains more than two points; then
$\Lambda$ is
a closed $H$-invariant subset of $\partial {\cal C}(S)$

Let $T$ be the space of triples of pairwise 
distinct points in $\Lambda$.
We showed in 
\cite{H08} that there is a continuous non-trivial
cocycle
$\rho:T\to L^2(H)$ for the action of 
$H<{\cal M\cal C\cal G}$ on $T$. Via $\beta$, 
this cocycle
then induces a non-trivial measurable 
$L^2(H)$-valued cocycle
on $(G\times X)/\Gamma\times B$. Since the action of
$G$ on $B$ is doubly ergodic, such a cocycle 
defines a non-trivial cohomology class
for $G$ with coefficients in the separable 
Hilbert space
$L^{[2]}((G\times X)/\Gamma,L^2(H))$ of
all measurable maps $(G\times X)/\Gamma\to L^2(H)$
with the additional property that 
for each such map $\alpha$, the function
$x\to \Vert\alpha(x)\Vert$ is square integrable
on $(G\times X)/\Gamma$.

For every $i\leq n$ define $G_i^\prime=\prod_{j\not=i}G_j$.
By assumption, the projection of $\Gamma$
to $G_i=G_i^\prime\backslash G$ is dense and 
therefore by Lemma 2.2.13 of \cite{Z84}, the
action of $\Gamma$ on $G_i$ by right translation is ergodic.
Then by Moore's ergodicity theorem, the action of
$G_i^\prime$ on $G/\Gamma$ is ergodic as well. 
Since by assumption the action of $\Gamma$ on 
$(X,\mu)$ is mildly mixing, the results of 
\cite{SW82} imply that  
the action of $G_i$ on $((G\times X)/\Gamma,\nu)$ 
is ergodic.
As in 
\cite{MMS04}, we conclude from the results
of Burger and Monod \cite{BM99,BM02} that there is
some $i\leq n$ and there
is an equivariant map $(G\times X)/\Gamma\to 
L^2(H)$ for the restriction of 
$\beta$ to $G_i^\prime\times (G\times X)/\Gamma$.
Since the action of $G_i^\prime$ on $(G\times X)/\Gamma$ is ergodic,
the restriction of $\beta$ to
$G_i^\prime$ is equivalent to a cocycle into 
the stabilizer of a non-zero element of 
$L^2(H)$ (see Lemma 5.2.11 of \cite{Z84}). 
Since the action of $H$ on itself
by left translation is simply transitive, this
stabilizer is a
finite subgroup of $H$.

We now follow the proof of Theorem 1.2 of
\cite{MS04} and find a minimal such finite subgroup $K$ of 
${\cal M\cal C\cal G}$.
The cocycle $\beta$ is cohomologous to a cocycle
$\beta^\prime$ into the normalizer $N$ of $K$ in ${\cal M\cal C\cal G}$,
and the same is true for the cocycle $\alpha$.
Proposition 3.7 of \cite{MS04} then shows that
$\beta^\prime$ is induced from a continuous homomorphism
$G\to N/K$ which factors through $G_i$. In other words,
the first alternative in the statement of the lemma is satisfied.
This completes the proof of the lemma.
\end{proof}

Theorem \ref{theorem2} from the introduction follows
from Lemma \ref{cocycle} and an analysis of cocycles
$\alpha$ which are cohomologous 
to a cocycle, again
denoted by $\alpha$, 
with values in the stabilizer $H$ of 
a multi-curve $c$. Now if $c$ has $d\leq 3g-3+m$ components,
then the group $H$ is a direct product of 
$\mathbb{Z}^d$ with the mapping 
class group ${\cal M}_0$ of the (possibly disconnected)
surface $S_0$ which we obtain from $S-c$ by
pinching each boundary circle to a point, 
and where the group $\mathbb{Z}^d$ 
is generated by the Dehn twists about the
components of $c$. Note that ${\cal M}_0$ is a finite
extension of the direct product of the mapping
class groups of the connected components of $S_0$.
Let $p:H\to {\cal M}_0$ be
the natural projection and
let ${\cal C}_0$
be the space of all complete geodesic laminations
on $S_0$. Then ${\cal C}_0$
is a compact ${\cal M}_0$-space and hence
we can apply our above procedure to the cocycle
$p\circ \alpha$, the compact space ${\cal C}_0$ and 
the group ${\cal M}_0$. After at most $3g-3+k$ steps
we conclude that indeed the cocycle $\alpha$
is cohomologous to a cocycle with values in subgroup
of ${\cal M\cal C\cal G}$ which is a finite
extension of a group of the form $\mathbb{Z}^\ell\times 
N/K$ where $\ell\leq 3g-3+k$, where
$K<{\cal M\cal C\cal G}$ is finite, $N$ is the normalizer
of $K$ in ${\cal M\cal C\cal G}$ and such that
there is a continuous homomorphism $G\to N/K$.
This completes the proof of Theorem \ref{theorem2}.

\section{Appendix: Splitting, Shifting  and Carrying}

In this appendix we establish some technical properties
of complete 
train tracks on $S$. This appendix only uses
results from the literature which are summarized in Section 2,
and is 
independent from the rest of the paper. Its main goal
is the proof of Proposition \ref{carrynear} 
which is
used in an essential way in Section 6.

By the definition given in Section 2, 
a complete train track on the non-exeptional surface $S$
of finite type is a train track which is maximal, generic
and birecurrent. By Lemma \ref{recchar}, 
such a train track $\tau$ carries 
a complete geodesic lamination $\lambda$. It 
can be split at any large branch 
to a train track $\sigma$ which carries $\lambda$ as well, and
$\sigma$ is complete and is carried by $\tau$.
Moreover, the choice
of a right or left split at $e$ is uniquely determined
by $\lambda$. On the other hand, $\tau$ can not be
split to \emph{every} complete train track $\eta$ which is carried
by $\tau$. Namely, to pass from $\tau$ to $\eta$ it
may be necessary to do some shifting moves as well (see
Chapter 2 of \cite{PH92} for a detailed discussion
of this fact).
Shifting involves making choices, and 
these choices are hard to control combinatorially.

We overcome this difficulty by obtaining 
a better quantitative understanding of the relation between
splitting, shifting  and carrying of complete train
tracks. This is done by putting a train track $\sigma$ which
is carried by a train track $\tau$ in a standard position with
respect to $\tau$ and use this standard position to define
a numerical invariant which can be controlled under splitting moves.
Part of the material presented here is motivated by the
results in Section 2.3 and Section 2.4 of \cite{PH92}. 
Recall from Section 2 the
definition of a foliated neighborhood
$A$ of a generic train track $\tau$ \cite{PH92}.

\begin{definition}\label{generalposition} 
Let $\tau$ be a generic train track with foliated neighborhood $A$
and collapsing map $F:A\to \tau$. 
A generic train track $\sigma$ 
is \emph{in general position with respect to $A$}
if the following three conditions are satisfied.
\begin{enumerate}
\item $\sigma$ is contained in the interior of $A$ and is
transverse to the ties.
\item No
switch of $\sigma$ is mapped by $F$ to a switch of $\tau$.
\item 
For every $x\in\tau$ the tie $F^{-1}(x)$ contains
at most one switch of $\sigma$. 
\end{enumerate}
\end{definition}

Note that by definition, a train track in general position
with respect 
to a foliated neighborhood 
$A$ of $\tau$ 
is always carried by $\tau$, with the restriction of the
collapsing map $F:A\to \tau$ as a carrying map. 
Moreover, every generic train track $\sigma$ which
is carried by $\tau$ can be isotoped to a train track
in general position with respect to $A$.

Let $\sigma$ be a generic train track in general position
with respect to the standard neighborhood $A$ of the 
complete train track $\tau$, with collapsing
map $F:A\to\tau$. Let $b$ be a branch of 
$\tau$ which is incident and large on a switch $v$
of $\tau$. A \emph{cutting arc}
for $\sigma$ and $v$ is an embedded arc $\gamma:[0,d]\to F^{-1}(b)$
which is transverse to the ties of $A$, which is disjoint from 
$\sigma$ except possibly at its endpoint and such that
the length of $F(\gamma)\subset b$ is maximal among all arcs
with these properties. The maximality condition implies that
either $F(\gamma[0,d])=b$, i.e. that $\gamma$ crosses through
the foliated rectangle $F^{-1}(b)$, or that $\gamma(b)$ is a switch
of $\sigma$ contained in the interior of $F^{-1}(b)$, and the
component of $S-\sigma$ which contains $v$ has
a cusp at $\gamma(d)$.

We use cutting arcs to investigate whether or not $\sigma$ is carried
by a split of $\tau$ at some large branch $e$.

\begin{lemma}\label{cuttingarc}
Let $\sigma$ be in general position with respect to the foliated
neighborhood $A$ of $\tau$, with collapsing map
$F:A\to \tau$. Let $e$ be a large branch which is incident on 
the switches $v,v^\prime$ of $\tau$. 
Let $\gamma:[0,d]\to F^{-1}(e),\gamma^\prime:[0,d]\to F^{-1}(e)$
be cutting arcs for $\sigma$ and $v,v^\prime$. 
If $\sigma$ is not carried
by a split of $\tau$ at $e$ then $\gamma(d),\gamma^\prime(d)$ 
are switches of 
$\sigma$ contained in the interior of $F^{-1}(e)$, and there
is a trainpath $\rho:[0,m]\to \sigma\cap F^{-1}(e)$ 
connecting $\rho(0)=\gamma(d)$ to $\rho(m)=\gamma^\prime(d)$.
\end{lemma}
\begin{proof} Let $e$ be a large branch of $\tau$ incident on the
switches $v$ and $v^\prime$. Let $\sigma$ be in general position
with respect to the foliated
neighborhood $A$ of $\tau$ with collapsing map $F:A\to \tau$.
Assume that $\sigma$ is not
carried by any split of $\tau$ at $e$
and let $\gamma:[0,d]\to F^{-1}(e),\gamma^\prime:[0,d]\to
F^{-1}(e)$ be cutting
arcs for $v,v^\prime$ and $\sigma$. 

We show first that $\gamma(d)$ is 
a switch of $\sigma$ contained in the interior of $F^{-1}(e)$.
For this we argue by contradiction and we
assume that $\gamma(d)\in F^{-1}(v^\prime)$. 
Since $\sigma$ is contained in
the interior of $A$ and $v^\prime$ is contained in the boundary of
$A$, we may assume that $\gamma(d)$ is contained in
$F^{-1}(v^\prime)-v^\prime -\sigma$. 
Then $\gamma(d)$ is contained in the boundary of a unique
foliated rectangle $R\subset A$ whose image under the collapsing
map $F$ is a branch of $\tau$ which is incident and small at 
the switch $v^\prime$.
Extend $\gamma$ a bit beyond $\gamma(d)$ 
to an arc $\tilde\gamma$ transverse to the ties of $A$ and disjoint
from $\sigma$ whose endpoint is contained in the interior of $R$.
The train track obtained by 
cutting the foliated neighborhood
$A$ open along $\tilde\gamma$ as shown in Figure D
and by collapsing the ties of the induced 
foliation of $A-\tilde\gamma$ carries $\sigma$ and is a split
of $\tau$ at $e$. This means that $\sigma$ is carried by
a split of $\tau$ at $e$ which is a contradiction.
As a consequence, $w=\gamma(d)$ and $w^\prime=
\gamma^\prime(d)$ are switches
of $\sigma$ contained in the interior of $F^{-1}(e)$ which are cusps
for the complementary components of $\sigma$ containing $v,v^\prime$.

Let $\nu:[0,\infty)\to \sigma$ be a trainpath
on $\sigma$ issuing from $\nu(0)=w$ such that
$\nu[0,1/2]$ is the half-branch of $\sigma$
which is incident on $w$ and large at $w$. Since $\sigma$
is transverse to the ties of $A$, there is a smallest
number $t>0$ such that $F(\nu(t))=v^\prime$. The arc
$\nu[0,t]$ is embedded in $\sigma$. 
Now if every such embedded arc
of class $C^1$ in $F^{-1}(e)\cap \sigma$
connects $w$ to the same component
of $F^{-1}(v^\prime)-v^\prime$ then $\sigma$ is carried
by a split of $\tau$ at $e$. Thus 
there are two embedded arcs $\nu,\nu^\prime$ of class $C^1$ in
$\sigma\cap F^{-1}(e-F(\gamma[0,d)))$
whose endpoints are contained in the two different components
of $F^{-1}(v^\prime)-v^\prime$.

Let $C$ be the 
connected component of $F^{-1}(e)-(\nu\cup\nu^\prime)$ which
contains the
switch $v^\prime$ of $\tau$. 
The closure of $C$ does not intersect $F^{-1}(v)$, and it contains
the switch $w^\prime$ of $\sigma$.
Every trainpath on $\sigma$
issuing from $w^\prime$ which intersects $F^{-1}(v)$ necessarily
intersects one of the arcs $\nu$ or $\nu^\prime$. 
Since the restriction of the differential $dF$ of $F$ to
$\sigma$ vanishes nowhere, this implies
that there is an embedded trainpath 
$\rho:[0,m]\to \sigma$ which is entirely
contained in $F^{-1}(e)$ and connects $w=\rho(0)$ to
$w^{\prime}=\rho(m)$. This shows the lemma.
\end{proof}

Next we show that a train track $\sigma\in {\cal V}({\cal T\cal T})$
which is carried by a train track $\tau\in {\cal V}({\cal T\cal T})$
and is not carried by any split of $\tau$ can be shifted and
isotoped to a train track $\sigma^\prime$ which is in standard
position relative to $\tau$.
For this we have to introduce
some more terminology.

Let $\sigma$ be any generic train track on $S$ and let
$\rho:[0,n]\to \sigma$ be an embedded trainpath. For $i\in
\{1,\dots,n-1\}$ the point $\rho(i)$ is a switch of $\sigma$
contained in the interior of $\rho[0,n]$. Since every switch
is trivalent, there is a unique
branch of $\tau$ which is incident on $\rho(i)$ and not contained
in $\rho[0,n]$. We call $\rho(i)$ a \emph{right switch} if this
branch lies to the right of $\rho$ with respect to the orientation
of $S$ and the orientation of $\rho$; otherwise $\rho(i)$ is
called a \emph{left switch}. Moreover, the switch $\rho(i)$ is
called \emph{incoming} if the half-branch $\rho[i-1/2,i]$ is small
at $\rho(i)$, and we call the switch $\rho(i)$ \emph{outgoing}
otherwise.

Define a \emph{special trainpath} on a 
generic train track $\sigma$ to be
a trainpath $\rho:[0,2k-1]\to \sigma$ of length $2k-1$ for some
$k\geq 1$ with the following properties.
\begin{enumerate}
\item
$\rho[0,2k-1]$ is embedded in $\sigma$.
\item
For each $j\leq k-1$ the branch $\rho[2j,2j+1]$ is large and the
branch $\rho[2j+1,2j+2]$ is small.
\item
With respect to the orientation of $S$ and the orientation of
$\rho$, right and left switches in $\rho[1,2k-2]$ alternate.
\end{enumerate}
The left part of Figure E shows a special trainpath of length 5.
\begin{figure}[hb]
\includegraphics{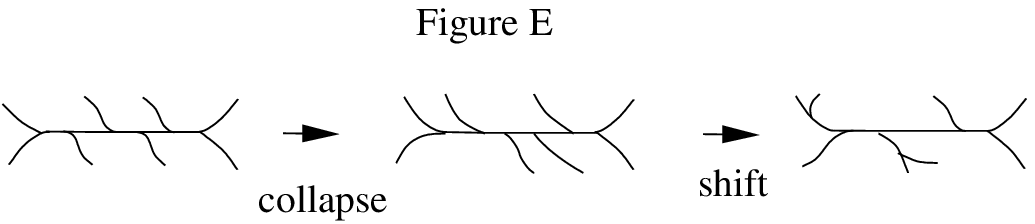}
\end{figure}

Let again 
$\sigma$ be a generic train track and let $\rho:[0,k]\to \sigma$ be any
embedded trainpath on $\sigma$. Define a \emph{standard
neighborhood} of $\rho$ in $\sigma$ to be a closed connected
neighborhood of $\rho[0,k]$ in the graph $\sigma$ which does not
contain a switch distinct from the switches contained in
$\rho[0,k]$. As an example, a standard neighborhood of a single
branch $e$ of $\sigma$ incident on two distinct switches
is the union of $e$ with four neighboring
half-branches.

Let 
$\tau\in {\cal V}({\cal T\cal T})$ be a complete
train track and let 
$A$ be a foliated neighborhood of $\tau$ with collapsing map 
$F:A\to \tau$.
Let $\sigma\subset A$ be a complete train track which is in
general position with respect to $A$.
Then the restriction of the collapsing map $F$ to $\sigma$ is
surjective. Thus if we denote for a branch $b$ of $\tau$ by
$\nu(b,\sigma)$ the minimal cardinality of $F^{-1}(x)\cap \sigma$
where $x$ varies through the points in $b$ then $\nu(b,\sigma)>0$ for
every branch $b$. Since no switch of $\sigma$ is mapped by $F$ to
a switch of $\tau$, there is a point $x$ contained in the interior
of $b$ for which the cardinality of $F^{-1}(x)\cap \sigma$ equals
$\nu(b,\sigma)$. Moreover, since any tie of $A$ contains
at most one switch of $\sigma$ by assumption, we may assume
that $F^{-1}(x)$ does not contain any switch of $\sigma$.
Note that if we deform $\sigma$ with an isotopy
of train tracks in general position relative to $A$, then 
the values $\nu(b,\sigma)$ do not change along the isotopy.
However, this may not be the case for general
deformations of $\sigma$ transverse to the ties of $A$.

The following technical lemma shows that a train track $\sigma\in
{\cal V}({\cal T\cal T})$ 
which is carried by a train track $\tau\in {\cal V}({\cal
T\cal T})$ and is not carried by any split of $\tau$ can be shifted
and isotoped to a train track $\sigma^\prime$ which is in standard
position relative to $\tau$.

\begin{lemma}\label{specialtrain}
Let $\tau \in {\cal V}({\cal T\cal T})$ 
and let $A$ be
a foliated neighborhood of $\tau$ with collapsing map
$F:A\to\tau$. Let $\sigma\in 
{\cal V}({\cal T\cal T})$ be in general position with
respect to $A$, let $e$ be a large
branch of $\tau$ and let $V\subset U$ be standard 
neighborhoods of $e$ in $\tau$
so that $V$ is contained in the interior of $U$.
If $\sigma$ is not carried by a split of $\tau$ at $e$,
then $\sigma$ can be shifted and isotoped to a train track
$\sigma^\prime$ with the following properties.
\begin{enumerate}
\item $\sigma^\prime$ is in general position with respect to $A$,
and we have
$\sigma^\prime\cap F^{-1}(\tau-U)= \sigma\cap
F^{-1}(\tau-U)$.
\item $F^{-1}(V)\cap \sigma^\prime$ is the
disjoint union of a standard neighborhood $V^\prime$ of a special
trainpath $\rho$ on $\sigma^\prime$ and a (possibly empty)
collection of simple arcs not containing any switch of
$\sigma^\prime$. The set $V^\prime$ is mapped by $F$ onto $V$.
\item
If $\nu(e,\sigma)=1$ then
$F^{-1}(V)\cap \sigma^\prime$ is a standard neighborhood of a
single large branch in $\sigma^\prime$.
\end{enumerate}
\end{lemma}
\begin{proof}
Let $\tau$ be a complete train track with foliated neighborhood
$A$ and collapsing map $F:A\to \tau$.
Let $\sigma\subset A$ be a complete
train track in general position with respect to $A$. 
Then $\sigma$ is carried by $\tau$ with 
the restriction of $F$ as
a carrying map. 
Let $e$ be a large branch of $\tau$ and
assume that $\sigma$ is not carried by any split 
of $\tau$ at $e$.

Let $v,v^\prime$ be the switches 
of $\tau$ on which the large branch $e$ is
incident. The points $v,v^\prime$ are contained in complementary
regions of $\sigma$. 
Let $\gamma:[0,d]\to F^{-1}(e)$ and 
$\gamma^\prime:[0,d^\prime]\to F^{-1}(e)$ 
be cutting arcs for $\sigma$ and $v,v^\prime$.
By Lemma \ref{cuttingarc} 
the points $w=\gamma(d)$ and
$w^\prime=\gamma(d^\prime)$ are switches of $\sigma$ contained
in the interior of $F^{-1}(e)$. Moreover, there is a trainpath
$\rho:[0,m]\to \sigma\cap F^{-1}(e)$ which connects
$\rho(0)=\gamma(d)$ to $\rho(m)=\gamma^\prime(d)$.
The half-branch $\rho[0,1/2]$ is large at $\rho(0)$, and
$\rho[m-1/2,m]$ is large at $\rho(m)$. Since
$\sigma$ does not have embedded bigons, $\rho$ is unique.
We call $\rho$ the \emph{cutting connector} for $\sigma$ and $e$.

Recall the definition of right and left and incoming and outgoing
switches along the trainpath $\rho$. Define the switch $\rho(0)$ to
be incoming and the switch $\rho(m)$ to be outgoing. 
Let $j\in
\{1,\dots,m\}$ be the smallest number with the property 
that $\rho(j)$
is outgoing. We claim that if $\rho(j)$ 
is right outgoing, then none of the switches 
$\rho(k)$ for $k>j$ is right
incoming. 
For this we assume to the contrary that for some $k>j$
the switch $\rho(k)$ is right incoming. 
Now every subarc of a trainpath on
$\sigma$ which is entirely contained in $F^{-1}(e)$ is mapped
$C^1$-diffeomorphically onto a subarc of $e$ by the map $F$. Let
$b_j,b_k$ be the branches of $\sigma$ which are incident on
$\rho(j),\rho(k)$ and whose 
interiors are disjoint from $\rho$. Then there
is a subarc of a trainpath on $\sigma$ which begins at
$\rho(k)$, connects $\rho(k)$ to a point in $F^{-1}(F\rho(0))$ 
and whose initial
segment is contained in the branch 
$b_k$. Since $\sigma$ does not have
embedded bigons, the interior of this arc is entirely contained
in the interior of the component $R$ of $F^{-1}(F(\rho[0,m]))-\rho[0,m]$
to the right of $\rho$ with respect to the orientation of $\rho$
and the orientation of $S$. 
Similarly, there is an embedded
subarc of a trainpath on $\sigma$ which begins at $\rho(j)$, connects
$\rho(j)$ to a point in $F^{-1}(F\rho(m))$ and whose initial segment
is contained in the branch $b_j$. The interior
of this arc is contained in $R$ as well.
Therefore these arcs have to intersect in the interior of $R$. 
However, $\sigma$ does not have embedded
bigons and hence this is impossible. Note also that if there is a point
$x\in e$ such that $F^{-1}(x)\cap \sigma$ consists of a single point,
then $x\in F(\rho[j-1,j])$ and none of the switches
$\rho(i)$ for $i´\geq j$ is incoming.

We can modify the train track $\sigma$ with a sequence
of shifts and isotopies in $m-1$ steps as follows. 
Let again $j\leq m$ be the smallest number such that the
switch $\rho(j)$ is outgoing. 
The branches incident on the switches
$\rho(1),\dots,\rho(j-1)$ can successively be moved backward along
$\rho$ with a sequence of shifts past the switch $w$
to a switch contained in $F^{-1}(e-F\rho[0,m])$ as shown on 
the right hand side of Figure E.
This operation can be done in such a way that
the resulting train track $\tilde \sigma$ is 
in general position with respect to $A$, that
$\nu(e,\sigma)=1$ if and only if $\nu(e,\tilde \sigma)=1$ and
that $\tilde \sigma$ coincides with $\sigma$ in $A-F^{-1}(e)$.
The length of the cutting connector $\tilde \rho$ 
for $v$ and $v^\prime$ in $\tilde \sigma$
equals $m-j+1$, and $\tilde \rho[k,k+1]=\rho[k+j-1,k+j]$
for $k\geq 1$.
The branch $\tilde \rho[0,1]$ in $\tilde \sigma$
which is incident and large at $w$ is a large branch. 
Moreover, if $\nu(e,\sigma)=1$ then there is a point
$x\in F(\tilde \rho[0,1])$ such that the cardinality
of $\tilde \sigma\cap F^{-1}(x)$ equals one.

If each of the switches $\tilde \rho(i)$ for $i\geq 1$
is outgoing then we can modify $\tilde \sigma$ with a
sequence of shifts so that the cutting connector in the
modified train track is a single large branch. Note that
this is the case if $\nu(e,\sigma)=1$. Otherwise assume
that $\tilde \rho(1)$ is right outgoing and let 
$k>1$ be the smallest number such that $\tilde \rho(k)$ 
is left incoming. Then each of the switches 
$\rho(i)$ for $1\leq i\leq k-1$ is right outgoing and hence
$\tilde \sigma$ can be modified with a sequence of shifts
in such a way that the 
path $\tilde \rho[1,k]$ corresponds to a single small branch
in the cutting connector for the modified train track.
Inductively, this just means that with a sequence of shifts
and an isotopy we can modify $\sigma$ to a train track $\sigma_1$
which is in general position with respect to $e$, which 
coincides with $\sigma$ in $A-F^{-1}(e)$ and such that
the cutting connector $\rho_1$ for $v$ and $v^\prime$ in $\sigma_1$
is a special trainpath. If $\nu(e,\sigma)=1$ then
$\rho_1$ consists of a single large branch, and there is 
a point $z$ in the interior of this branch such that
the cardinality of $F^{-1}(F(z))\cap \sigma_1$ equals one.

Since $\sigma_1$ is mapped by $F$ onto $\tau$,
with respect to the orientation of $S$ and the orientation of
a cutting arc $\gamma$ for $\sigma_1$ and $v$,
for each point $u\in F(\gamma)$ there is a unique point
$r(u)\in \sigma_1\cap F^{-1}(u)$ (or a unique point $\ell(u)\in
\sigma_1\cap F^{-1}(u)$) 
which is to the right (or to the left) of $\gamma$
and which is closest to $\gamma\cap F^{-1}(u)$ with this property
as measured along the arc $F^{-1}(u)$. Then $u\to r(u)$ 
and $u\to \ell(u)$ are
embedded arcs of class $C^1$ in $\sigma_1$ which are disjoint from
$\gamma$ except at their endpoints. 

Let $V\subset U$ be standard neighborhoods of $e$ so that the closure
of $V$ is contained in the interior of $U$.
We can isotope $\sigma_1$ by sliding the switches contained in the 
interior of the arcs $r$ and $\ell$ backward along
$r,\ell$ outside of $F^{-1}(V)$. This isotopy can
be chosen to be supported in 
$F^{-1}(U)-F^{-1}(e-F(\gamma[0,d]))$ and to be transverse
to the ties of $A$. Moreover, we may assume that the image
of $\sigma_1$ under this isotopy is in general position with
respect to $A$. Since
for every $u\in e-F(\gamma[0,d))$
the isotopy does not
change the cardinality of $F^{-1}(u)\cap \sigma_1$, if
$\nu(e,\sigma_1)=1$ then the same is true for the image
of $\sigma_1$ under this isotopy.

As a consequence, after modifying $\sigma_1$ with an isotopy
supported in $F^{-1}(U)$ we obtain a train track $\sigma_2$
whose intersection with $F^{-1}(V)$ contains
a connected component which is a standard 
neighborhood $V^\prime$ of
a special trainpath $\rho_2$, and $\rho_2$ is the cutting connector
for $v$ and $v^\prime$ in $\sigma_2$.

Now if there is a connected component $\beta\not=V^\prime$
of $F^{-1}(V)\cap \sigma_2$, then $\nu(e,\sigma)>1$ and
$\beta\cap F^{-1}(e)$  is contained in a connected
component of $F^{-1}(e)-\rho_2$
whose boundary intersects both
$F^{-1}(v)$ and $F^{-1}(v^\prime)$. 
It is then immediate that we can deform $\sigma_2$ with an isotopy 
supported in $F^{-1}(U)$ to a 
train track $\sigma^\prime$ which is in general
position with respect to $A$ and has the property that every connected
component of $\sigma^\prime\cap F^{-1}(e)$ which does not coincide
with the standard neighborhood $V^\prime$ of $\rho_2$ is a single
arc not containing any switches. 
In other words, $\sigma^\prime$ satisfies the requirements in the lemma.
\end{proof}

The following observation is 
a partial converse of Lemma \ref{specialtrain}. 
It is
used to determine whether or not a
given complete train track $\sigma\in {\cal V}({\cal T\cal T})$ 
which is carried by a complete
train track $\tau\in {\cal V}({\cal T\cal T})$ is carried by a split of
$\tau$ at a given large branch $e$.

\begin{lemma}\label{foliated}
Let $A$ be a foliated neighborhood of a
complete train track $\tau$ with collapsing map $F:A\to \tau$. Let
$\sigma\in {\cal V}({\cal T\cal T})$ be in general position
with respect to $A$.
Let $e$ be a large branch of $\tau$; if
$\nu(e,\sigma)=1$ then $\sigma$ is not carried by any split of
$\tau$ at $e$. \end{lemma}
\begin{proof}
Let $\tau\in {\cal V}({\cal T\cal T})$ 
be a complete train
track with foliated neighborhood $A$ and collapsing map $F:A\to
\tau$. Let $\sigma\in {\cal V}({\cal T\cal T})$ 
be in general position with respect to $A$.
Let $e$ be a large branch of $\tau$ and assume that
$\nu(e,\sigma)=1$. Then 
there is a point $x$ in the interior of $e$ with the property that
$F^{-1}(x)\cap \sigma$ is a single point $z$ 
contained in the interior of a
branch of $\sigma$.

Let $v$ be switch on which the branch $e$ is incident and
let $\gamma:[0,d]\to F^{-1}(e)$ be a cutting arc for $\sigma$ and $v$.
Since $\sigma$ is complete, 
the restriction of the map $F$ to $\sigma$ is surjective.
Thus there is an embedded arc 
$\nu:[0,1]\to\sigma\cap F^{-1}(e)$ 
(or $\nu^\prime:[0,1]\to F^{-1}(e)$)
which is mapped by $F$ onto $e$ and 
which begins at a point $\nu(0)\in F^{-1}(v)$ (or 
$\nu^\prime(0)\in F^{-1}(v)$)
to the right (or to the left) 
of $v$ with respect to the orientation of $\gamma$
and the orientation of $S$. 
Since $F^{-1}(x)\cap \sigma$ consists of the single point $z$, the arcs 
$\nu,\nu^\prime$ have to 
intersect. Thus $\gamma(d)$ is a switch of $\sigma$ contained
in the interior of $e$, and  
$x\in e-F(\gamma[0,d])$. Similarly,
a cutting arc $\gamma^\prime:[0,d]\to F^{-1}(e)$ for $\sigma$ and 
$v^\prime$ ends at 
a switch of $\sigma$ contained in the interior of $F^{-1}(e)$ 
and such that $x\in e-F(\gamma^\prime[0,d])$.

Now every embedded arc $\alpha:[0,1]\to \sigma$ of class $C^1$ 
with $\alpha(0)=\gamma(d)$ and $F\alpha[0,1]=e-F\gamma[0,1)$ 
passes through $z$, and the same holds true for 
an embedded arc $\alpha^\prime:[0,1]\to \sigma$ of 
class $C^1$ with $\alpha^\prime(0)=\gamma^\prime(d)$ and 
$F\alpha^\prime[0,1]=e-\gamma^\prime[0,1)$. 
But this implies that there is an embedded trainpath in 
$\sigma$ which is contained in $F^{-1}(e)$ and connects
$\gamma(d)$ to $\gamma^\prime(d)$. Then 
$\sigma$ is not carried by a split of 
$\tau$ at $e$.
\end{proof}

Call two train tracks $\tau,\sigma\in {\cal V}({\cal T\cal T})$
\emph{shift equivalent} if $\sigma$ can be obtained
from $\tau$ by a sequence of shifts. This clearly
defines an equivalence relation on ${\cal V}({\cal T\cal T})$.
The number of points in each equivalence class is bounded
from above by a universal constant only depending
on the topological type of the surface $S$.
The following lemma characterizes complete train tracks 
$\sigma\in {\cal V}({\cal T\cal T})$ 
which are shift equivalent to a complete train track 
$\tau\in {\cal V}({\cal T \cal T})$.

\begin{lemma}\label{shiftequi}
Let $\tau\in {\cal V}({\cal T\cal T})$ 
and let $A$ be
a foliated neighborhood of $\tau$ with collapsing map $F:A\to
\tau$. Let $\sigma \in {\cal V}({\cal T\cal T})$ 
be in general position with respect to $A$.
If $\nu(e,\sigma)=1$ for every large branch $e$ of $\tau$, then
$\sigma$ is shift equivalent to $\tau$.\end{lemma}
\begin{proof}
Let $\tau\in {\cal V}({\cal T\cal T})$ with 
foliated neighborhood $A$ and collapsing map $F:A\to
\tau$. Let $\sigma\in 
{\cal V}({\cal T\cal T})$ be in general position with respect
to the $A$. Assume that
$\nu(e,\sigma)=1$ for every large branch of $\tau$; we have to
show that $\sigma$ is shift equivalent to $\tau$.

By Lemma \ref{foliated}, 
$\sigma$ is not carried by any split
of $\tau$. 
Lemma \ref{specialtrain} 
then shows that 
after modifying $\sigma$ by a sequence
of shifts and an isotopy we may assume that for every large branch
$e$ of $\tau$ there is a standard neighborhood $U$ of $e$ such
that $F^{-1}(U)\cap \sigma$ is a standard neighborhood of a large
branch in $\sigma$.

Following \cite{PH92}, define a 
\emph{large one-way trainpath} on $\tau$
to be a trainpath $\rho:[0,m]\to \tau$ such that $\rho[0,1/2]$ and
$\rho[m-1/2,m]$ are large half-branches and that for every $i<m$ the
half-branch $\rho[i-1/2,i]$ is small. 
Every switch $v$ in $\tau$ is the starting point of a
unique large one-way trainpath, and this path is embedded. Define
the \emph{height} $h(v)$ of a switch $v$ of $\tau$ to be the
length of the large one-way trainpath starting at $v$. The
switches of height 1 are precisely the switches on which a large
branch is incident.

Write $\sigma_1=\sigma$.
We construct inductively for each $m\geq 1$ a 
complete train track
$\sigma_m$ which can be connected to $\sigma_{m-1}$ by a sequence
of shifts and an isotopy and with the following properties.
\begin{enumerate}
\item $\sigma_m$ is in general position with
respect to $A$.
\item Every large one-way trainpath $\rho$ on $\tau$ of
length at most $m$ admits a standard neighborhood $U$ with the
property that $F^{-1}(U)$ is a standard neighborhood $V$ of a
large one-way trainpath $b(\rho)$ in $\sigma_m$ which is isotopic to
$\rho$.
\end{enumerate}

Assume that for some $m-1\geq 1$ we constructed the train track
$\sigma_{m-1}$. Let $v$ be a switch in $\tau$ of height $m$ and
let $\rho:[0,m]\to \tau$ be the large one-way trainpath issuing
from $v$. Let $\gamma:[0,d]\to F^{-1}(\rho[0,1])$ be a cutting
arc for $\sigma$ and $v$.
We claim that $\gamma(d)$ is a switch of $\sigma$ 
contained in the interior of $F^{-1}(\rho[0,1])$. For
this assume to the contrary that this is not the case. Then
$\gamma(d)$ is contained in the
singular tie $F^{-1}(\rho(1))$. Since by assumption singular
ties in $A$ do not contain switches of $\sigma_{m-1}$, the point
$\gamma(d)$ is not a
switch of $\sigma_{m-1}$. Thus by possibly modifying $\gamma$ with
a small isotopy we may assume that 
$\gamma(d)\not\in \sigma_{m-1}$. 

Since the height of $v$ equals $m\geq 2$, 
$\rho[0,1]$ is a mixed branch in $\tau$ and
therefore the switch $\rho(1)$ divides the tie $F^{-1}(\rho(1))$
into two closed connected subarcs with intersection $\{\rho(1)\}$,
where one of the subarcs is the
intersection of $F^{-1}(\rho(1))$ with the boundary of the
rectangle $F^{-1}(\rho[0,1])$; we denote this arc by $c$.
Since $\sigma_{m-1}$ is contained in the interior of the foliated
neighborhood $A$ of $\tau$ and $\gamma(d)\in c$,
by possibly modifying $\gamma$ with another isotopy
we may assume that $\gamma(d)$ is
contained in the interior of $c$. 

Now $F$ maps $\sigma_{m-1}$ onto
$\tau$ and hence with respect to the orientation of $S$ and the
orientation of $\gamma$ there is a point in $\sigma_{m-1}\cap
F^{-1}(\rho(0))$  
to the right of $\rho(0)=\gamma(0)$, and there is a point in
$\sigma_{m-1}\cap F^{-1}(\rho(0))$ to the left of $\gamma(0)$. 
Since $\gamma$ is disjoint from $\sigma_{m-1}$ this implies
that there is a point $y\in c\cap \sigma_{m-1}$ to the
right of $\gamma(d)$, and there is a point $z\not=y\in c\cap 
\sigma_{m-1}$ to the left of $\gamma(d)$.  
However, by the
induction hypothesis, the intersection of $\sigma_{m-1}$ with $c$
consists of 
a single point and hence we arrive at a contradiction. 

As a consequence, the cutting arc $\gamma$ terminates at a switch $w$ of
$\sigma_{m-1}$ which is contained in the interior of
$F^{-1}(\rho[0,1])$. As in the proof of 
Lemma \ref{specialtrain} 
we can modify
$\sigma_{m-1}$ by an isotopy in such a way that for the resulting
train track, again denoted by $\sigma_{m-1}$, there are two small
half-branches which are incident on $w$ and are mapped by $F$ onto
$F(\gamma[0,d])$.

Let $\xi:[0,n]\to
\sigma_{m-1}$ be the large one-way trainpath issuing from the
switch $w$. Since $\sigma_{m-1}$ is transverse to the ties of $A$
and its intersection with the boundary component $c$ of the
rectangle $F^{-1}(\rho[0,1])$ consists of a unique point, the arc
$\xi[0,n]$ can not be contained in $F^{-1}(\rho[0,1])$. Then there
is a minimal number $t>0$ such that $\xi(t)\in c$. If $\ell\geq
0$ is the nonnegative integer such that $t\in [\ell,\ell+1)$, then
for $0<i\leq \ell$ the half-branch $\xi[i-1/2,i]$ is small.
Thus as in the proof of Lemma \ref{specialtrain}
we can modify
$\sigma_{m-1}$ by a sequence of shifts and an isotopy in such a
way that in the modified train track $\sigma^\prime$ the arc
$\xi[0,\ell+1]$ consists of a unique branch. Moreover, there is a standard
neighborhood $W$ of the trainpath $\rho$ in $\tau$ with the
property that $F^{-1}(\rho)\cap \sigma^\prime$ is a standard
neighborhood of a large trainpath on $\sigma^\prime$ which is
isotopic to $W$.

We can now repeat this construction for all large one-way
trainpaths on $\tau$ of length $m$. The resulting train track
$\sigma_m$ satisfies properties 1) and 2) above for $m$.

If $m_0$ denotes the maximal height of any large trainpath on
$\tau$, then after $m_0$ steps we obtain a train track
$\sigma_{m_0}$ which can be connected to $\sigma$ by a sequence of
shifts and an isotopy, is contained in $A$ and transverse to the
ties. Let $B\subset \tau$ be a closed subset which we
obtain from $\tau$ by removing
from each small branch $b$ of $\tau$ an open subset
whose closure is contained in the interior of $b$. 
By construction, the train
track $\sigma_{m_0}$ contains a closed subset $D$ which is
isotopic to $B$ and coincides with $F^{-1}(B)\cap \sigma_{m_0}$.
Since $\sigma_{m_0}$ and $\tau$ are both complete, they have the
same number of switches. Thus every switch of $\sigma_{m_0}$ is
contained in $D$ and therefore $\sigma_{m_0}$ is isotopic to
$\tau$. This finishes the proof of the lemma. 
\end{proof}

Now we are ready to show the main result
of this appendix which
relates splitting to carrying in
a quantitative way. For this define a 
\emph{splitting, shifting and collapsing sequence} to be a sequence
$\{\alpha_i\}_{0\leq i\leq m}$ of complete train tracks
such that for every $i$ the train track $\alpha(i+1)$
can be obtained from $\alpha(i)$ by a single split,
a single shift or a single collapse. Recall from the introduction
or from Subsection 2.1 the definition of a complete geodesic lamination.

\begin{proposition}\label{carrynear}
There is a number $\chi>0$ with the following property.
Let $\sigma\prec\tau\in {\cal V}({\cal T\cal T})$ and
let $\lambda$ be a complete geodesic lamination carried by
$\sigma$. Then
$\tau$ is splittable to a train track
$\tau^\prime\in {\cal V}({\cal T\cal T})$ which carries
$\lambda$ and which can be obtained from $\sigma$
by a splitting, shifting and collapsing 
sequence of length at most $\chi$
consisting of complete train tracks which carry $\lambda$.
\end{proposition}
\begin{proof} The number of complete train tracks which can
be obtained from a fixed train track by a sequence of shifts
is uniformly bounded. 
Hence it is enough to show that 
whenever $\sigma\prec\tau$ and $\lambda\in {\cal C\cal L}$
is carried by $\sigma$, then $\tau$ is splittable
to a train track $\sigma^\prime$ which carries $\lambda$ and
can be obtained from $\sigma$ by a uniformly bounded number
of modifications where each modification consists
in a number of shifts and a single split or a single collapse.

Thus let $\tau$ be a complete train track and let $A$ be a foliated
neighborhood of $\tau$ with collapsing map $F:A\to \tau$.
Let $\sigma\subset A$ be a complete
train track in general position with respect to $A$.
Let $\lambda$ be a complete geodesic
lamination which is carried by $\sigma$. We use the collapsing
map $F$ to define the quantities $\nu(b,\sigma)$ for 
a branch $b$ of $\tau$ as in the beginning of this section.
If $e$ is a large branch of 
$\tau$ with $\nu(e,\sigma)\geq 2$ and if 
$\sigma$ is carried by a split $\tilde\tau$ of $\tau$ at $e$ then 
for every standard neighborhood $V$ of $e$ there is 
a foliated neighborhood $\tilde A$ of $\tilde \tau$ which coincides
with $A$ outside of $F^{-1}(V)$ and such that $\sigma$ is 
in general position with respect to $\tilde A$.
Thus if $b\not=e$ is any branch different from $e$ with
$\nu(b,\sigma)=1$ then we may assume that $\nu(\tilde b,\sigma)=1$
for the branch $\tilde b$ of $\tilde \tau$ corresponding to 
$b$ under the natural identification of the branches of
$\tau$ with the branches of $\tilde\tau$. 

Let $m\geq 0$ be the number of branches $b$ of
$\tau$ with $\nu(b,\sigma)=1$.
Let $\tau_0$ be a train
track which can be obtained from $\tau$ by a splitting sequence,
which carries $\sigma$ and such that
no split of $\tau_0$ carries
$\sigma$. By the above observation, there is a 
foliated neighborhood $A_0$ of $\tau_0$ with collapsing
map $F_0:A_0\to\tau_0$ such that
$\sigma$ is in general position with respect to $A_0$ 
and that the number $m_0$ of branches $b$ of $\tau_0$ with
$\nu(b,\sigma)=1$ is not smaller than $m$.

By Lemma \ref{specialtrain},
$\sigma$ can be modified with a sequence of
shifts and an isotopy to a train track $\sigma_0$ which
is in general position with respect to $A_0$
and such that the following properties hold.
\begin{enumerate}
\item Every large branch $e$ of $\tau_0$
admits a standard neighborhood $V(e)$ in $\tau_0$ 
such that $F_0^{-1}(V(e))\cap \sigma_0$ is the disjoint union of 
a standard neighborhood $V_0$ of a special trainpath $\rho$ on
$\sigma_0$ and a (possibly empty) collection of simple
arcs not containing any switch of $\sigma_0$. 
The set $V_0$ is mapped by $F_0$ onto $V(e)$.
\item If $b\in \tau_0$ is any branch with $\nu(b,\sigma)=1$ 
then $\nu(b,\sigma_0)=1$.
\end{enumerate}

If $\nu(e,\sigma_0)=1$ for every large
branch $e$ of $\tau_0$ then by Lemma \ref{shiftequi},
$\sigma_0$ is shift equivalent to $\tau_0$ and we are done.
Otherwise we choose a large branch $e$ of 
$\tau_0$ with $\nu(e,\sigma_0)\geq 2$. We 
modify $\tau_0,\sigma_0$ to train tracks
$\tau_1,\sigma_1$ with the following properties.
\begin{enumerate}
\item[(a)] $\tau_1$ can be obtained from $\tau_0$ by at most one 
split at $e$.
\item[(b)] $\sigma_1$ carries $\lambda$ and
can be obtained from $\sigma_0$ by a sequence of collapses, shifts and splits
of uniformly bounded length.
\item[(c)] There is a foliated neighborhood $A_1$ of $\tau_1$ 
such that $\sigma_1$ is in general position with respect to
$A_1$, is transverse to the ties
and such that the number of branches $b$ of $\tau_1$ with
$\nu(b,\sigma_1)=1$ 
is not smaller than $m_0+1$.
\end{enumerate}

Namely,
using the above notations, let 
$2\ell-1\geq 1$ be the length of the special trainpath 
$\rho\subset F_0^{-1}(V(e))\cap \sigma_0$.
Note that $\ell\geq 1$ is bounded from above by a constant
only depending on the topological type of $S$.
If $\ell=1$ then $\rho$ consists of a single large
branch $e^\prime$ in $\sigma_0$, and
$F_0^{-1}(V(e))\cap \sigma_0$ 
is the disjoint union of a standard neighborhood of the large
branch $e^\prime$ and a collection of simple arcs.
There is a unique choice of 
a right or left split of $\sigma_0$ at $e^\prime$ such that
the split track $\sigma_1$ carries $\lambda$. If this split
is a right (or left) split, then $\sigma_1$ is carried
by the train track $\tau_1$ which can be obtained from
$\tau_1$ by a right (or left) split at $e$. 
If $\tilde e$ is the diagonal branch of the split in
$\tau_1$, then for a suitable choice of a foliated neighborhood
of $\tau_1$ and a suitable collapsing map we have 
$\nu(\tilde e,\sigma_1)=1$. Moreover, via
the natural identification of the branches of $\tau_0$
with the branches of $\tau_1$ we have
$\nu(\tilde b,\sigma_{1})\leq \nu(b,\sigma_{0})$ 
for every branch $b$ in $\tau_{0}$ and the corresponding
branch $\tilde b$ in $\tau_1$. In particular, the number
of branches $b\in \tau_1$ with $\nu(b,\sigma_1)=1$ 
is at least $m_0+1$ and the train tracks $\sigma_1\prec\tau_1$
satisfy the properties (a)-(c) above.

If the length $2\ell -1$ of the special trainpath $\rho$ on
$\sigma_0$ is at least 3 then $\rho[1,2]$ is a small branch which
can be collapsed as shown in 
Figure E. The train track $\sigma_0^\prime$ which is obtained from
$\sigma_0$ by this collapse is 
carried by $\tau_0$. In particular, this train track is transversely
recurrent and hence complete.
Moreover, for a suitable choice
of a carrying map
we have $\nu(b,\sigma_0^\prime)\leq \nu(b,\sigma_0)$ for
every branch $b$ of $\tau_0$. 

The train track $\sigma_0^\prime$ is not carried by any
split of $\tau_0$. Thus 
by Lemma \ref{specialtrain} and as shown in Figure E,
$\sigma_0^\prime$ can be shifted and isotoped 
to a train track $\eta$ which satisfies 
$\nu(b,\eta)=1$ if and only if $\nu(b,\sigma_0^\prime)=1$ 
for every branch $b$
of $\tau_0$ and such that there is a standard neighborhood $V$ of $e$
with the property that $F_0^{-1}(V)\cap \eta$ is the disjoint union
of a special trainpath of length $2\ell-3$ and a collection of 
simple arcs. Repeat this procedure with $\eta\prec\tau_0$. After
$\ell-1$ modifications of $\sigma_0$ consisting each
of a single collapse and a sequence of shifts and isotopies
we obtain a train track $\eta^\prime$ 
which is in general position with
respect to the foliated neighborhood $A_0$ 
and such that 
the preimage of a standard neighborhood $V(e)$ of 
$e$ intersects $\eta^\prime$ in  
a single large branch $e^\prime$ 
and a collection of simple arcs. Moreover, 
we have
$\nu(b,\eta^\prime)\leq \nu(b,\sigma_0)$ for every branch
$b$ of $\tau_0$. If $\nu(e,\eta^\prime)=1$ then 
the number of branches $b$ of $\tau_0$ with $\nu(b,\eta^\prime)=1$
is strictly bigger than the number of branches of 
$\tau_0$ with $\nu(b,\sigma_0)=1$ and
the pair $\sigma_1=\eta^\prime\prec\tau_1=\tau_0$ satisfies the above
properties.
Otherwise 
modify $\eta^\prime$ and $\tau_0$ with a single split at $e$ 
as above
to obtain a pair of train tracks 
which fulfill the requirements (a)-(c)
of the inductive construction.

Reapply the above procedure with the train tracks
$\sigma_1\prec\tau_1$. After at most $p$ steps where
$p$ is the number of branches of a complete train track on 
$S$ we obtain a pair of 
train tracks $\sigma^\prime\prec \tau^\prime$ with the following
properties.
\begin{enumerate}
\item $\tau^\prime$ carries $\lambda$ and 
is obtained from $\tau$ by a splitting sequence.
\item $\sigma^\prime$ is in general position with respect
to a foliated neighborhood $A^\prime$ of $\tau^\prime$.
\item $\nu(e,\sigma^\prime)=1$ for every
large branch $e$ of $\tau^\prime$. 
\item $\sigma^\prime$ carries $\lambda$ and 
can be obtained from $\sigma$ by a splitting, shifting
and collapsing sequence of uniformly bounded length
consisting of train tracks which carry $\lambda$.
\end{enumerate}
By Lemma \ref{shiftequi},
$\sigma^\prime$ and $\tau^\prime$ are
shift equivalent. This completes the proof of 
the proposition.
\end{proof}

\bigskip

{\bf Acknowledgement:} Part of this work was carried out during a
visit of the MSRI in Berkeley. I thank the Institute for its
hospitality. I am also grateful to Curt McMullen for useful comments.

\noindent
MATHEMATISCHES INSTITUT DER UNIVERSIT\"AT BONN\\ 
BERINGSTRASSE 1, 53115 BONN, GERMANY\\
e-mail: ursula@math.uni-bonn.de


\begin{thebibliography}{CEG87}


\bibitem[A96]{A96} S.~Adams, {\em Reduction of cocycles with
hyperbolic target}, Erg. Th. \& Dyn. Syst. 16 (1996), 1111--1145.

\bibitem[AR00]{AR00} C.~Anantharaman-Delaroche, J.~Renault,
{\sl Amenable groupoids}, Monographie 36 de l`Enseignement
Math., Gen\`eve 2000.

\bibitem[BCH94]{BCH94} P.~Baum, A.~Connes, N.~Higson,
{\em Classifying space for proper $G$-actions and
$K$-theory of group $C^*$-algebras}, Contemp. Math. 167
(1994), 241--291.


\bibitem[Bo97]{Bo97} F.~Bonahon, {\em Geodesic laminations on
surfaces}, Contemp. Math. 269 (1997), 1--37.

\bibitem[Bw08]{Bw03} B.~Bowditch, {\em Tight
geodesics in the curve complex}, Invent. Math 171 (2008),
281--300.


\bibitem[BH99]{BH99} M.~Bridson, A.~Haefliger, {\sl Metric
spaces of non-positive curvature}, Springer Grund\-leh\-ren 319,
Springer, Berlin 1999.




\bibitem[BM99]{BM99} M.~Burger, N.~Monod,
{\em Bounded cohomology of lattices
in higher rank Lie groups}, J.~Eur. Math. Soc. 1 (1999), 199--235.

\bibitem[BM02]{BM02} M.~Burger, N.~Monod,
{\em Continuous bounded cohomology and applications
to rigidity theory}, Geom. Funct. Anal. 12 (2002),
219--280.



\bibitem[B92]{B92} P.~Buser, {\sl Geometry and spectra
of compact Riemann surfaces}, Birkh\"auser, Boston 1992.

\bibitem[CEG87]{CEG87} R.~Canary, D.~Epstein, P.~Green,
{\em Notes on notes of Thurston}, in ``Analytical and geometric
aspects of hyperbolic space'', edited by D.~Epstein, London Math.
Soc. Lecture Notes 111, Cambridge University Press, Cambridge 1987.

\bibitem[CB88]{CB88} A.~Casson with S.~Bleiler, {\sl Automorphisms
of surfaces after Nielsen and Thurston}, Cambridge University
Press, Cambridge 1988.


\bibitem[FLM01]{FLM01} B.~Farb, A.~Lubotzky, Y.~Minski, 
{\em Rank one phenomena
for mapping class groups}, Duke Math. J. 106 (2001), 581-597.

\bibitem[FLP91]{FLP91} A.~Fathi, F.~Laudenbach, V.~Po\'enaru, {\sl Travaux de
Thurston sur les surfaces,} Ast\'erisque 1991.

\bibitem[H06]{H05} U.~Hamenst\"adt, {\em Train tracks
and the Gromov boundary of the complex of curves},
in "Spaces of Kleinian groups" 
(Y. Minsky, M. Sakuma, C. Series, eds.), 
London Math. Soc. Lec. Notes 329 (2006), 187-207.


\bibitem[H08]{H08} U.~Hamenst\"adt,
{\em Bounded cohomology and isometry groups of hyperbolic spaces},
to appear in J. Eur. Math. Soc. 


\bibitem[Har86]{Har86} J.~Harer, {\it The virtual cohomological
dimension of the mapping class group of an oriented surface},
Invent. Math. 84 (1986), 157-176.

\bibitem[Hat91]{Hat91} A.~Hatcher, {\it On triangulations of
surfaces}, Topology Appl. 40 (1991), 189-194.



\bibitem[H00]{H00} N.~Higson, {\em Biinvariant $K$-theory
and the Novikov conjecture},
Geom. Funct. Anal. 10 (2000), 563--581.

\bibitem[HR00]{HR00} N.~Higson, J.~Roe, 
{\em Amenable group actions and the Novikov
conjecture}, J. reine angew. Math. 519 (2000), 143--153.


\bibitem[HM79]{HM79} J.~Hubbard, H.~Masur,
{\em Quadratic differentials and foliations}, Acta Math. 142
(1979), 221--274.

\bibitem[I02]{I02} N.~V.~Ivanov, {\sl Mapping class groups},
Chapter 12 in Handbook of Geometric Topology (Editors
R.J.~Daverman and R.B.~Sher), Elsevier Science (2002), 523-633.

\bibitem[Ka03]{Ka03} V.~Kaimanovich, {\em Double
ergodicity of the Poisson boundary and applications
to bounded cohomology}, Geom. Funct. Anal. 13 (2003),
852--861.

\bibitem[Ka04]{Ka04} V.~Kaimanovich, {\em Boundary
amenability of hyperbolic spaces}, Contemp. Math 347 (2004),
83--114.

\bibitem[KM96]{KM96} V.~Kaimanovich, H.~Masur,
{\em The Poisson boundary of the mapping class group},
Invent. Math. 125 (1996), 221--264.

\bibitem[KMS86]{KMS86} S.~Kerckhoff, H.~Masur, J.~Smillie,
{\em Ergodicity of billiard flows and quadratic
differentials}, Ann. Math. 124 (1986), 293-311.

\bibitem[Ki05]{Ki05} Y.~Kida, {\em The mapping class
group from the viewpoint of measure equivalence theory},
arXiv:math.GR/0512230.

\bibitem[Kl99]{Kl99} E.~Klarreich, {\em The boundary
at infinity of the curve complex and the relative
Teichm\"uller space}, unpublished manuscript, Ann Arbor 1999. 

\bibitem[M82]{M82} H.~Masur, {\em Interval exchange transformations
and measured foliations}, Ann. Math. 115 (1982), 169-201.

\bibitem[MM99]{MM99} H.~Masur, Y.~Minsky, {\em Geometry of the
complex of curves I: Hyperbolicity}, Invent. Math. 138 (1999),
103-149.

\bibitem[MP89]{MP89} J.~McCarthy, A.~Papadopoulos, {\em Dynamics on
Thurston's sphere of projective measured foliations}, 
Comm. Math. Helv. 64 (1989), 133--166.



\bibitem[MMS04]{MMS04} I.~Mineyev, N.~Monod, Y.~Shalom,
{\em Ideal bicombings for hyperbolic groups and
applications}, Topology 43 (2004),1319--1344.

\bibitem[MV03]{MV03} G.~Mislin, A.~Valette, {\sl Proper
group actions and the Baum-Connes conjecture}, 
Birkh\"auser, Basel 2003.

\bibitem[MS04]{MS04} N.~Monod, Y.~Shalom, {\em Cocycle
superrigidity and bounded cohomology for negatively curved
spaces}, J.~Diff. Geom. 67 (2004), 395--456.



\bibitem[M03]{M03} L.~Mosher, {\em Train track expansions of measured
foliations}, unpublished manuscript, 2003.


\bibitem[PH92]{PH92} R.~Penner with J.~Harer, {\sl Combinatorics
of train tracks}, Ann. Math. Studies 125, Princeton University
Press, Princeton 1992.

\bibitem[SW82]{SW82} K.~Schmidt, P.~Walters, {\em Mildly mixing
actions of locally compact groups}, Proc. London Math. Soc. 45 (1982),
506--518.

\bibitem[T79]{T79} W.~Thurston, {\em Three-dimensional geometry
and topology}, unpublished manuscript, 1979.

\bibitem[V86]{V86} W.~Veech, {\em The Teichm\"uller geodesic
flow}, Ann. Math. 124 (1986), 441--530.

\bibitem[Y98]{Y98} G.~Yu, {\em The coarse Baum-Connes conjecture for
groups with finite asymptotic dimension}, Ann. Math. 147 (1998),
325--355. 

\bibitem[Y00]{Y00} G. Yu, {\em The coarse Baum-Connes
conjecture for spaces which admit a uniform embedding
into Hilbert space}, Invent. Math. 139 (2000), 201--240.

\bibitem[ZB04]{ZB04} X.~Zhu, F.~Bonahon, {\em The metric space of
geodesic laminations on a surface I}, Geom. Top. 8 (2004), 539--564.

\bibitem[Z84]{Z84} R.~Zimmer, {\sl Ergodic theory and
semisimple groups}, Birkh\"auser, Boston 1984.

\end{thebibliography}
\end{document}